\documentclass[draft,12pt]{article}

\usepackage{amsfonts,amsmath,amsthm,amscd,amssymb,latexsym,cite,verbatim,texdraw,floatflt,caption2,pb-diagram}

\usepackage[mathscr]{eucal}
\usepackage{xcolor}

\usepackage{mathrsfs}

\newcommand{\tsfrac}[2]{{\textstyle \frac{#1}{#2}}}

\newcommand{\banf}{{\sf Proof.}\ }
\newcommand{\bend}{\hspace*{0ex} \hfill \hbox{\vrule height
    1.5ex\vbox{\hrule width 1.4ex \vskip 1.4ex\hrule  width 1.4ex}\vrule
    height 1.5ex}}

\renewcommand{\abstract}{\textbf{Abstract.} }

\numberwithin{equation}{section}

\newtheorem{proposition}{Proposition}
\newtheorem*{proposition*}{Proposition}
\newtheorem{theorem}{Theorem}
\newtheorem*{theorem*}{Theorem}
\newtheorem{remark}{Remark}

\newtheorem*{definition*}{Definition}
\newtheorem{lemma}{Lemma}
\newtheorem*{lemma*}{Lemma}
\newtheorem{corollary}{Corollary}
\newtheorem*{corollary*}{Corollary}

\newcommand{\sep}{/\kern-2pt/ }

\usepackage{color}
\usepackage{graphicx}

\allowdisplaybreaks

\begin{document}

\title{Approximation of  generalized Wiener classes of functions
of several variables in different metrics}

\author{Andrii Shydlich \\
{\small Institute of Mathematics of the National Academy  of Sciences of Ukraine}}

\maketitle

\begin{abstract}
 The paper presents new and known results on estimates of important linear
 and nonlinear approximation characteristics of generalized Wiener classes of
 functions of several variables in different metrics.
\end{abstract}


\section{Introduction} 

 Let $d$ be a fixed positive integer $(d\in {\mathbb N})$, let $\mathbb R^d$
 and $\mathbb Z^d$  be the sets of ordered collections ${k}:=(k_1,\ldots,k_d)$ of $d$ real
 and  integer  numbers correspondingly. Let also $\mathbb T^d\!\!:=[0,2\pi)^d$
 denote $d$-dimensional torus, and let  $L_p:=L_p(\mathbb T^d)$, \mbox{$1\le p\le \infty,$} be the space
 of all Lebesgue-measurable on $\mathbb R^d$
$2\pi$-periodic in each variable functions $f$  with finite norm
 \[
 \|f\|_{_{\scriptstyle L_p}}:=\|f\|_{_{\scriptstyle L_p(\mathbb T^d)}}=\left\{\begin{matrix}\Big((2\pi)^{-d}\displaystyle{\int_{\mathbb T^d}}|f({ x})|^p{\rm d}{x}\Big)^{\frac 1p},
 \quad \hfill & 1\le p< \infty,\\
 \operatorname{ess\,sup}_{{x}\in \mathbb T^d}|f({x})|,\quad \hfill &  p=\infty.\end{matrix}\right.
 \]
  Set $({k,x}):=k_1x_1+k_2x_2+\ldots+k_dx_d$, $e_k(x):={\rm e}^{{\rm i}(k,x)}$ and for any $f\in L_1$,
  we denote the Fourier coefficients of $f$ by
  $\widehat f({k}):=(2\pi)^{-d}{\displaystyle\int}_{\!\!\!\!\mathbb T^d}f({x})\overline{e}_k(x){\rm d}{x}$, ${k}\in\mathbb Z^d,
$
where $\overline{z}$ is the complex conjugate of $z$.

Further, let  ${\mathcal S}^p:={\mathcal S}^p({\mathbb T}^d)$, $0<p\le\infty$, be  the space   of all
functions $f\in L_1$ with the finite $\ell_p$-(quasi-)norm
\begin{equation}\label{Eq:Intr01}
 \|f\|_{_{\scriptstyle {\mathcal S}^p}}:=\|\{\widehat{f} ({ k})\}_{k\in {\mathbb Z}^d}\|_{\ell_p({\mathbb Z}^d)}
 = \left\{\begin{matrix}
 \big(\sum_{k\in {\mathbb Z}^d}  |\widehat f({ k})|^p\big)^{\frac 1p},\quad 0<p<\infty,\\ \sup_{k\in {\mathbb Z}^d}|\widehat f({ k})|,\quad\quad\quad p=\infty.\end{matrix}\right.
\end{equation}
Approximation characteristics of the spaces ${\mathcal S}^p$ of one and several variables were
  actively studied in the papers of Stepanets, his students, and followers
  (see, for example \cite{Stepanets_2001, Stepanets_Serdyuk_UMZh2002, Vakarchuk_2004},
  \cite[Chap.~11]{Stepanets_2005},  \cite{Vakarchuk_Shchitov_2006, Stepanets_UMZh2006_1, Timan_Shavrova_2007, Savchuk_Shidlich_2014, Abdullayev_Ozkartepe_Savchuk_Shidlich_2019, Abdullayev_Serdyuk_Shidlich_2021}, etc.).
  It is also worth mentioning a series of papers
  \cite{Chaichenko_Savchuk_Shidlich_2020, Abdullayev_Chaichenko_Shidlich_2021, Chaichenko_Shidlich_Shulyk_2022, Serdyuk_Shidlich_2022, Chaichenko_Shidlich_2024, Volosivets_2025}
  devoted mainly to direct and inverse approximation theorems  using various linear methods in
  different spaces of  Besicovitch-Musielak-Orlicz-Stepanets  type  of periodic
  and almost periodic functions, which can be considered as certain generalizations
  of the spaces  ${\mathcal S}^p({\mathbb T}^1)$. An overview of results in this direction can be
  also found in  \cite{Serdyuk_Shidlich_2025}.

 The spaces ${\mathcal S}^p$ are also known as Wiener spaces or Wiener-type spaces,
 due to their connection with the classical Wiener algebra in the case where $p=1$,
 and are often denoted by ${\mathcal A}_p$.
 They can be seen as periodic versions of so-called Barron spaces \cite{Barron_1993, Voigtlaender_2022}
 and share similar properties regarding non-linear approximation.
 For $p=1$ the spaces ${\mathcal S}^p$ consist  of functions whose Fourier series are absolutely convergent.
 The questions of absolute convergence  and summability of general trigonometric
 series and Fourier series were, in particular,  studied
 by Wiener \cite{Wiener_1933}, Sz\'{a}sz \cite{Szasz_1942}, Stechkin
 \cite{Stechkin_1951, Stechkin_1953,Stechkin_1955_IZV, Stechkin_1955},
 Sunouchi \cite{Sunouchi_1951, Sunouchi_1952},  etc.
 The main results in this direction, the properties of such spaces, and their
 generalization were described in Kahane's monograph  \cite{Kahane_1970}. Let us also mention
 a significant series of results by M\'{o}ricz  \cite{Moricz_2006, Moricz_2008, Moricz_2008_1, Moricz_2010},
 etc., which establish a relationship  between the  absolute convergence of (multiple)
  Fourier series and the structural smoothness of functions, specifically by
  introducing  enlarged  Lipschitz and Zygmund classes defined by the mixed modulus of continuity.

 Let  $\Psi=\{\Psi_k\}_{k \in {\mathbb Z}^d}$ be a sequence of complex numbers, $\Psi_k\not=0$, and
\begin{equation}\label{Eq:Intr001}
{\mathcal F}_{q}^{\Psi}:={\mathcal F}_{q}^{\Psi}({\mathbb T}^d):=\Big\{f\in L_1\ :\ \    \|\{\widehat{f} ({k})/\Psi_k\}_{k\in {\mathbb Z}^d}\|_{\ell_q({\mathbb Z}^d)}\le 1\Big\},\quad 0<q\le \infty.
\end{equation}
Classes ${\mathcal F}_{q}^{\Psi}$ are also called weighted Wiener classes or generalized Wiener classes.

Next, let $\psi=\psi(t)$, $t\ge 1$, be a positive nonincreasing function. Consider the classes ${\mathcal F}_{q}^{\Psi}$ in the case when the sequence
 $\Psi$ satisfies the condition
  \begin{equation}\label{Eq:Intr02}
   |\Psi_0|=\psi(1)\qquad {\rm and}\qquad  |\Psi_k|=\psi(|k|_{r}),
 \end{equation}
 where  $|k|_{r}:=\|k\|_{\ell_r^d}$, $0<r\le \infty$, is the $\ell_r$-(quasi-)norm defined for  $k\in {\mathbb Z}^d$ similarly to \eqref{Eq:Intr01}.  In this case, we denote ${\mathcal F}_{q}^{\Psi}={\mathcal F}_{q}^{\Psi}({\mathbb T}^d)={\mathcal F}_{q,r}^{\psi}({\mathbb T}^d)=:{\mathcal F}_{q,r}^{\psi}$.

 The main goal of the paper is to find the exact order estimates for
  linear and nonlinear approximations of classes ${\mathcal F}_{q,r}^{\psi}$
 in various metrics depending on the rate of decay to zero of functions $\psi$.
 Similar estimates play important role for estimates of the corresponding approximations
 of  functional classes  such as Sobolev classes, Besov classes, etc.
 (see, for example, \cite{DeVore_Temlyakov_1995, Temlyakov_1998, Li_2010, Nguyens_2022}),
 and can also serve as upper bounds for (nonlinear) sampling errors in $L_p$
 (see \cite{Jahn_Ullrich_Voigtlaender_2023, Moeller_Stasyuk_Ullrich_2024, Moeller_Stasyuk_Ullrich_2025}, etc.).

  If  $\psi(t)=t^{-s}$, $s\in {\mathbb N}$ and $r=\infty$, then
  ${\mathcal F}_{q,\infty}^{\psi}=:{\mathcal F}_{q,\infty}^{s}$ is a set of functions
  whose $s$th partial derivatives have absolutely convergent Fourier series.
  It is also called periodic (isotropic) Wiener class. If $q=2$,  ${\mathcal F}_{q,\infty}^{s}$
  is equivalent  (modulo constants) to the unit ball of the Sobolev class $W^s_2$.
  Approximative characteristics of the classes ${\mathcal F}_{q,r}^{\psi}$ for different
  $r\in (0,\infty]$ and for the various functions $\psi$ were investigated by many authors
  (see, e.g., \cite{DeVore_Temlyakov_1995, Temlyakov_1998, Li_2010, Shidlich_2011, Shidlich_2013,Shidlich_2014, Romanyuk_2015, Shidlich_2016, Kolomoitsev_Lomako_Tikhonov_2023}).
  In particular,  DeVore and Temlyakov \cite{DeVore_Temlyakov_1995} found the exact order estimates
  for the best $m$-term trigonometric approximations of the classes ${\mathcal F}_{q,\infty}^{s}$, $s>0$,
  in the spaces $L_p$.  Temlyakov \cite{Temlyakov_1998} obtained the exact order estimates for approximations
  of these classes by $m$-term greedy polynomials  in  $L_p$.  In the case where $\psi(t)$ is a
  positive function that decreases to zero no faster than some power function,
  the best $m$-term one-sided trigonometric approximations and
  approximations by  $m$-term one-sided Greedy-liked polynomials of the classes
  ${\mathcal F}_{q,\infty}^{\psi}$ were studied in \cite{Li_2010}.  Linear approximation in
  Wiener type spaces was studied in  \cite{Kolomoitsev_Lomako_Tikhonov_2023}.

 If the sequence  $\Psi$ satisfies the conditions similar to \eqref{Eq:Intr02},
 where $\psi$ is a power function: $\psi(t)=t^{-s}$, and instead of the functional
 $|k|_{r}$ we consider the functional  $|k|_{mix}:=\prod_{i=1}^d(1+|k_i|)^r$, then the
 classes ${\mathcal F}_{q}^{\Psi}$ are denoted by ${\rm S}^{\psi}_q {\mathcal A}$ and called
 multivariate weighted Wiener classes with mixed smoothness. Approximative characteristics of
 such classes were studied in  \cite{Nguyens_2022, Jahn_Ullrich_Voigtlaender_2023, Moeller_Stasyuk_Ullrich_2024, Moeller_Stasyuk_Ullrich_2025}, etc.

 Note that nonlinear approximation  of the generalized Wiener classes was studied
 in \cite{Stepanets_2001},  \cite[Chap.~11]{Stepanets_2005},
 \cite{Nguyens_2022, DeVore_Petrova_Wojtaszczyk_2024}, etc. In particular,
 Stepanets \cite{Stepanets_2001},  \cite[Chap.~11]{Stepanets_2005}  and
 V.\,K.~Nguyen and V.\,D.~Nguyen \cite{Nguyens_2022} found the exact values of the
 best $m$-term trigonometric approximations of the classes ${\mathcal F}_{q}^{\Psi}$
 in the  spaces ${\mathcal S}^p$ for all  $0<p,q\le \infty$.


\section{Approximative characteristics}

 Let  ${\mathscr X}$ be one of the spaces $L_p$, $1\le p\le \infty$, or ${\mathcal S}^p$, $0<p\le\infty$, $m\in {\mathbb N}$, let $\gamma_m$  be a collection   of $m$ different vectors of ${\mathbb Z}^d$,
 and let $f$ be any  function from ${\mathscr X}$. The quantity
   \begin{equation}\label{Eq:ApCh01}
          E_{\gamma_m}(f)_{_{\scriptstyle {\mathscr X}}} =
          \inf_{c_k\in {\mathbb C}} \Big\|\,f - \sum \limits_{k \in {\gamma_m}}c_k\,e_k  \Big\|_{_{\scriptstyle {\mathscr X}}}
    \end{equation}
 is called the best approximation of the function $f$   by $m$-term polynomials corresponding to the collection $\gamma_m$ in the space ${\mathscr X}$.

 Next, let  $S_{\gamma_m}(f)=\sum_{k \in \gamma_m}\widehat f({ k}) e_k$ be the Fourier sum
 corresponding to the collection $\gamma_m$, and
    \begin{equation}\label{Eq:ApCh02}
  {\mathscr E}_{\gamma_m}(f)_{_{\scriptstyle {\mathscr X}}}=\|\,f - S_{\gamma_m}(f)
  \|_{_{\scriptstyle {\mathscr X}}}
  \end{equation}
 be the approximation of the function $f$ by the Fourier sum $S_{\gamma_m}(f)$ in  ${\mathscr X}$.

 If ${\mathfrak N}$ is a subset of the space ${\mathscr X}$, then
 $E_{\gamma_m}(\mathfrak N)_{_{\scriptstyle {\mathscr X}}}$ and ${\mathscr E}_{\gamma_m}(\mathfrak N)_{_{\scriptstyle {\mathscr X}}}$ denote the exact upper bounds
 of the quantities  \eqref{Eq:ApCh01} and \eqref{Eq:ApCh02} over the set ${\mathfrak N}$, i.e.,
    \begin{equation}\label{Eq:ApCh03}
          E_{\gamma_m}(\mathfrak N)_{_{\scriptstyle {\mathscr X}}} = \sup\limits_{f\in\mathfrak  N} E_{\gamma_m}(f)_{p}\qquad \mbox{\rm and}\qquad
          {\mathscr E}_{\gamma_m}(\mathfrak N)_{_{\scriptstyle {\mathscr X}}} = \sup\limits_{f\in\,\mathfrak N} {\mathscr E}_{\gamma_m}(f)_{_{\scriptstyle {\mathscr X}}}.
     \end{equation}
Denote by $\Gamma_m$ a set of all collections   of $m$ different vectors of ${\mathbb Z}^d$. The quantities
\begin{equation}\label{Eq:ApCh04}
   {\mathscr D}_m({\mathfrak N})_{_{\scriptstyle {\mathscr X}}} = \inf\limits_{\gamma_m\in \Gamma_m}
  E_{\gamma_m}({\mathfrak N})_{_{\scriptstyle {\mathscr X}}}\qquad \mbox{\rm  and}\qquad
   {\mathscr D}_m^\perp({\mathfrak N})_{_{\scriptstyle {\mathscr X}}} = \inf\limits_{\gamma_m\in \Gamma_m}
  {\mathscr E}_{\gamma_m}({\mathfrak N})_{_{\scriptstyle {\mathscr X}}}
 \end{equation}
 are called the basis width and projection width (or Fourier-width) of order $m$ of the set ${\mathfrak N}$
 in   ${\mathscr X}$.

Further, for  $f\in {\mathscr X}$, let $\{{k}_l\}_{l=1}^\infty=\{{k}_l(f)\}_{l=1}^\infty$ denote a rearrangement of vectors of ${\mathbb Z}^d$ such that
\begin{equation}\label{Eq:ApCh05}
|\widehat{f}({k}_1)|\ge |\widehat{f}({k}_2)|\ge \ldots.
\end{equation}
In general case, this rearrangement is not unique. In such  case, we take any rearrangement satisfying \eqref{Eq:ApCh05}.

We define $\Sigma_m$ to be  the class of all complex trigonometric polynomials of the form
 $
 T\,{=}\sum_{{ k}\in \gamma_m} c_k e_k,
 $
where  $\gamma_m$ is a  collection from the set $\Gamma_m$.

 In addition to \eqref{Eq:ApCh01}--\eqref{Eq:ApCh04},
 consider the quantities
\begin{equation}\label{Eq:ApCh07}
\|f-G_m(f)\|_{_{\scriptstyle {\mathscr X}}}=\Big\|f(\cdot)-\sum\limits_{l=1}^m \widehat{f}({k}_l)e_{{k}_l}\Big\|_{_{\scriptstyle {\mathscr X}}},
\end{equation}
 \begin{equation}\label{Eq:ApCh08}
 \sigma_m^\perp(f)_{_{\scriptstyle {\mathscr X}}}=\inf\limits_{\gamma_m\in \Gamma_m}\Big\|f-\sum\limits_{{ k}\in \gamma_m}\widehat{f}({ k})e_k\Big\|_{_{\scriptstyle {\mathscr X}}}=\inf\limits_{\gamma_m\in \Gamma_m} {\mathscr E}_{\gamma_m}(f)_{_{\scriptstyle {\mathscr X}}},
 \end{equation}
and
\begin{equation}\label{Eq:ApCh09}
\sigma_m(f)_{_{\scriptstyle {\mathscr X}}}=\inf\limits_{T\in \Sigma_m}\|f-T\|_{_{\scriptstyle {\mathscr X}}}=\inf\limits_{\gamma_m\in \Gamma_m} {E}_{\gamma_m}(f)_{_{\scriptstyle {\mathscr X}}}.
\end{equation}
The quantities \eqref{Eq:ApCh09} and  \eqref{Eq:ApCh08} are respectively called  the best $m$-term trigonometric and the best $m$-term orthogonal trigonometric approximations of the function $f$ in the space ${\mathscr X}$. The quantity \eqref{Eq:ApCh07} is called the approximation  of the function $f$ by $m$-term greedy  polynomials in the space ${\mathscr X}$.

For a set  ${\mathfrak N}\subset {\mathscr X}$, we put
 \begin{equation}\label{Eq:ApCh10}
  \sigma_m^\perp({\mathfrak N})_{_{\scriptstyle {\mathscr X}}}=\sup\limits_{f\in {\mathfrak N}}\sigma_m^\perp(f)_{_{\scriptstyle {\mathscr X}}} \qquad \mbox{\rm and}\qquad
\sigma_m({\mathfrak N})_{_{\scriptstyle {\mathscr X}}}=\sup\limits_{f\in {\mathfrak N}}\sigma_m(f)_{_{\scriptstyle {\mathscr X}}}.
\end{equation}
In general case, the quantities \eqref{Eq:ApCh07} depend on the choice of the rearrangement satisfying \eqref{Eq:ApCh05}. So, for the unique definition, we put
\begin{equation}\label{Eq:ApCh11}
 G_m({\mathfrak N})_{_{\scriptstyle {\mathscr X}}}=
 \sup\limits_{f\in {\mathfrak N}}\inf\limits_{\{k_l(f)\}_{l=1}^\infty}
 \Big\|f(\cdot)-\sum\limits_{l=1}^m \widehat{f}(k_l(f))e_{k_l(f)}\Big\|_{_{\scriptstyle {\mathscr X}}}.
\end{equation}
In \eqref{Eq:ApCh11}, for any function $f\in {\mathfrak N}$, we consider the infimum on all rearrangements, satisfying \eqref{Eq:ApCh05}, but it should be noted that results, formulated in this paper, are also true for any other rearrangements, satisfying \eqref{Eq:ApCh05}.


 Research of the quantities of the form \eqref{Eq:ApCh07}--\eqref{Eq:ApCh09} goes back to
 the paper of Stechkin \cite{Stechkin_1955}. Order estimates of these quantities on different
 classes of functions of one and several variables were obtained by many authors.
 In particular,  the bibliography of papers with the similar results can be found in
 \cite{DeVore_1998, Temlyakov1993, Temlyakov_2011, Romanyuk_2012, Dung_Temlyakov_Ullrich_2018}.


 It follows from \eqref{Eq:ApCh01}--\eqref{Eq:ApCh11}  that
 \begin{equation}\label{Eq:ApCh12}
   \sigma_m(f)_{_{\scriptstyle L_p}}\le \sigma_m^\perp (f)_{_{\scriptstyle L_p}}\le \|f-G_m(f)\|_{_{\scriptstyle L_p}}\qquad \forall f\in L_p,
  \end{equation}
\begin{equation}\label{Eq:ApCh13}
\sigma_m(f)_{_{\scriptstyle {\mathcal S}^p}}=\sigma_m^\perp (f)_{_{\scriptstyle {\mathcal S}^p}}=\|f-G_m(f)\|_{_{\scriptstyle {\mathcal S}^p}} \qquad \forall f\in {\mathcal S}^p,
\end{equation}
 and   for any $\gamma_m\subset {\mathbb Z}^d$
 and ${\mathfrak N}\subset {\mathscr X}$
 \begin{equation}\label{Eq:ApCh15}
  \sigma_m({\mathfrak N})_{_{\scriptstyle {\mathscr X}}}\le
   {\mathscr D}_m({\mathfrak N})_{_{\scriptstyle {\mathscr X}}}
   \le {E}_{\gamma_m} ({\mathfrak N})_{_{\scriptstyle {\mathscr X}}}, \quad
    \sigma_m^\perp({\mathfrak N})_{_{\scriptstyle {\mathscr X}}}\le {\mathscr D}_m^\perp({\mathfrak N})_{_{\scriptstyle {\mathscr X}}}\le{\mathscr E}_{\gamma_m} ({\mathfrak N})_{_{\scriptstyle {\mathscr X}}}.
 \end{equation}



 \section{Approximative characteristics of the classes ${\mathcal F}_{q}^{\Psi}$ in the spaces ${\mathcal S}^p$}\label{Sec:AC_Sp}

\subsection{Exact values of best approximations and basis widths} Let  $\Psi=\{\Psi_k\}_{k \in {\mathbb Z}^d}$ be a sequence of complex numbers, $\Psi_k\not=0$,  such that  there exists a non-increasing rearrangement $\bar\Psi= \{\bar \Psi_j\}_{j=1}^\infty$ of the number system   $\{|\Psi_k|\}_{k \in {\mathbb Z}^d}$.  It is clear that in this case the system $\{|\Psi_k|\}_{k \in {\mathbb Z}^d}$ is bounded, i.e.
   \begin{equation}\label{Eq:ExV_BaW01}
   |\Psi_k|\le K\qquad \forall k \in {\mathbb Z}^d.
 \end{equation}
 Hereinafter,  $K$, $c$, $K_0,\ldots$ are positive constants that do not depend on the corresponding variable ($k$, $t$, etc.).

 A sufficient condition guaranteeing the existence of the rearrangement $\bar\Psi= \{\bar \Psi_j\}_{j=1}^\infty$ is the condition
 \begin{equation}\label{Eq:ExV_BaW02}
   \lim\limits_{|k|\to\infty} |\Psi_k|=0,
 \end{equation}
 however such rearrangement also exists, for example,  when  $\{|\Psi_k|\}_{k \in {\mathbb Z}^d}$ is a constant.

  In the case where ${\mathfrak N}={\mathcal F}_{q}^{\Psi}$ and ${\mathscr X}={\mathcal S}^p$,
  the values of the characteristics \eqref{Eq:ApCh03} and \eqref{Eq:ApCh04} were found by
  Stepanets \cite[Ch. XI]{Stepanets_2005}, \cite{Stepanets_UMZh2006_1} for all  $0<p,q<\infty$.  In order to formulate this result,
  for any collection $\gamma_m$ of $m$ different vectors of ${\mathbb Z}^d$,
  by $\Psi_{\gamma_m}=\{\Psi_{\gamma_m}(k)\}_{k\in {\mathbb Z}^d}$  denote  a system of numbers such that
 \begin{equation}\label{Eq:ExV_BaW03}
  \Psi_{\gamma_m}(k)=\left\{\begin{matrix}0,\quad\hfill & k\in\gamma_m,\\
  \Psi_k,\quad\hfill & k\overline{\in} \gamma_m,\end{matrix}\right.
  \end{equation}
 and  by $\bar{\Psi}_{\gamma_m}=\{\bar{\Psi}_{\gamma_m}(j)\}_{j=1}^\infty$ denote  a non-increasing rearrangement
  of the system $\{|\Psi_{\gamma_m}(k)|\}_{k\in {\mathbb Z}^d}$.

 {\bf Theorem A (\!\!\cite[Ch. XI]{Stepanets_2005}, \cite{Stepanets_UMZh2006_1}).}
 {\it Let $0<p,q<\infty$, $m\in {\mathbb N}$ and  $\Psi=\{\Psi_k\}_{k \in {\mathbb Z}^d}$ be a sequence of complex numbers, $\Psi_k\not=0$,
 such that  there exists a non-increasing rearrangement $\bar\Psi= \{\bar \Psi_j\}_{j=1}^\infty$ of the number system   $\{|\Psi_k|\}_{k \in {\mathbb Z}^d}$.

 {\rm (i)} In the case $0<p\le q<\infty$,
 \begin{equation}\label{Eq:ExV_BaW04}
   E_{\gamma_m}({\mathcal F}_{q}^{\Psi})_{_{\scriptstyle {\mathcal S}^p}}=
  {\mathscr E}_{\gamma_m}({\mathcal F}_{q}^{\Psi})_{_{\scriptstyle {\mathcal S}^p}} =\bar{\Psi}_{\gamma_m}(1) \qquad \forall \gamma_m\in \Gamma_m,
     \end{equation}
 and
      \begin{equation}\label{Eq:ExV_BaW05}
 {\mathscr D}_m({\mathcal F}_{q}^{\Psi})_{_{\scriptstyle {\mathcal S}^p }}
 ={\mathscr D}_m^\perp({\mathcal F}_{q}^{\Psi})_{_{\scriptstyle {\mathcal S}^p }}=\bar{\Psi}_{m+1}.
     \end{equation}

 {\rm (ii)} In the case when $0<q<p<\infty$ and $\sum_{k\in {\mathbb Z}^d} |\Psi_k|^\frac{pq}{q-p}<\infty$,
 \begin{equation}\label{Eq:ExV_BaW06}
   E_{\gamma_m}({\mathcal F}_{q}^{\Psi})_{_{\scriptstyle {\mathcal S}^p}}=
  {\mathscr E}_{\gamma_m}({\mathcal F}_{q}^{\Psi})_{_{\scriptstyle {\mathcal S}^p}} =\Big( \sum_{k=1}^{\infty}  \bar{\Psi}_{\gamma_m}^{\frac{p\,q}{q-p}}(k) \Big)^{\frac{q-p}{p\,q}} \qquad \forall \gamma_m\in \Gamma_m,
     \end{equation}
  and
      \begin{equation}\label{Eq:ExV_BaW07}
 {\mathscr D}_m({\mathcal F}_{q}^{\Psi})_{_{\scriptstyle {\mathcal S}^p}}={\mathscr D}_m^\perp({\mathcal F}_{q}^{\Psi})_{_{\scriptstyle {\mathcal S}^p}}
    =\Big( \sum_{k=m+1}^{\infty} \bar{\Psi}_k^{\frac{p\,q}{q-p}}
  \Big)^{\frac{q-p}{p\,q}}.
     \end{equation}
  Moreover, for any collection $\gamma_m^*=\{k_1,\ldots, k_m\}$   from the set $\Gamma_m$ such that
 \begin{equation}\label{Eq:ExV_BaW08}
    \gamma_m^*=\{k_j\in {\mathbb Z}^d:\quad |\Psi_{k_j}|=\bar{\Psi}_j,\qquad j=1,2,\ldots,m\},
 \end{equation}
 the following relation holds  in both cases:
   \begin{equation}\label{Eq:ExV_BaW09}
 {\mathscr D}_m({\mathcal F}_{q}^{\Psi})_{_{\scriptstyle {\mathcal S}^p}}={\mathscr D}_m^\perp({\mathcal F}_{q}^{\Psi})_{_{\scriptstyle {\mathcal S}^p}}=E_{\gamma_m^*}({\mathcal F}_{q}^{\Psi})_{_{\scriptstyle {\mathcal S}^p}}=
  {\mathscr E}_{\gamma_m^*}({\mathcal F}_{q}^{\Psi})_{_{\scriptstyle {\mathcal S}^p}}.
   \end{equation}}

Next, we formulate a theorem that complements this statement for cases of infinite values of parameters $p$ and $q$.

 \begin{theorem}\label{Th:ApChFS01}
Let $m\in {\mathbb N}$ and  $\Psi=\{\Psi_k\}_{k \in {\mathbb Z}^d}$ be a sequence of complex numbers, $\Psi_k\not=0$,
 such that  there exists a non-increasing rearrangement $\bar\Psi= \{\bar \Psi_j\}_{j=1}^\infty$ of the system   $\{|\Psi_k|\}_{k \in {\mathbb Z}^d}$.

{\rm (i)} If $0<q<p=\infty$ or $p=q=\infty$,  then relations \eqref{Eq:ExV_BaW04} and \eqref{Eq:ExV_BaW05} hold.

 {\rm (ii)} If $0<p<q=\infty$ and the series $\sum_{k\in {\mathbb Z}^d} |\Psi_k|^p$ converges,  then
 \begin{equation}\label{Eq:ExV_BaW10}
   E_{\gamma_m}({\mathcal F}_{q}^{\Psi})_{_{\scriptstyle {\mathcal S}^p}}=
  {\mathscr E}_{\gamma_m}({\mathcal F}_{q}^{\Psi})_{_{\scriptstyle {\mathcal S}^p}} =\Big( \sum_{k=1}^{\infty}  \bar{\Psi}_{\gamma_m}^{p}(k) \Big)^{\frac{1}{p}} \qquad \forall \gamma_m\in \Gamma_m,
     \end{equation}
  and
      \begin{equation}\label{Eq:ExV_BaW11}
 {\mathscr D}_m({\mathcal F}_{q}^{\Psi})_{_{\scriptstyle {\mathcal S}^p}}={\mathscr D}_m^\perp({\mathcal F}_{q}^{\Psi})_{_{\scriptstyle {\mathcal S}^p}}
    =\Big( \sum_{k=m+1}^{\infty} \bar{\Psi}_k^{p}\Big)^{\frac{1}{p}}.
     \end{equation}
Moreover, for any collection $\gamma_m^*\in\Gamma_m$ satisfying \eqref{Eq:ExV_BaW08}, relation  \eqref{Eq:ExV_BaW09}  holds. \end{theorem}

The proof of this statement, as well as the proofs of all other statements in this article, will be given in Section \ref{Proofs}.




 \subsection{Exact values of best $n$-term approximation}\label{EV_Sp}
 {\bf \ref{EV_Sp}.1.} Exact values of   the quantities
 $\sigma_m({\mathcal F}_{q}^{\Psi})_{_{\scriptstyle {\mathcal S}^{p}}}$  were obtained in \cite{Stepanets_2001}, \cite[Ch. XI]{Stepanets_2005} ((i), (ii))  and \cite{Nguyens_2022}
  ((iii)-(v)). They follow from the following statement.

 {\bf Theorem B (\!\!\cite{Stepanets_2001}, \cite[Ch. XI]{Stepanets_2005}, \cite{Nguyens_2022}).}
 {\it Let $0<p,q\le \infty$, $m\in {\mathbb N}$ and  $\Psi=\{\Psi_k\}_{k \in {\mathbb Z}^d}$ be a sequence of complex numbers
 such that  there exists a non-increasing rearrangement $\bar\Psi= \{\bar \Psi_j\}_{j=1}^\infty$ of the number system
  $\{|\Psi_k|\}_{k \in {\mathbb Z}^d}$.

 {\rm (i)} If $0<q\le p<\infty$, then
  \begin{equation}\label{Eq:ExV_BnA01}
  \sigma_m({\mathcal F}_{q}^{\Psi})_{_{\scriptstyle {\mathcal S}^p}}=\sup\limits
   _{l>m}\frac{(l-m)^\frac 1p}{(\sum _{j=1}^l\bar \Psi^{-q}_j)^{\frac 1{q}}}.
\end{equation}

 {\rm (ii)} If $0<p<q<\infty$ and the series $\sum_{k\in {\mathbb Z}^d} |\Psi_k|^\frac{pq}{q-p}$ converges,  then
 \begin{equation}\label{Eq:ExV_BnA02}
  \sigma_m({\mathcal F}_{q}^{\Psi})_{_{\scriptstyle {\mathcal S}^p}}= \bigg(  (l_m-m)^{\frac{q}{q-p}}\Big(\sum\limits_{j=1}^{l_m} \bar{\Psi}^{-q}_j\Big)^{\frac{p}{p-q}}
  +    \sum\limits_{j=l_m+1}^{\infty} \bar{\Psi}^\frac{pq}{q-p}_j\bigg)^{\frac{q-p}{pq}}  ,
  \end{equation}
where   the number $l_m$ is defined by
 \begin{equation}\label{Eq:ExV_BnA03}
   \bar{\Psi}^{-q}_{l_m} \leq
   \tsfrac{1}{l_m-m}\,\sum\limits_{j=1}^{l_m} \bar{\Psi}^{-q}_j < \bar{\Psi}^
   {-q}_{l_m+1}.
 \end{equation}

{\rm (iii)} If $0<q<p=\infty$,  then
 \begin{equation}\label{Eq:ExV_BnA04}
\sigma_m({\mathcal F}_{q}^{\Psi})_{_{\scriptstyle {\mathcal S}^p}}=
 \Big(\sum\limits_{j=1}^{m+1} \bar{\Psi}^{-q}_j\Big)^{-\frac 1q}.
  \end{equation}

 {\rm (iv)} If $0<p<q=\infty$ and the series $\sum_{k\in {\mathbb Z}^d} |\Psi_k|^p$ converges,  then
 \begin{equation}\label{Eq:ExV_BnA05}
\sigma_m({\mathcal F}_{q}^{\Psi})_{_{\scriptstyle {\mathcal S}^p}}=
 \Big(\sum\limits_{j=m+1}^{\infty} \bar{\Psi}^{p}_j\Big)^{\frac 1p}.
  \end{equation}

   {\rm (v)} If $p=q=\infty$,  then
 \begin{equation}\label{Eq:ExV_BnA06}
\sigma_m({\mathcal F}_{q}^{\Psi})_{_{\scriptstyle {\mathcal S}^p}}=
  \bar{\Psi}_{m+1}.
  \end{equation}
 }

 From the formulation of Theorems A, B, and 1 we can see that the values of
 ${\mathscr D}_m({\mathcal F}_{q}^{\Psi})_{_{\scriptstyle {\mathcal S}^p}}$,
 ${\mathscr D}_m^\perp({\mathcal F}_{q}^{\Psi})_{_{\scriptstyle {\mathcal S}^p}}$ and
 $\sigma_m({\mathcal F}_{q}^{\Psi})_{_{\scriptstyle {\mathcal S}^p}}$ significantly depend
 on the behavior of the sequence $\Psi$. Therefore, in order to obtain their estimates (as $m\to\infty$)
 for a specific sequence $\Psi$, it is necessary to additionally investigate the behavior of the functionals
 given in the right-hand sides of the corresponding relations.

 In this paper, by studying such functionals and using Theorems A, B, and 1, we find the asymptotic behavior of the above quantities for classes ${\mathcal F}_{q}^{\Psi}$ in the case when the sequence $\Psi$  satisfies  condition \eqref{Eq:Intr02},  where $\psi=\psi(t)$, $t\ge 1$, is a positive non-increasing function, i.e., for classes ${\mathcal F}_{q,r}^{\psi}$.

  Note that in the terms of similar functionals, solutions of other problems of
  approximation theory are formulated (see, e.g.,  \cite{Pinkus_1985} (Ch.~6),
  \cite{Stepanets_Shydlich_2003}, \cite{Fang_Qian_2006}, \cite{Fang_Qian_2007},
  \cite{Stepanets_Shidlich_2010}). Therefore, the study of such functionals  is interesting  in itself.

 {\bf \ref{EV_Sp}.2.} Let us provide a few auxiliary facts and estimates regarding the structure of sequences of the form
 \eqref{Eq:Intr02}.
 It is clear that for any  the sequence $\Psi$ satisfying condition \eqref{Eq:Intr02} with non-increasing positive function $\psi$, there exists a non-increasing rearrangement $\bar\Psi= \{\bar \Psi_j\}_{j=1}^\infty$ of the number system
  $\{|\Psi_k|\}_{k \in {\mathbb Z}^d}$ and it has the   stepwise form:
  \begin{equation}\label{Eq:ExV_BnA07}
   \bar\Psi_j=\psi(s),\quad j\in (V_{s-1},V_{s}],\ \ s=0,1,2,\ldots,
 \end{equation}
where
 \begin{equation}\label{Eq:ExV_BnA08}
   V_s=V_s(d):=\#\{k\in {\mathbb Z}^d:\|k\|_{\ell_r^d}\le  s\},\qquad  s=0,1,\ldots,\qquad  V_{-1}:=0.
 \end{equation}
Therefore, to study functionals that depend on such sequences $\Psi$,
it is necessary to have convenient estimates for the numbers $V_s$.

 Let  $B_r^d=\{x\in {\mathbb R}^d:\ \|x\|_{\ell_r^d}\le 1\}$  be the unit $\ell_r$-ball
 of the space ${\mathbb R}^d$.  It is known (see, e.g., \cite{Wang_2005}) that the volume of
 the ball $B_r^d$ is calculated by the formula
   \begin{equation}\label{Eq:ExV_BnA09}
 M_{r,d}:=\operatorname{vol}(B_r^d)=\tsfrac{\big(2\Gamma(1+\frac 1r)\big)^d}{\Gamma(1+\frac dr)},\qquad r\in (0,\infty],
 \end{equation}
and in particular, $M_{1,d}=\tsfrac {2^d}{d!}$, $M_{2,d}=\tsfrac {\pi^{d/2}}{\Gamma(\frac d2+1)}$ and $M_{\infty,d}=2^d$.

 For a real number $a$, we denote   $(a)_+=\max\{0,a\}$.  For positive sequences  $a(s)$ and $b(s)$
 (or functions $a(s)$ and $b(s)$, $s\ge 1$), the expression   ``$a(s)\asymp b(s)$''
 means that there are constants $0<K_1<K_2$ such that for any $s\in\mathbb{N}$ (or $s\ge 1$),
 $\ a(s)\le K_2b(s)$ (in this case, we write ``$a(s) \ll b(s)$'') and $a(s)\ge K_1b(s)$
 (in this case, we write  ``$a(s) \gg b(s)$'').

  The following lemma follows from known general results on  asymptotic behavior of
  the number of lattice points in scaled convex bodies (see, e.g., \cite[Ch. 14]{Gruber_2007}).

 \begin{lemma} Let $r\in(0,\infty]$ and $d\in\mathbb N$.
Then  for any $s=0,1,2,\ldots$,
 \begin{equation}\label{Eq:ExV_BnA10}
M_{r,d} \big((s-c_{r,d})_+\big)^d \le V_s \le M_{r,d} (s+c_{r,d})^d,
  \end{equation}
where $c_{r,d}:=\tsfrac {d^\frac 1r}{2}$ and
  \begin{equation}\label{Eq:ExV_BnA11}
   \nu_s:=V_s-V_{s-1}\asymp s^{d-1}.
  \end{equation}

\end{lemma}

\banf   Relation \eqref{Eq:ExV_BnA11} follows from  \eqref{Eq:ExV_BnA10}.
 Therefore, it suffices to prove \eqref{Eq:ExV_BnA10}.
For each $k\in\mathbb Z^d$, let $Q_k := k + [-\tfrac12,\tfrac12]^d$ be the unit cube centered at $k$. The family $\{Q_k\}_{k\in\mathbb Z^d}$ forms a tiling of
$\mathbb R^d$, and $\operatorname{vol}(Q_k)=1$ for all $k$. By definition, the values $V_s$ can be expressed as the volume of the union of these cubes:
\[
V_s = \operatorname{vol}\Big(\mathop{\bigcup}\limits_{\|k\|_{\ell_r^d}\le s} Q_k \Big).
\]
Note  the value $c_{r,d}$ is the maximum distance from the center to any point in the
unit cube $[- \frac 12,\frac 12]^d$,  i.e.,
$
c_{r,d}= \sup_{y\in[- \frac 12,\frac 12]^d} \|y\|_{\ell_r^d}.
$

\textit{Upper bound.} Suppose $\|k\|_{\ell_r^d}\le s$. For any point $x \in Q_k$, we have
$x-k \in [-\frac{1}{2}, \frac{1}{2}]^d$ and by  the triangle inequality
\[
\|x\|_{\ell_r^d} \le \|k\|_{\ell_r^d} + \|x-k\|_{\ell_r^d} \le s + c_{r,d}.
\]
This implies that the union of cubes is contained within the scaled ball of the radius $s+c_{r,d}$, i.e.,
\[
\bigcup_{\|k\|_{\ell_r^d}\le s} Q_k \subset (s+c_{r,d}) B_r^d.
\]
Taking the volume of both sides of this relation, we obtain
\[
V_s \le \operatorname{vol}((s+c_{r,d})B_r^d) = M_{r,d} (s+c_{r,d})^d.
\]

\textit{Lower bound.} It is sufficient to consider the case  when $s>c_{r,d}$.
Let $x \in (s-c_{r,d})B_r^d$. Since the cubes $Q_k$ tile $\mathbb{R}^d$,
 $x$ must belong to $Q_k$ for some $k \in \mathbb{Z}^d$. For this $k$, we have
\[
\|k\|_{\ell_r^d} \le \|x\|_{\ell_r^d} + \|k-x\|_{\ell_r^d} \le (s-c_{r,d}) + c_{r,d} = s.
\]
This shows that every point in the smaller ball is covered by a cube whose center is in the set counted by $V_s$, i.e.,
\[
(s-c_{r,d})B_r^d \subset \bigcup_{\|k\|_{\ell_r^d}\le s} Q_k.
\]
This similarly yields
\[
M_{r,d} (s-c_{r,d})^d \le V_s.
\]
\bend

\begin{remark} Denoting for a given positive integer $s$ by $n_s$ a number such that
 \begin{equation}\label{Eq:ExV_BnA12}
   V_{n_{s}-1}<s\le  V_{n_{s}},
  \end{equation}
we see that
    \begin{equation}\label{Eq:ExV_BnA13}
 \big((s/M_{r,d})^\frac 1d-c_{r,d}\big)_+\le n_s <
  (s/M_{r,d})^\frac 1d+c_{r,d}+1.
  \end{equation}

 \end{remark}

  \subsection{Order estimates of best $n$-term approximations  and basis widths of the
  classes ${\cal F}_{q,r}^{\psi}$ in the spaces ${\mathcal S}^p$}\label{Oder_Estimates_Sp}

{\bf \ref{Oder_Estimates_Sp}.1.} Estimates of the approximation characteristics of classes ${\mathcal F}_{q,r}^{\psi}$ and methods for obtaining them  depend significantly on the asymptotic behavior of functions $\psi$. Depending on this, we further distinguish the following different subsets of   $\psi$ with common properties and formulate the corresponding results.

  First, denote by  $B$  the set of all positive  non-increasing functions $\psi(t)$, $t\ge 1$,
 which  satisfy the so-called $\Delta_2$-condition,  i.e., for all  $t\ge 1$
   \begin{equation}\label{Eq:OE_Sp01}
     1<\tsfrac{\psi(t)}{\psi(2 t)}\le K_3.
    \end{equation}

Natural representatives of the set $B$ are functions of the form: $\psi(t)\equiv c$,
$\psi(t)=\ln^\varepsilon (t+{\rm e})$ for $\varepsilon<0$,  $\psi(t)=t^{-r}\ln^\varepsilon (t+{\rm e})$ for $r>0$ and $\varepsilon \in {\mathbb R}$, etc.

\begin{theorem}\label{Th:_OE_Sp1}
 Assume that  $0<p,q,r\le\infty$,
$\psi\in B$ and in the case  $p<q$, moreover, for all\ \  $t$, larger than a certain number  $t_0$, $\psi(t)$ is convex and satisfies the condition
\begin{equation}\label{Eq:OE_Sp02}
\tsfrac{t|\psi'(t)|}{\psi(t)}\ge \beta, \qquad \psi'(t):=\psi'(t+),
\end{equation}
 with a certain $\beta>d(\frac 1p-\frac 1q)$. Then
 \begin{equation}\label{Eq:OE_Sp03}
   {\mathscr D}_m({\cal F}_{q,r}^{\psi})_{_{\scriptstyle {\mathcal S}^p }}
   \asymp \left\{\begin{matrix}\psi(m^{\frac 1d}),\qquad \hfill 0<q\le p\le \infty, \\{\psi(m^{\frac 1d})}{m^{\frac 1p-\frac 1q}},\qquad \hfill 0<p<q\le \infty.\end{matrix}\right.
\end{equation}
and
 \begin{equation}\label{Eq:OE_Sp04}
\sigma_m({\mathcal F}_{q,r}^{\psi})_{_{\scriptstyle {\mathcal S}^p}}\asymp
{\psi(m^{\frac 1d})}{m^{\frac 1p-\frac 1q}},\qquad 0<p,q\le\infty.
\end{equation}
\end{theorem}

\begin{remark}\label{Rem:_OE_Sp1}
Note that condition \eqref{Eq:OE_Sp02}    guarantees that  $\sum_{k\in {\mathbb Z}^d} |\Psi_k|^\frac{pq}{q-p}<\infty$.
 \end{remark}

\banf Indeed, in this case, for all  $\tau\ge t_0$,
 \begin{equation}\label{Eq:OE_Sp06}
 \tsfrac{|\psi'(\tau)|}{\psi(\tau)}\ge \tsfrac{\beta}{\tau}.
\end{equation}
Integrating each part of this relation in the range from $t_0$ to $t$, $t\,{>}t_0$, we obtain
the estimate $\psi(t)\,{\ll}\, t^{-\beta}$, where $\beta> \frac {d(q-p)}{pq}$.

Therefore,  taking into account \eqref{Eq:ExV_BnA07} and \eqref{Eq:ExV_BnA11}, we conclude that
 \begin{equation}\label{Eq:OE_Sp08}
 \sum_{k\in {\mathbb Z}^d} |\Psi_k|^\frac{pq}{q-p}=\sum\limits_{n=1}^\infty{\nu_n}\psi^\frac{pq}{q-p}(n) \ll \sum\limits_{n=1}^{\infty}{n^{d-1}}\psi^\frac{pq}{q-p}(n)
 \ll  \sum\limits_{n=1}^{\infty} n^{d-1} n^{-\frac{pq \beta}{q-p}} < \infty.
\end{equation} \bend


 {\bf \ref{Oder_Estimates_Sp}.2.} To obtain similar estimates in the case where the
 functions $\psi$ tend to zero faster than any power function, we give the following definitions.

 Consider the set  ${\mathfrak M}$
 of all positive convex functions $\psi(t)$,  $t\ge 1$ such that
 \begin{equation}\label{Eq:OE_Sp09}
 \lim\limits_{t\to \infty} \psi(t)=0.
 \end{equation}
 Denote by ${\mathfrak M}'_{\infty}$ the set of all functions $\psi\in {\mathfrak M}$, satisfying conditions
 \begin{equation}\label{Eq:OE_Sp10}
   \alpha(\psi,t):=\tsfrac{\psi(t)}{t|\psi'(t)|}\downarrow 0, \qquad \psi'(t):=\psi'(t+),
   \end{equation}
  and $\frac{\psi(t)}{|\psi'(t)|}\uparrow \infty$.  Denote by ${\mathfrak M}_{\infty}^c$
  the set of all functions $\psi\in {\mathfrak M}$  satisfying   (\ref{Eq:OE_Sp10}) and
 \begin{equation}\label{Eq:OE_Sp11}
   K_4\le \tsfrac{\psi(t)}{|\psi'(t)|}\le K_5\qquad \forall t\ge 1.
  \end{equation}
  Finally, denote by  ${\mathfrak M}''_{\infty}$ the set of all   $\psi\in {\mathfrak M}$  satisfying condition  $\frac{\psi(t)}{|\psi'(t)|}\downarrow 0$.

Natural representatives of the sets ${\mathfrak M}'_{\infty}$,  ${\mathfrak M}_{\infty}^c$ and  ${\mathfrak M}''_{\infty}$ are functions of the form ${\rm exp}(-a t^s)$, $a>0$,  in cases where $s\in (0,1)$, $s=1$ and  $s>1$, respectively.

 {\bf \ref{Oder_Estimates_Sp}.3.} For  any positive integers  $d$ and $m$ and any $0<r\le\infty$, we denote
 \begin{equation}\label{Eq:OE_Sp12}
 \widetilde{n}_m:=\widetilde{n}_m(r,d)=(m/M_{r,d})^\frac 1d,
 \end{equation}
where $M_{r,d}$ is the number defined by the relation \eqref{Eq:ExV_BnA09}.

\begin{theorem}\label{Th:_OE_Sp2}
Assume that  $0<p,q,r\le\infty$ and $\psi\in {\mathfrak M}'_{\infty}\cup {\mathfrak M}^c_{\infty}$. Then
\begin{equation}\label{Eq:OE_Sp13}
   {\mathscr D}_m({\cal F}_{q,r}^{\psi})_{_{\scriptstyle {\mathcal S}^p }}
   \asymp \left\{\begin{matrix}\psi(\widetilde{n}_m),\qquad \hfill 0<q\le p\le \infty, \\
   \psi(\widetilde{n}_m) \big(m \alpha(\psi,\widetilde{n}_m)\big)^{\frac 1p-\frac 1q},\qquad \hfill 0<p<q\le \infty,\end{matrix}\right.
\end{equation}
where  $\widetilde{n}_m$ is defined by  \eqref{Eq:OE_Sp12},  and
 \begin{equation}\label{Eq:OE_Sp14}
\sigma_m({\cal F}_{q,r}^{\psi})_{_{\scriptstyle {\mathcal S}^p}}\asymp
  \psi(\widetilde{n}_m) \big(m \alpha(\psi,\widetilde{n}_m)\big)^{\frac 1p-\frac 1q},\qquad 0<p,q\le\infty.
  \end{equation}

\end{theorem}

 Since for any $\psi\in {\mathfrak M}$ satisfying \eqref{Eq:OE_Sp10},
 inequality \eqref{Eq:OE_Sp06} holds for any $\beta>0$ and sufficiently large $\tau$, then
  it can be similarly proven that for any $\psi$ from the sets ${\mathfrak M}'_{\infty}$,
 ${\mathfrak M}^c_{\infty}$ or ${\mathfrak M}''_{\infty}$, the series
 $\sum_{k\in {\mathbb Z}^d} |\Psi_k|^s$ with any $s>0$  are convergent.

In the case where $d=1$, the classes ${\cal F}_{q,r}^{\psi}({\mathbb T}^1)=:{\cal F}_{q}^{\psi}$ do not depend on $r$,  the constant $M_{r,1}=2$. Therefore, for any function $\psi\in {\mathfrak M}_{\infty}'\cup{\mathfrak M}_{\infty}^c$,
 \[
   {\mathscr D}_m({\cal F}_{q}^{\psi})_{_{\scriptstyle {\mathcal S}^p({\mathbb T}^1) }}
   \asymp \left\{\begin{matrix}\psi(\tsfrac m2),\qquad \hfill 0<q\le p\le \infty, \\
   {\psi(\tsfrac m2)}{(m \alpha(\psi,\tsfrac m2))^{\frac 1p-\frac 1q}},\qquad \hfill 0<p<q\le \infty,\end{matrix}\right.
   \eqno(\ref{Eq:OE_Sp13}')
 \]
and
 \[
\sigma_m({\cal F}_{q}^{\psi})_{_{\scriptstyle {\mathcal S}^p({\mathbb T}^1)}}\asymp
  {\psi(\tsfrac m2)}{(m \alpha(\psi,\tsfrac m2))^{\frac 1p-\frac 1q}},\qquad 0<p,q\le\infty.
  \eqno(\ref{Eq:OE_Sp14}')
 \]



 {\bf \ref{Oder_Estimates_Sp}.4.} Consider the case when $\psi$ belongs to the set ${\mathfrak M}''_{\infty}$.
 First of all, note that the following statement follows from relations \eqref{Eq:ExV_BaW05}, \eqref{Eq:ExV_BnA06}
 and \eqref{Eq:ExV_BnA07}:

\begin{proposition}\label{Prop:_OE_Sp2}
   Assume that $0<r\le \infty$, $m\in {\mathbb N}$ and $\psi=\psi(t)$, $t\ge 1$, is a
   positive non-increasing function.  Then for any   $0<q\le p\le \infty$ the following relation holds:
  \begin{equation}\label{Eq:OE_Sp15}
   {\mathscr D}_m({\cal F}_{q,r}^{\psi})_{_{\scriptstyle {\mathcal S}^p }}=\sigma_m({\cal F}_{\infty,r}^{\psi})_{_{\scriptstyle {\mathcal S}^\infty}}=\psi(n_{m+1}),
  \end{equation}
  where the number $n_{m+1}$ is defined by \eqref{Eq:ExV_BnA12} for $s=m+1$.
\end{proposition}

Let, as above, the numbers $V_s=V_s(d)$, $s=0,1,\ldots$, be defined by \eqref{Eq:ExV_BnA08}.


\begin{theorem}\label{Th:_OE_Sp3}
  Assume that $0<r\le \infty$, $m\in [V_{s-1},V_s)$, $s\in {\mathbb N}$ and  $\psi\in {\mathfrak M}''_{\infty}$.

   {\rm (i)} Let $0<p<q\le\infty$  and the function $\psi$ also satisfies the condition
     \begin{equation}\label{Eq:OE_Sp16}
      \lim\limits_{t\to\infty} \tsfrac{t^{\beta}\psi(t+1)}{\psi(t)}=0
    \end{equation}
   with $\beta=\tsfrac{d-1}{p}$. Then
    \begin{equation}\label{Eq:OE_Sp17}
    \sigma_m({\cal F}_{q,r}^{\psi})_{_{\scriptstyle {\mathcal S}^p}}\asymp
    \psi(s)\tsfrac{(V_{s}-m)^{\frac 1p}} {m^{\frac{d-1}{qd}}}.
    \end{equation}

   {\rm (ii)} Let $0<q\le p<\infty$  and $\psi$  satisfies condition \eqref{Eq:OE_Sp16} with $\beta =\tsfrac{d-1}{q}$. Then
   \begin{equation}\label{Eq:OE_Sp18}
  \sigma_m ({\mathcal F}_{q,r}^{\psi})_{_{\scriptstyle {\mathcal S}^p}}\asymp
   \tsfrac {\psi (s)}{(m+1-V_{s-1})^{\frac 1q-\frac 1p}},
   \end{equation}
    provided that $m=V_{s-1}$ or that the following additional condition is satisfied:
   \begin{equation}\label{Eq:OE_Sp19}
  p (V_{s}- m)\ge q(V_{s}-V_{s-1}),
  \end{equation}
   and   relation \eqref{Eq:OE_Sp17} holds,    provided that the   condition \eqref{Eq:OE_Sp19} is not satisfied.

   {\rm (iii)} If $0<q<p=\infty$ and  $\psi$  satisfies  \eqref{Eq:OE_Sp16} with   $\beta =\tsfrac{d-1}{q}$,
   then   \eqref{Eq:OE_Sp18} holds.

    {\rm (iv)}  If $0<p<q\le\infty$ and  $\psi$  satisfies \eqref{Eq:OE_Sp16}  with  $\beta =(d-1)\big(\tsfrac 1p-\tsfrac 1q\big)$, then
   \begin{equation}\label{Eq:OE_Sp20}
        {\mathscr D}_m({\cal F}_{q,r}^{\psi})_{_{\scriptstyle {\mathcal S}^p }} \asymp
        \psi(s) (V_{s}-m)^{\frac 1p-\frac 1q}.
    \end{equation}

\end{theorem}

\begin{remark}\label{Rem:_OE_Sp3}
 For any  $\psi\in {\mathfrak M}''_\infty$,  condition \eqref{Eq:OE_Sp16} is equivalent to the following:
   \begin{equation}\label{Eq:OE_Sp21}
     \lim\limits_{t\to \infty}\Big(\tsfrac{|\psi'(t)|}{\psi(t)}- \beta \ln t\Big)=+\infty.
    \end{equation}
    Therefore, for $d=1$, condition \eqref{Eq:OE_Sp16} holds for any $\psi\in {\mathfrak M}''_\infty$.
 \end{remark}

\banf  Indeed, if $\psi\in {\mathfrak M}''_\infty$ and  \eqref{Eq:OE_Sp21} holds, then for any $M>0$ and
sufficiently large $t$, we have
   \[
 \ln  \tsfrac{\psi(t)}{\psi(t+1)}=\int_t^{t+1}\tsfrac{|\psi'(\tau)|}{\psi(\tau)}{\rm d}\tau>
 \tsfrac{|\psi'(t)|}{\psi(t)}>\beta \ln t+M.
  \]
 Therefore, relation \eqref{Eq:OE_Sp16} is satisfied.

On the other hand side, if  $\psi\in {\mathfrak M}''_\infty$ and  \eqref{Eq:OE_Sp16}  is satisfied, then
for any $M>0$ and sufficiently large $t$, we have
\[
 M\le \ln \tsfrac {\psi(t)}{\psi(t+1)}-\beta \ln (t+1)
  =\int_t^{t+1}\tsfrac{|\psi'(\tau)|}{\psi(\tau)}{\rm d}\tau
  -\beta \ln (t+1)\le
  \tsfrac{ |\psi'(t+1)|}{\psi(t+1)}  -\beta \ln (t+1),
 \]
and relation  \eqref{Eq:OE_Sp21} holds. \bend


\begin{remark}\label{Rem:_OE_Sp4}
   From  \eqref{Eq:OE_Sp17} and \eqref{Eq:ExV_BnA10},
   it follows that under the conditions when relation \eqref{Eq:OE_Sp17} holds, we have
    \begin{equation}\label{Eq:OE_Sp22}
    \tsfrac{\psi(s)} {m^{\frac{d-1}{qd}}} \ll\sigma_m({\cal F}_{q,r}^{\psi})_{_{\scriptstyle {\mathcal S}^p}}\ll
     \psi(s) m^{\frac{d-1}{d}(\frac 1p-\frac 1q)}.
    \end{equation}

\end{remark}

Note that in the case when  $r=\infty$, for any $s \in {\mathbb N}$ we have
$V_s=(2s+1)^d$. Therefore, if $m\in [V_{s-1}, V_s)$, the number  $s$ is defined by the equality $s=[\frac{(m+1)^{\frac 1d}}2]$.  Here and below, $[x]$ denotes the greatest integer in   $x$.

If $d=1$, the classes  ${\cal F}_{q,r}^{\psi}=:{\cal F}_{q}^{\psi}$  does not depend on
 $r$, and for any  $m\in [V_{s-1}, V_s)$ we have $m=V_{s-1}=V_{s}-1$ and $s=[\tsfrac{m+1}2]$.
 Therefore, for any   $\psi\in {\mathfrak M}_{\infty}''$,
 \begin{equation}\label{Eq:OE_Sp23}
 \sigma_m({\cal F}_{q}^{\psi})_{_{\scriptstyle {\mathcal S}^p({\mathbb T}^1)}}\asymp
 {\mathscr D}_m({\cal F}_{q}^{\psi})_{_{\scriptstyle {\mathcal S}^p({\mathbb T}^1)}}\asymp \psi([\tsfrac{m+1}2]),\qquad 0<p,\,q\le \infty.
 \end{equation}

If  $d>1$, the   obtained estimates depend significantly on the placement of the number $m$ on the half-segment $[V_{s-1}, V_{s})$. Considering some specific subsequences $m(s)$ in Theorem \ref{Th:_OE_Sp3}, we obtain the following corollary.


\begin{corollary}\label{Cor:_OE_Sp3}
  Assume that $d>1$, $m\in [V_{s-1},V_s)$, $s\in {\mathbb N}$, the parameters   $p$, $q$, $r$ and $\psi$ satisfy condition of Theorem \ref{Th:_OE_Sp3}.

    {\rm (i)} Let $m=m(s)=V_{s}-c_s$, $1\le c_s\le c$. Then  for   $0<p,\,q\le\infty$
    \begin{equation}\label{Eq:OE_Sp24}
    \sigma_m({\cal F}_{q,r}^{\psi})_{_{\scriptstyle {\mathcal S}^p}}\asymp
    \tsfrac{\psi(s)}{m^{\frac{d-1}{qd}}},
    \end{equation}
   and for    $0<p<q\le\infty$
   \begin{equation}\label{Eq:OE_Sp25}
        {\mathscr D}_m({\cal F}_{q,r}^{\psi})_{_{\scriptstyle {\mathcal S}^p }} \asymp
        \psi(s).
    \end{equation}

    {\rm (ii)}  Let  the sequence $m=m(s)$ be such that
    \begin{equation}\label{Eq:OE_Sp26}
      V_{s}-m(s)\asymp V_{s}-V_{s-1}.
    \end{equation}
    Then  for any  $0<p<q\le\infty$
     \begin{equation}\label{Eq:OE_Sp27}
    \sigma_m({\cal F}_{q,r}^{\psi})_{_{\scriptstyle {\mathcal S}^p}}\asymp
    {\psi(s)}{m^{\frac{d-1}{d}(\frac 1p-\frac 1q)}}
    \end{equation}
    and
    \begin{equation}\label{Eq:OE_Sp28}
        {\mathscr D}_m({\cal F}_{q,r}^{\psi})_{_{\scriptstyle {\mathcal S}^p }}  \asymp {\psi(s)}{m^{\frac{d-1}{d}(\frac 1p-\frac 1q)}}.
    \end{equation}
   If   $0<q\le p<\infty$, then relation  \eqref{Eq:OE_Sp27} holds, in particular,
   when $m(s)-V_{s-1}\to \infty$ as $s\to \infty$.
   If $m=m(s)=V_{s-1}+c_s$, $0\le c_s\le c$,  then  for    $0<q< p<\infty$
     \begin{equation}\label{Eq:OE_Sp29}
    \sigma_m({\cal F}_{q,r}^{\psi})_{_{\scriptstyle {\mathcal S}^p}}\asymp
     \psi(s).
    \end{equation}
    In the case $0<p=q<\infty$, if $m=m(s)=V_{s-1}+c_s$,   $1\le c_s\le c$,  then
    relation \eqref{Eq:OE_Sp27} holds,    and if $m=m(s)=V_{s-1}$, then relation \eqref{Eq:OE_Sp29} holds.
\end{corollary}

\begin{remark}\label{Rem:_OE_Sp5}
  Comparing the obtained estimates, we conclude the following:
 \begin{itemize}
 \item[$\bullet$]  If  the function  $\psi$ belongs to the set $B$ and satisfies condition \eqref{Eq:OE_Sp02}
      with a certain $\beta>d(\frac 1p-\frac 1q)$ or to the set
      ${\mathfrak M}'_{\infty}\cup {\mathfrak M}^c_{\infty}$, then for $0<p\le q\le \infty$
       \begin{equation}\label{Eq:OE_Sp30}
       \sigma_m({\cal F}_{q,r}^{\psi})_{_{\scriptstyle {\mathcal S}^p}}\asymp {\mathscr D}_m({\cal F}_{q,r}^{\psi})_{_{\scriptstyle {\mathcal S}^p }}.
       \end{equation}
 \item[$\bullet$] If $\psi\in B\cup {\mathfrak M}'_{\infty}$ and  $0<q< p\le \infty$,  then
      \begin{equation}\label{Eq:OE_Sp31}
      \sigma_m({\cal F}_{q,r}^{\psi})_{_{\scriptstyle {\mathcal S}^p}}=o\big( {\mathscr D}_m({\cal F}_{q,r}^{\psi})_{_{\scriptstyle {\mathcal S}^p }}\big),\qquad m\to \infty.
      \end{equation}
 \item[$\bullet$] If  $\psi\in {\mathfrak M}^c_{\infty}$ and $0<q< p\le \infty$, then \eqref{Eq:OE_Sp31}
      holds for $d>1$, and
      \begin{equation}\label{Eq:OE_Sp32}
      \sigma_m({\cal F}_{q}^{\psi})_{_{\scriptstyle {\mathcal S}^p({\mathbb T}^1)}}\asymp {\mathscr D}_m({\cal F}_{q}^{\psi})_{_{\scriptstyle {\mathcal S}^p ({\mathbb T}^1)}}.
      \end{equation}
 \item[$\bullet$] If  $\psi\in {\mathfrak M}''_{\infty}$, then for any $0<p,\,q\le \infty$
   relation \eqref{Eq:OE_Sp32} holds. If  $d>1$, then  relation \eqref{Eq:OE_Sp30} holds, in particular,
    when the sequence $m=m(s)$ satisfies  \eqref{Eq:OE_Sp26}
        and $0<p<q\le\infty$,  when $p=q=\infty$, or when $m=m(s)=V_{s-1}+c_s$, $0\le c_s\le c$, and    $0<q< p<\infty$.

 \end{itemize}

 \end{remark}

Note that for finite values $0<p,q<\infty$, the statements of Theorems \ref{Th:_OE_Sp1},
\ref{Th:_OE_Sp2} and \ref{Th:_OE_Sp3}  were proven in \cite{Shidlich_2011, Shidlich_2013, Shidlich_2014,Shidlich_2016}.
 However, given the relative inaccessibility of these sources for English-speaking readers,
 we present these statements in this work with proof, slightly changing their presentation and structure.

 In the case when  $\psi(t)=1$, relation \eqref{Eq:OE_Lp06} also follows  from
 Lemma 4.4 \cite{Moeller_Stasyuk_Ullrich_2025}.   For $d=1$ and $0<p,\,q<\infty$, statements similar to Theorems \ref{Th:_OE_Sp1},
  \ref{Th:_OE_Sp2} and \ref{Th:_OE_Sp3} were proven in \cite{Shydlich_UMZh2009}.
  Exact order estimates of integral analogues of the quantities
  $\sigma_m({\mathcal F}_{q,r}^{\psi})_{_{\scriptstyle {\mathcal S}^p}}$ were found in
     \cite{Stepanets_Shidlich_IZV_2010}.

  \section{Approximation characteristics of the classes ${\cal F}_{q,r}^{\psi}$
  in the spaces $L_p$}\label{Oder_Estimates_Lp}

 {\bf \ref{Oder_Estimates_Lp}.1.}
 In the case where $2\le p<\infty$, based on the Hausdorff-Young theorem
 (see, for example,  \cite[\S 0.1]{Temlyakov1993}), for any   $f\in L_p$, the following inequality holds:
  \begin{equation}\label{Eq:OE_Lp01}
 \|f\|_{_{\scriptstyle L_{p}}}\le  \|f\|_{_{\scriptstyle {\mathcal S}^{p'} }}.
 \end{equation}
Here and below, for any $1<p<\infty$, we set $p':=\frac{p}{p-1}$  and $p':=\infty$ when $p=1$.

If $1\le p<2$, then for any   $f\in L_p({\mathbb T}^d)$
   \begin{equation}\label{Eq:OE_Lp02}
 \|f\|_{_{\scriptstyle L_{p} }}\le  \|f\|_{_{\scriptstyle L_{2} }}= \|f\|_{_{\scriptstyle {\mathcal S}^{2}}}.
 \end{equation}
 Thus, from the estimates for  the approximate quantities in the spaces ${\mathcal S}^p$ obtained in
 Subsection \ref{Oder_Estimates_Sp}, estimates from above of the similar quantities
 in the spaces $L_p$ also follow. Here, we   consider some of the cases
 in which the corresponding estimates from below are obtained.

{\bf \ref{Oder_Estimates_Lp}.2. } As mentioned above, in the case where $\psi$ is a power function, i.e.,
 $\psi(t)=t^{-s}$, $s>0$, for all $1\le p\le \infty$
 the exact order estimates of the quantities $\sigma_m({\mathcal F}_{q,\infty}^{\psi})_{_{\scriptstyle L_p}}$  and $G_m({\mathcal F}_{q,\infty}^{\psi})_{_{\scriptstyle L_p}}$ were obtained in
 \cite{DeVore_Temlyakov_1995} and \cite{Temlyakov_1998}, correspondingly.
 In particular, from Theorems 6.1 \cite{DeVore_Temlyakov_1995} and 3.1 \cite{Temlyakov_1998},
 it follows that for all $s>d(1-\frac 1q)_+$,
 \begin{equation}\label{Eq:OE_Lp03}
\sigma_m({\mathcal F}_{q,\infty}^{s})_{_{\scriptstyle L_p}}\asymp  m^{-\frac sd-\frac 1q+\frac 12}, \  \ \ \ 1\le p\le \infty,
 \end{equation}
and
\begin{equation}\label{Eq:OE_Lp04}
G_m({\mathcal F}_{q,\infty}^{s})_{_{\scriptstyle L_p}} \asymp \left\{\begin{matrix} m^{-\frac sd-\frac 1q+\frac 12},\  \ \ \ 1\le p<2, \\ m^{-\frac sd-\frac 1q+1-\frac 1p},\  \ \ \ 2\,\le\, p\,<\infty.\end{matrix}\right.
\end{equation}
From the following Theorem \ref{Th:_OE_Lp1}, in particular, it follows that
for  $\sigma_m({\mathcal F}_{q,\infty}^{\psi})_{_{\scriptstyle L_p}}$ and
$G_m({\mathcal F}_{q,\infty}^{\psi})_{_{\scriptstyle L_p}}$, the estimates of forms
(\ref{Eq:OE_Lp03}) and (\ref{Eq:OE_Lp04})   are satisfied for a wider set   of the functions $\psi$.

\begin{theorem}\label{Th:_OE_Lp1}
  Assume that $0< r,\,q\le \infty$, $1\le p<\infty$, $\psi\in B$ and in the case
  $\tsfrac p{p-1}<q$, moreover, for all\ \  $t$, larger than a certain number  $t_0$,  $\psi(t)$ is \
  convex and satisfies the condition \eqref{Eq:OE_Sp02}
  with $\beta>d(\frac 12-\frac 1q)$ when $1< p\le 2$ and $\beta>d(1-\frac 1p-\frac 1q)$ when $2\le p<\infty$.
  Then
 \begin{equation}\label{Eq:OE_Lp05}
  G_m({\mathcal F}_{q,r}^{\psi})_{_{\scriptstyle L_p}}\asymp
  \sigma_m^\perp({\mathcal F}_{q,r}^{\psi})_{_{\scriptstyle L_p}}\asymp
  \left\{\begin{matrix}  {\psi(m^{\frac 1d})}{m^{\frac 12-\frac 1q}},\qquad \ \ \ 1\le p\le 2,
  \\ {\psi(m^{\frac 1d})}{m^{1-\frac 1p-\frac 1q}},\qquad  2 \le p<\infty\end{matrix},\right.
 \end{equation}
for all $1\le p\le 2$
\begin{equation}\label{Eq:OE_Lp06}
\sigma_m({\mathcal F}_{q,r}^{\psi})_{_{\scriptstyle L_p}}\asymp   {\psi(m^{\frac 1d})}{m^{\frac 12-\frac 1q}},
\end{equation}
and  for  all  $2<p<\infty$
 \begin{equation}\label{Eq:OE_Lp07}
 {\psi(m^{\frac 1d})}{m^{\frac 12-\frac 1q}}\ll \sigma_m({\mathcal F}_{q,r}^{\psi})_{_{\scriptstyle L_p}}\ll \psi(m^{\frac 1d}) m^{1-\frac 1p-\frac 1q}.
 \end{equation}
\end{theorem}

In the case $2<p\le \infty$, the following theorem is true.

\begin{theorem}\label{Th:_OE_Lp2}
  Assume that $1\le r\le \infty$, $2< p\le\infty$, $0<q< \infty$,
  the function $\psi$ belongs to the set $B$ and   for all\ \  $t$, larger than a
  certain number  $t_0$, $\psi(t)$ is convex and satisfies condition \eqref{Eq:OE_Sp02}
   with $\beta>d(1-\frac 1q)_+$. Then relation \eqref{Eq:OE_Lp06} holds.
 \end{theorem}

Note that conditions in Theorems \ref{Th:_OE_Lp1} and \ref{Th:_OE_Lp2} guarantee the embedding ${\mathcal F}_{q,r}^{\psi}{\subset} L_p$.

Putting $r=\infty$ and $\psi(t)=t^{-s}$,  from Theorems \ref{Th:_OE_Lp1} and \ref{Th:_OE_Lp2}  we obtain the following corollary:

\begin{corollary}\label{Cor: Lp1}
 Assume that $1\le p<\infty$, $0<q<\infty$, $s$ is a positive number, which in the case
 $\tsfrac p{p-1}<q$, satisfies the inequality $s>\beta$, where $\beta$ is defined
 in Theorem \ref{Th:_OE_Lp1}. Then   relation
 relation \eqref{Eq:OE_Lp04} holds for    $1\le p<\infty$,  and
  relation  \eqref{Eq:OE_Lp03} holds for    $1\le p\le 2$.
 If $s>d(1-\frac 1q)_+$, then relation \eqref{Eq:OE_Lp03} holds for  $1\le p\le\infty$.
\end{corollary}

This statement complements the above-mentioned results   of \cite{DeVore_Temlyakov_1995} and \cite{Temlyakov_1998}   in the following sense:

\begin{itemize}

 \item[$\bullet$] From Corollary \ref{Cor: Lp1} it follows that in the case $1 <q \le \tsfrac p{(p-1)}$, relation \eqref{Eq:OE_Lp03}  (for $1\le p\le 2$) and  relation \eqref{Eq:OE_Lp04}  (for $1\le p<\infty$) also hold for all $s>0$.

 \item[$\bullet$] If $1<p\le 2$ and $q>\tsfrac p{p-1}$, then relations \eqref{Eq:OE_Lp03} and \eqref{Eq:OE_Lp04}  also hold for all  $s$ such that $d(\frac 12-\frac 1q)<s\le d(1-\frac 1q)$.

 \item[$\bullet$] If $2<p<\infty$ and $q>\tsfrac p{p-1}$, then relation \eqref{Eq:OE_Lp04}  also holds for all  $s$ such that $d(1-\frac 1p-\frac 1q)<s\le d(1-\frac 1q)$.

      \item[$\bullet$] In the case
       $2<p\le \infty$, conditions on $s$ in Corollary \ref{Cor: Lp1} (for validity of relation \eqref{Eq:OE_Lp03}) are the same as in Theorem 6.1 \cite{DeVore_Temlyakov_1995}.
 \end{itemize}

  Note also that if $0<q\le \tsfrac p{p-1}$, then the conditions of Theorem \ref{Th:_OE_Lp1} are satisfied,
  for example, for the functions $\psi(t)=t^{-s}\ln^\varepsilon (t+{\rm e})$, where $s> 0$,
  $\varepsilon \in {\mathbb R}$,   $\psi(t)=\ln^\varepsilon (t+{\rm e})$,  $\varepsilon<0$,
  and  $\psi(t)\equiv c$. If $1<\tsfrac p{p-1}<q$ and $1<p\le 2$, then the conditions of
  Theorem \ref{Th:_OE_Lp1} are satisfied  for  $\psi(t)=t^{-s}\ln^\varepsilon (t+{\rm e})$,
  where $\varepsilon \in {\mathbb R}$ and $s>d(\frac 12-\frac 1q)$.
  If $1<\tsfrac p{p-1}<q$ and $2<p<\infty$, then the conditions of Theorem \ref{Th:_OE_Lp1} are satisfied
  for $\psi(t)=t^{-s}\ln^\varepsilon (t+{\rm e})$, where $\varepsilon \in {\mathbb R}$ and $s>d(1-\frac 1p-\frac 1q)$.
 The conditions of Theorem \ref{Th:_OE_Lp2} are satisfied  for $\psi(t)=t^{-s}\ln^\varepsilon (t+{\rm e})$, where $\varepsilon \in {\mathbb R}$ and $s>d(1-\frac 1q)_+$.


{\bf \ref{Oder_Estimates_Lp}.3.} Consider the case when the function $\psi$ decreases faster that any power
function.

 \begin{theorem}\label{Th:_OE_Lp_3}
 Assume that $m\in {\mathbb N}$, $0<r,\,q\le\infty$, $1\le p<\infty$, $n\in [V_{m-1}, V_m)$ and
 the function $\psi$ belongs to the set ${\mathfrak M}''_\infty$ and
 satisfies \eqref{Eq:OE_Sp16} with $\beta >\max\{\tsfrac{d-1}{p'},\tsfrac{d-1}{q}\}$.

{\rm (i)} If $m=m(s)=V_{s}-c_s$, $1\le c_s\le c$, then for any $0<q\le\infty$
    \begin{equation}\label{Eq:OE_Lp08}
    \sigma_m^\perp({\cal F}_{q,r}^{\psi})_{_{\scriptstyle L_p}}\asymp
    G_m({\cal F}_{q,r}^{\psi})_{_{\scriptstyle L_p}}\asymp
    \tsfrac{\psi(s)}{m^{\frac{d-1}{qd}}}.
    \end{equation}

{\rm (ii)} If $m=m(s)=V_{s-1}+c_s$, $0\le c_s\le c$, and   $0<q< p'<\infty$ or if $m=m(s)=V_{s-1}$ and
  $0<p'=q<\infty$, then
     \begin{equation}\label{Eq:OE_Lp09}
    \sigma_m^\perp({\cal F}_{q,r}^{\psi})_{_{\scriptstyle L_p}}\asymp
    G_m({\cal F}_{q,r}^{\psi})_{_{\scriptstyle L_p}}\asymp  \psi(s).
    \end{equation}
\end{theorem}

 As noted above, if $d=1$, the classes  ${\cal F}_{q,r}^{\psi}={\cal F}_{q}^{\psi}$
 does not depend on $r$, and $s=[\tsfrac{m+1}2]$ for any  $m\in [V_{s-1}, V_s)$.
 Therefore, the following corollary follows from Theorem \ref{Th:_OE_Lp_3}.

\begin{corollary}\label{Cor: Lp2}
 For any $1\le p<\infty$,  $0<q\le\infty$ and any function  $\psi\in {\mathfrak M}_{\infty}''$
 \begin{equation}\label{Eq:OE_Lp10}
  \sigma_m^\perp({\cal F}_{q}^{\psi})_{_{\scriptstyle L_{p}({\mathbb T}^1)}}\asymp
  G_m({\cal F}_{q}^{\psi})_{_{\scriptstyle L_p({\mathbb T}^1)}}\asymp   \psi([\tsfrac{m+1}2]).
\end{equation}
\end{corollary}

The similar estimate can be proven in the case when $d=1$ and  $\psi\in {\mathfrak M}_\infty^c$.

\begin{theorem}\label{Th:_OE_Lp_4}
 For any $1\le p<\infty$,  $0<q\le\infty$ and any function  $\psi\in {\mathfrak M}_{\infty}^c$
 \begin{equation}\label{Eq:OE_Lp11}
  \sigma_m^\perp({\cal F}_{q}^{\psi})_{_{\scriptstyle L_{p}({\mathbb T}^1)}}\asymp
  G_m({\cal F}_{q}^{\psi})_{_{\scriptstyle L_p({\mathbb T}^1)}}
  \asymp {\mathscr D}_m^\perp({\cal F}_{q,r}^{\psi})_{_{\scriptstyle L_{p}}}\asymp
  \psi(\tsfrac{m}2)\asymp  \psi([\tsfrac{m+1}2]).
\end{equation}
\end{theorem}

The following statement also holds for the quantities ${\mathscr D}_m^\perp({\cal F}_{q,r}^{\psi})_{_{\scriptstyle L_{p}}}$.
 \begin{theorem}\label{Th:_OE_Lp_5}
 Assume that $m\in {\mathbb N}$, $0<r,\,q\le\infty$, $1\le p<\infty$ and  $m\in [V_{s-1}, V_s)$.

{\rm (i)} If  $0<q\le  p'$, then for any positive non-increasing function $\psi$, we have
  \begin{equation}\label{Eq:OE_Lp12}
    {\mathscr D}_m^\perp({\cal F}_{q,r}^{\psi})_{_{\scriptstyle L_{p}}}=\psi(s).
  \end{equation}

{\rm (ii)}  If $1<p'<q\le \infty$ and  $\psi\in{\mathfrak M}_\infty''$ and  condition
\eqref{Eq:OE_Sp16} is satisfied  with  $\beta =\tsfrac {(d-1)(q-p')}{p'q}$, then for
the sequence $m=m(s)=V_{s}-c_s$, $1\le c_s\le c$, the following estimate holds:
  \begin{equation}\label{Eq:OE_Lp13}
   {\mathscr D}_m^\perp({\cal F}_{q,r}^{\psi})_{_{\scriptstyle L_{p}}}\asymp\psi(s).
 \end{equation}

 \end{theorem}

\begin{corollary}\label{Cor: Lp3}
 Let $1\le p<\infty$. For
 any positive non-increasing function $\psi$ and   $0<q\le  \tsfrac {p}{p-1}$,
$$
{\mathscr D}_m^\perp({\cal F}_{q}^{\psi})_{_{\scriptstyle L_{p}({\mathbb T}^1)}}= \psi([\tsfrac{m+1}2]).
\eqno(\ref{Eq:OE_Lp12}')
$$
If $1<p'<q\le \infty$,  then for any  $\psi\in {\mathfrak M}_{\infty}''$,
  $$
{\mathscr D}_m^\perp({\cal F}_{q}^{\psi})_{_{\scriptstyle L_{p}({\mathbb T}^1)}}\asymp   \psi([\tsfrac{m+1}2]).
\eqno(\ref{Eq:OE_Lp13}')
$$
\end{corollary}

In the case when $\psi\in B$, estimates for ${\mathscr D}_m^\perp({\cal F}_{q,r}^{\psi})_{_{\scriptstyle L_{p}}}$
are given by the following  statement.

\begin{theorem}\label{Th:_OE_Lp_6} Assume that conditions of Theorem  \ref{Th:_OE_Lp1} are satisfied. Then
\begin{equation}\label{Eq:OE_Lp14}
{\mathscr D}_m^\perp({\cal F}_{q,r}^{\psi})_{_{\scriptstyle L_{p}}}\asymp
\left\{\begin{matrix}\psi(m^{\frac 1d}),\quad \hfill &  0<q\le \tsfrac p{p-1},\ 1\le p<\infty,\cr
{\psi(m^{\frac 1d})}{m^{\frac 12-\frac 1q}},\quad \hfill &  \tsfrac p{p-1}<q\le\infty,\ 1< p<2, \cr
{\psi(m^{\frac 1d})}{m^{1-\frac 1p-\frac 1q}},\quad \hfill &  \tsfrac p{p-1}<q\le\infty,\ 2\le p<\infty.\end{matrix}\right.
\end{equation}
\end{theorem}

\vskip -5mm

 \begin{remark}\label{Rem:_Lp_1}
 Analyzing the results of Section \ref{Oder_Estimates_Lp},  we conclude the following:
 \begin{itemize}
 \item[$\bullet$] If $d=1$  and  $\psi\in {\mathfrak M}''_{\infty}\cup {\mathfrak M}_\infty^c$, then
 for  any $1\le p<\infty$ and $0<q\le \infty$
 \[
 \sigma_m^\perp({\cal F}_{q}^{\psi})_{_{\scriptstyle L_{p}({\mathbb T}^1)}}\asymp G_n({\cal F}_{q}^{\psi})_{_{\scriptstyle L_p({\mathbb T}^1)}} \asymp  {\mathscr D}_m^\perp({\cal F}_{q}^{\psi})_{_{\scriptstyle L_{p}({\mathbb T}^1)}}.
 \]
  \item[$\bullet$] If $d>1$  and  condition of Theorems \ref{Th:_OE_Lp_3}
  and \eqref{Th:_OE_Lp_5} are satisfied, then the similar estimate
  \begin{equation}\label{Eq:OE_Lp15}
  \sigma_m^\perp({\cal F}_{q,r}^{\psi})_{_{\scriptstyle L_{p}}}\asymp
  G_m({\cal F}_{q,r}^{\psi})_{_{\scriptstyle L_{p}}}\asymp
  {\mathscr D}_m^\perp({\cal F}_{q,r}^{\psi})_{_{\scriptstyle L_{p}}}
  \end{equation}
 holds  when  $m=m(s)=V_{s-1}+c_s$, $0\le c_s\le c$, and   $0<q< p'<\infty$, or  when $m=m(s)=V_{s-1}$ and
  $0<p'=q<\infty$.   If $m=m(s)=V_{s}-c_s$, $1\le c_s\le c$, then for any $0<q\le\infty$
  \begin{equation}\label{Eq:OE_Lp16}%
  \sigma_m^\perp({\cal F}_{q,r}^{\psi})_{_{\scriptstyle L_{p}}}\asymp {G}_n({\cal F}_{q,r}^{\psi})_{_{\scriptstyle L_{p}}}=o\Big({\mathscr D}_m^\perp({\cal F}_{q,r}^{\psi})_{_{\scriptstyle L_{p}}}\Big).
 \end{equation}
  \item[$\bullet$] In the case when conditions of Theorem  \ref{Th:_OE_Lp1} are satisfied,
 for all   $p'\le q\le\infty$ relation \eqref{Eq:OE_Lp15} holds, and  for all  $0<q\le p'$, $1\le p<\infty$
 relation \eqref{Eq:OE_Lp16} holds.

 \end{itemize}

\end{remark}

  Note that for finite values $0<q<\infty$, the statements of Theorems \ref{Th:_OE_Lp1}-\ref{Th:_OE_Lp_6}
  were  were proven in \cite{Shidlich_2011, Shidlich_2013, Shidlich_2014, Shidlich_2016}. In the case when
  $\psi(t)=1$, $0<q\le 1$ and $2\le p<\infty$, the upper estimate for the quantity
  $\sigma_m({\mathcal F}_{q,r}^{\psi})_{_{\scriptstyle L_p}}$ can be also obtained from
  Theorem 4.4 \cite{Moeller_Stasyuk_Ullrich_2024}.

\section{Proof of theorems}\label{Proofs}

\subsection{Proof of Theorem \ref{Th:ApChFS01}}

  Consider  the case when $0<q<p=\infty$ or $p=q=\infty$. Due to \eqref{Eq:ApCh01}, \eqref{Eq:ApCh02},  \eqref{Eq:Intr01}, \eqref{Eq:Intr001} and \eqref{Eq:ExV_BaW03}, for any $f\in {\mathcal F}_{q}^{\Psi}$ and $\gamma_m\in \Gamma_m$, we have
      \begin{equation}\label{Eq:PrTh1_01}
  E_{\gamma_m}(f)_{_{\scriptstyle {\mathcal S}^p}} =
  {\mathscr E}_{\gamma_m}(f)_{_{\scriptstyle {\mathcal S}^p}}=\|\,f - S_{\gamma_m}(f)
  \|_{_{\scriptstyle {\mathcal S}^p}}=\|\{ \Psi_{\gamma_m}(k)\widehat{f} ({ k})/\Psi_k\}_{k\in {\mathbb Z}^d}\|_{l_p({\mathbb Z}^d)}.
     \end{equation}
Therefore, for any $f\in {\mathcal F}_{q}^{\Psi}$ and $\gamma_m\in \Gamma_m$,   the following upper estimate holds:
      \begin{equation}\label{Eq:PrTh1_02}
  E_{\gamma_m}(f)_{_{\scriptstyle {\mathcal S}^\infty}}
   =\sup_{k\in {\mathbb Z}^d} \Big|\tsfrac{\Psi_{\gamma_m}(k)\widehat{f} ({ k})}{\Psi_k}\Big|
   \le \bar{\Psi}_{\gamma_m}(1)\sup_{k\in {\mathbb Z}^d} \Big|\tsfrac{\widehat{f} ({ k})}{\Psi_k}\Big|
  $$
  $$
  \le
  \bar{\Psi}_{\gamma_m}(1)\|\{ \widehat{f} ({ k})/\Psi_k\}_{k\in {\mathbb Z}^d}\|_{l_q({\mathbb Z}^d)}
  \le \bar{\Psi}_{\gamma_m}(1).
     \end{equation}
Let $k_0\in {\mathbb Z}^d$ be a  vector such that $|\Psi_{k_0}|=\bar{\Psi}_{\gamma_m}(1)$. The function $f_0(x):=\Psi_{k_0}e_{k_0}(x)$ belongs to the set ${\mathcal F}_{q}^{\Psi}$  and
 \begin{equation}\label{Eq:PrTh1_03}
  E_{\gamma_m}(f_0)_{_{\scriptstyle {\mathcal S}^p}} =
  {\mathscr E}_{\gamma_m}(f_0)_{_{\scriptstyle {\mathcal S}^p}} =|\Psi_{k_0}|=\bar{\Psi}_{\gamma_m}(1).
\end{equation}
 Combining \eqref{Eq:PrTh1_02} and \eqref{Eq:PrTh1_03}, we see that  in this case relation \eqref{Eq:ExV_BaW04} holds.

 Now, assume that $0<p<q=\infty$ and the series $\sum_{k\in {\mathbb Z}^d} |\Psi_k|^p$ converges. By virtue of \eqref{Eq:PrTh1_01}, for any $f\in {\mathcal F}_{\infty}^{\Psi}$ and $\gamma_m\in \Gamma_m$, we have
 \begin{equation}\label{Eq:PrTh1_04}
  E_{\gamma_m}(f)_{_{\scriptstyle {\mathcal S}^p}}^p =
   \sum_{k\in {\mathbb Z^d}} \big|\tsfrac{\Psi_{\gamma_m}(k)\widehat{f} ({ k})}{\Psi_k}\big|^p\le
   \sup_{k\in {\mathbb Z}^d} \Big|\tsfrac{\widehat{f} ({ k})}{\Psi_k}\Big|^p\sum_{k\in {\mathbb Z^d}} |\Psi_{\gamma_m}(k)|^p \le \sum_{j=1}^{\infty}  \bar{\Psi}_{\gamma_m}^{p}(j).
  \end{equation}
  Since
  \[
  \sum_{k=1}^{\infty}  \bar{\Psi}_{\gamma_m}^{p}(k)\le \sum_{k=1}^{\infty}  \bar{\Psi}_{k}^{p}=
  \sum_{k\in {\mathbb Z}^d} |\Psi_k|^p<\infty,
  \]
  for any $\varepsilon>0$ there exists a number $N_\varepsilon$ such that
  $
  \sum_{j=N_\varepsilon+1}^{\infty}  \bar{\Psi}_{\gamma_m}^{p}(j)<\varepsilon.
  $
  Consider any collection $\gamma_{N_\varepsilon}^*=\{k_1,\ldots, k_{N_\varepsilon}\}$   from the set $\Gamma_{N_\varepsilon}\setminus \gamma_m$ such that
 \[
    \gamma_{N_\varepsilon}^*=\{k_j\in {\mathbb Z}^d:\quad |\Psi_{k_j}|=\bar{\Psi}_{\gamma_m}(j),
    \quad j=1,2,\ldots,m\}
 \]
 Then the function $f_\varepsilon (x):=\sum_{k\in \gamma_{N_\varepsilon}^*} \Psi_k e_k(x)$ belongs to the set ${\mathcal F}_{\infty}^{\Psi}$ and
 \begin{equation}\label{Eq:PrTh1_05}
  E_{\gamma_m}(f_\varepsilon)_{_{\scriptstyle {\mathcal S}^p}}^p =
   \sum_{k\in \gamma_{N_\varepsilon}^*} |\Psi_k|^p=\sum_{j=1}^{N_\varepsilon}  \bar{\Psi}_{\gamma_m}^{p}(j)>\sum_{j=1}^{\infty}  \bar{\Psi}_{\gamma_m}^{p}(j)-\varepsilon.
  \end{equation}
 Combining  \eqref{Eq:PrTh1_04} and  \eqref{Eq:PrTh1_05} and taking into account the arbitrariness of $\varepsilon$, we see that \eqref{Eq:ExV_BaW10} is indeed true.

 Finally,  considering the infima of all quantities  \eqref{Eq:ExV_BaW04} and \eqref{Eq:ExV_BaW10} over all $\gamma_m\in \Gamma_m$, we see that in the corresponding cases, relations \eqref{Eq:ExV_BaW07} and \eqref{Eq:ExV_BaW11} are true, and for any collection $\gamma_m^*\in\Gamma_m$ satisfying \eqref{Eq:ExV_BaW08}, relation  \eqref{Eq:ExV_BaW09}  holds. \bend


\subsection{Auxiliary estimates}\label{Aux_St_01}

Before proving the theorems on   estimates of best $n$-term approximations  and basis widths   in the spaces ${\mathcal S}^p$, we give several auxiliary statements.

 {\bf \ref{Aux_St_01}.1.} First, we show how the belonging of a function to the sets
 $B$, ${\mathfrak M}'_{\infty}$, ${\mathfrak M}^c_{\infty}$ or ${\mathfrak M}''_{\infty}$
 affects its values at different points.

As follows from \eqref{Eq:OE_Sp01}, for any $\psi\in B$ and $c>1$, we have
 \begin{equation}\label{Eq:Aux000}
 \psi(t)\asymp \psi (ct).
  \end{equation}

In the case when $\psi\in {\mathfrak M}'_{\infty}\cup {\mathfrak M}^c_{\infty}$, we can obtain the following statement.

\begin{proposition}\label{Aux_Prop1}
Let $d\ge 1$ and $c>1$. Then for any   $\psi\in {\mathfrak M}'_{\infty}\cup {\mathfrak M}^c_{\infty}$
\begin{equation}\label{Aux_Eq01}
 \psi(t^{1/d})\asymp \psi((t+c)^{1/d})\qquad \mbox{and}\qquad  \tsfrac{\psi(t^ {1/d})}{|\psi'(t^ {1/d})|}  \asymp \tsfrac{\psi(ct^{1/d})}{|\psi'(ct^{1/d})|}.
\end{equation}
\end{proposition}

\banf  If $\psi$ belongs  to  ${\mathfrak M}'_{\infty}$ or ${\mathfrak M}^c_{\infty}$, then $\tsfrac{|\psi'(t)|}{\psi(t)}\le K_6$  for any $t\ge 1$.
Then
 \[
  \ln\tsfrac{\psi(t^{1/d})}{\psi((t+c)^{1/d})}= \int_{t^{1/d}}^{(t+c)^{1/d}}\tsfrac{|\psi'(\tau)|}{\psi(\tau)}d\tau \le  K_6((t+c)^\tsfrac 1d-t^\tsfrac 1d)\le K_6c.
 \]
Taking into account this relation and monotonicity of $\psi$, we see that the first relation in \eqref{Aux_Eq01} is true.

For $\psi\in {\mathfrak M}^c_{\infty}$ the validity of the second relation in \eqref{Aux_Eq01}  is obvious.
For $\psi\in {\mathfrak M}'_{\infty}$ it follows from monotonicity of the functions $\alpha(\psi,t)=\tsfrac{\psi(t)}{t|\psi'(t)|}$ and $\tsfrac{\psi(t)}{|\psi'(t)|}$:
\[
1\ge \tsfrac{\psi(t^ {1/d})/|\psi'(t^ {1/d})|}
{\psi(ct^{1/d}/|\psi'(ct^{1/d})|}= \tsfrac{1}{c}\cdot
\tsfrac{\alpha(\psi,t^{1/d})}{\alpha(\psi,ct^{1/d})}\ge \tsfrac{1}{c}.
\]
\bend

Finally,  in the case when  $\psi\in {\mathfrak M}''_{\infty}$, for any positive $c$ we have
  \[
  \ln\tsfrac{\psi(t)}{\psi(t+c)}= \int_{t}^{t+c}\tsfrac{|\psi'(\tau)|}{\psi(\tau)}d\tau\ge \tsfrac{|\psi'(t)|}{\psi(t)}\uparrow \infty.
  \]

 {\bf \ref{Aux_St_01}.2.} Next, we formulate several known facts for functions from the
 sets ${\mathfrak M}'_{\infty}$, ${\mathfrak M}^c_{\infty}$ and ${\mathfrak M}''_{\infty}$.
 For this purpose,  following Stepanets \cite[Ch.~3]{Stepanets_2005}
 (see also \cite{Stepanets_1999}), for any $\psi\in {\mathfrak M}$, consider the functions $\eta(t)=\eta(\psi,t)$ and $\mu(t)=\mu(\psi,t)$ such that
  \begin{equation}\label{Aux_Eq02_01}
 \eta(t)=\psi^{-1}\big(\tsfrac{\psi(t)}{2}\big) \qquad \mbox{and}\qquad \mu(t)=\tsfrac {t}{\eta(t)-t},\qquad t\ge 1,
   \end{equation}
 where $\psi^{-1}$  is the inverse function of $\psi$, as well as the sets
 \[
 {\mathfrak M}_\infty^+=\{\psi\in {\mathfrak M}\,:\,\mu(\psi,t)\uparrow \infty\} \qquad \mbox{and}\qquad F=\{\psi\in {\mathfrak M}\,:\,\eta'(\psi,t)\le K\}.
\]

 \begin{remark}\label{Aux_Rem1}
  It follows from Theorems 12.1 and 13.1 \cite[Ch.~3]{Stepanets_2005}
  (see also \cite[Th.~1 and 2]{Stepanets_1999}) that all  sets ${\mathfrak M}'_\infty$,
  ${\mathfrak M}^c_\infty$ and ${\mathfrak M}''_\infty$ belong to the set ${\mathfrak M}_\infty^+\subset F$.
  By virtue of Remarks 13.1 and 13.2 from \cite[Ch.~3]{Stepanets_2005}
  (see also \cite[Remarks 1 and 2]{Stepanets_1999}), for any function $\psi\in F$
  (in particular, for any $\psi\in {\mathfrak M}'_\infty\cup{\mathfrak M}^c_\infty\cup{\mathfrak M}''_\infty$)
  and for all $t\ge 1$, the following relations hold:
  \begin{equation}\label{Aux_Eq02}
   K_7(\eta(\psi,t)-t)\le \psi(t)/|\psi'(t)|=t\alpha(\psi,t)\le K_8(\eta(\psi,t)-t)
   \end{equation}
   and
  \begin{equation}\label{Aux_Eq03}
   2(\eta(\psi,t)-t)\le \eta(\psi,\eta(\psi,t))-\eta(\psi,t)\le K_9(\eta(\psi,t)-t).
   \end{equation}
  \end{remark}

 {\bf \ref{Aux_St_01}.3.} Now we prove statements that give estimates for some integrals and sums  containing functions from the sets $B$, ${\mathfrak M}'_{\infty}$, ${\mathfrak M}^c_{\infty}$ or ${\mathfrak M}''_{\infty}$

 \begin{proposition}\label{Aux_Prop2}
Let $d\ge 1$ and $0<s<\infty$. Then for any   $\psi\in B$
\begin{equation}\label{Aux_Eq003_1}
    \int_1^{l}\tsfrac{t^{d-1}{\rm d}t}{\psi^s(t)}
     \asymp \sum\limits_{k=1}^{l}\tsfrac{k^{d-1}}{\psi^s(k)}\asymp      \tsfrac{l^{d}} {\psi^s(l)}.
  \end{equation}
If, in addition, for all $t$ greater than a certain number $t_0$, $\psi(t)$ is convex and satisfies  condition \eqref{Eq:OE_Sp02} with  $\beta>d/s$, then
\begin{equation}\label{Aux_Eq003_2}
     \int_l^\infty t^{d-1}\psi^s(t){\rm d}t
     \asymp \sum\limits_{k=l+1}^{\infty}{k^{d-1}}\psi^s(k)\asymp          l^{d}\psi^s(l).
  \end{equation}
  \end{proposition}

\banf Since the function $\tsfrac{t^{d-1} }{\psi^s(t)}$   is monotonically increasing,  we have
   \begin{equation}\label{Aux_Eq12}
     I_l:=\int_1^{l}\tsfrac{t^{d-1}{\rm d}t}{\psi^s(t)}\le \sum\limits_{n=1}^{l}\tsfrac{n^{d-1}}{\psi^s(n)}\le  I_{l+1},
     \end{equation}
and for any $\psi \in B$, the following necessary estimate is true:
 \[
 \tsfrac {l^{d}}{\psi^s(l)} \ll \tsfrac {(l/2)^{d}}{\psi^s(l/2)}  \ll
  \sum_{l/2\le n\le l} \,\tsfrac {n^{d-1}}{\psi^s(n)}\ll
  \sum_{k=1}^{l} \,\tsfrac {n^{d-1}}{\psi^s(n)}\ll
  \tsfrac {l^{d}}{\psi^s(l)}.
  \]

 Now,  assume that for all $t$ greater than a certain number $t_0$, the function $\psi$ is also convex and satisfies  condition \eqref{Eq:OE_Sp02} with  $\beta> \tsfrac ds$. Then the function $t^{d}\psi^s(t)$ decreases to zero for $t>t_0$. Therefore, for $l>t_0$, $l\in {\mathbb N}$,
   \begin{equation}\label{Aux_Eq1201}
  \int_l^\infty t^{d-1}\psi^s(t){\rm d}t=:J_{l} \ge \sum\limits_{n=l+1}^{\infty}{n^{d-1}}\psi^s(n)\ge J_{l+1}.
  \end{equation}
 Since  $\psi\in B$, then
    \begin{equation}\label{Aux_Eq1202}
   J_{l} \ge \int_l^{2l} t^{d-1}\psi^s(t){\rm d}t\gg l^{d}\psi^s(l).
  \end{equation}
 Applying \eqref{Eq:OE_Sp02} and integrating by parts, we obtain
 \[
J_{l}\le \tsfrac 1{\beta } \int_l^\infty t^{d}\psi^{s-1}(t)|\psi'(t)|dt
=\tsfrac {l^d\psi^s(l)}{ s \beta } +\tsfrac{d}{ s \beta }
J_{l}.
\]
Then in view of  \eqref{Eq:OE_Sp02}, we see that
    \begin{equation}\label{Aux_Eq1203}
J_{l}\le \tsfrac{1}{ s \beta -d}l^d\psi^s(l)  \ll l^d\psi^s(l)
      \end{equation}
Combining \eqref{Aux_Eq1202}, \eqref{Aux_Eq1203} and \eqref{Aux_Eq1201}, we obtain \eqref{Aux_Eq003_2}.

  \bend

\begin{proposition}\label{Aux_Prop3}
Let $d\ge 1$ and $0<s<\infty$. Then for any   $\psi\in {\mathfrak M}'_{\infty}\cup {\mathfrak M}^c_{\infty}$
\begin{equation}\label{Aux_Eq04}
    \int_1^{l}\tsfrac{t^{d-1}{\rm d}t}{\psi^s(t)}
     \asymp \sum\limits_{n=1}^{l}\tsfrac{n^{d-1}}{\psi^s(n)}\asymp
     \tsfrac{l^{d}\alpha(\psi,l)} {\psi^s(l)}
  \end{equation}
and
\begin{equation}\label{Aux_Eq05}
     \int_l^\infty t^{d-1}\psi^s(t){\rm d}t
     \asymp \sum\limits_{n=l+1}^{\infty}{n^{d-1}}\psi^s(n)\asymp          l^{d}\psi^s(l)\alpha(\psi,l).
  \end{equation}
  \end{proposition}

\banf
First let us prove relation \eqref{Aux_Eq04}. Consider the function
    \begin{equation}\label{Aux_Eq06}
    f(t)=f(\psi,t):= \tsfrac{\psi^s(t)}{t^{d-1}},\quad t\ge 1.
    \end{equation}
Its derivative   has the form
 \[
 f'(t)= \Big(\tsfrac{\psi^s(t)}{t^{d-1}}\Big)'=-\tsfrac{\psi^s(t)}{t^{d-1}}
  \Big(s\tsfrac{|\psi'(t)|}{\psi(t)}+\tsfrac{d-1}{t}\Big),\qquad t\ge 1.
\]
Therefore, for any  $\psi\in {\mathfrak M}'_{\infty}\cup{\mathfrak M}^c_{\infty}$, we have
\begin{equation}\label{Aux_Eq07}
  \tsfrac 1{t\alpha(f,t)}:= \tsfrac {|f'(t)|}{f(t)}=
   \tsfrac{|\psi'(t)|}{\psi(t)}\Big(s+(d-1)\alpha(\psi,t)\Big)\asymp
 \tsfrac{|\psi' (t)|}{\psi(t)}\asymp  \tsfrac 1{t\alpha(\psi,t)}
\end{equation}
and $\alpha(f,t)={\tsfrac{f(t)}{t|f'(t)|}}\downarrow 0$. Using \eqref{Aux_Eq06} and   \eqref{Aux_Eq12}, integrating by parts, we obtain
 \[
  I_l = \int_1^{l}\tsfrac{{\rm d}t}{f(t)}=\tsfrac {l}{f(l)}-\tsfrac 1{f(1)}-\!\int_1^ {l}
 \tsfrac {{\rm d}t}{\alpha(f,t)f(t)}\!\ge \tsfrac {l}{f(l)}-\tsfrac 1{f(1)}-\tsfrac {I_l}{\alpha(f,l)},
 \]
 and
 \begin{equation}\label{Aux_Eq08}
 I_l\ge \tsfrac {\alpha(f;l)}{1+\alpha(f;l)}
 \Big ( \tsfrac {l}{f(l)}-\tsfrac 1{f(1)}\Big )
 \gg \tsfrac{l\alpha(f,l)}{f(l)}.
\end{equation}

Further, if   $\psi\in {\mathfrak M}'_{\infty}$, then by  \eqref{Aux_Eq07},  we have $\tsfrac{f(t)} {|f'(t)|}\uparrow \infty$ and  ${f(t)}\ge \tsfrac{{f\,'}^2(t)}{f''(t)}$, $f''(t):=f''(t+)$,
 for almost all $t\ge 1$. Hence,
\begin{equation}\label{Aux_Eq09}
I_l\le \int_1^{l}
\tsfrac{f''(t){\rm d}t}{{f\,'}^2(t)}=\tsfrac 1{|f\,'(l)|}-\tsfrac
1{|f\,'(1)|}\le \tsfrac 1 {|f'(l)|}= \tsfrac{l\alpha(f,l)}{f(l)}.
 \end{equation}
Due to \eqref{Aux_Eq08} and \eqref{Aux_Eq09}, taking into account \eqref{Aux_Eq06} and \eqref{Aux_Eq07}, we see   for any $\psi\in {\mathfrak M}'_\infty$
\begin{equation}\label{Aux_Eq10}
     I_l= \int_1^{l}\tsfrac{t^{d-1}{\rm d}t}{\psi^s(t)}=\int_1^{l}\tsfrac{{\rm d}t}{f(t)}  \asymp \tsfrac{l\alpha(f,l)}{f(l)}\asymp
     \tsfrac{l^{d}\alpha(\psi,l)} {\psi^s(l)}.
     \end{equation}

 If  $\psi\in {\mathfrak M}_\infty^c$, by  \eqref{Aux_Eq07} and \eqref{Eq:OE_Sp11}, we have
 $ K_{10}\le \tsfrac{f(t)}{|f'(t)|}\le K_{11}$, $t\ge 1$, and
  \begin{equation}\label{Aux_Eq11}
 I_l \le K_{11}
 \int_1^{l} \tsfrac{|f'(t)|{\rm d}t}{f^2(t)}
    =K_{11}
 \Big(\tsfrac {1}{f(l)}-\tsfrac {1}{f(1)}\Big)\asymp\tsfrac {1}{f(l)}\asymp \tsfrac{l\alpha(f,l)}{f(l)}.
\end{equation}
Based on \eqref{Aux_Eq08}, \eqref{Aux_Eq11} \eqref{Aux_Eq06}, and \eqref{Aux_Eq07}, we see for $\psi\in {\mathfrak M}_\infty^c$,  relation  \eqref{Aux_Eq10} also holds.

  Combining   \eqref{Aux_Eq10} and \eqref{Aux_Eq12}, we get the second estimate  in  \eqref{Aux_Eq04}.

  Now let us prove relation \eqref{Aux_Eq05}. For this purpose, we  first show that the integral $J_l$ defined in \eqref{Aux_Eq1201} satisfies relation
     \begin{equation}\label{Aux_Eq13}
     J_l \asymp
     l^{d}\psi^s(l)\alpha(\psi,l).
     \end{equation}
Since $\alpha(\psi,t)\downarrow 0$ for    $\psi\in {\mathfrak M}'_{\infty}\cup {\mathfrak M}^c_{\infty}$, integrating by parts, we obtain
 \[
  J_l=-\tsfrac{l^d\psi^s(l)}{d}+\tsfrac {s}d\int_l^\infty \tsfrac{t^{d-1}\psi^s(t)}{\alpha(\psi,t)}{\rm d}t  \ge -\tsfrac{l^d\psi^s(l)}{d}+\tsfrac {s}{d\alpha(\psi,l)}J_l
  \]
  and for sufficiently large  $l$ (such that $\alpha(\psi,l)<\tsfrac sd$),
   \begin{equation}\label{Aux_Eq13_001}
J_l\le \tsfrac{\alpha(\psi,l)}{s-d\alpha(\psi,l)}l^d\psi^s(l)\ll l^d\psi^s(l)\alpha(\psi,l).
     \end{equation}
On the other hand-side, by virtue of monotonicity of $\psi$ and \eqref{Aux_Eq02}, we have
   \[
     J_l\ge \int_l^{\eta(\psi,t)} t^{d-1}\psi^s(t){\rm d}t\ge
     l^{d-1}\psi^s(\eta(\psi,l))(\eta(\psi,l)-l)\gg l^d\psi^s(l)\alpha(\psi,l).
    \]
Therefore, relation \eqref{Aux_Eq13} is true indeed.

Further, the derivative of the function  $h(t):=t^{d-1}\psi^s(t)$  has the form
 \[
 h'(t)= t^{d-2}\psi^s(t)
 \Big(d-1-s\tsfrac {t|\psi'(t)|}{\psi(t)}\Big) \qquad \forall t>1.
\]
Hence, based on \eqref{Eq:OE_Sp10}, we see that   the function $h(t)$ decreases for sufficiently large $t$ (such that $(d-1)\alpha(\psi,t)<s$), and therefore, relation \eqref{Aux_Eq1201} holds,
which, based on \eqref{Aux_Eq13} and \eqref{Aux_Eq01}, gives the necessary second estimate in \eqref{Aux_Eq05}.
  \bend


\begin{proposition}\label{Aux_Prop4}
Let $d\ge 1$ and $0<s<\infty$. Then for any   $\psi\in {\mathfrak M}''_{\infty}$
\begin{equation}\label{Aux_Eq14}
    \sum\limits_{n=1}^{l}\tsfrac{n^{d-1}}{\psi^s(n)}\asymp    \tsfrac{l^{d-1}} {\psi^s(l)}
  \end{equation}
and
\begin{equation}\label{Aux_Eq15}
       \sum\limits_{n=l+1}^{\infty}{n^{d-1}}\psi^s(n)\asymp          l^{d-1}\psi^s(l+1).
  \end{equation}
  \end{proposition}

  \banf
   For any $\psi\in {\mathfrak M}_\infty''$, we have $\psi(t)\le K_{12}|\psi'(t)|$, $t\ge 1$. Therefore,
 \[
 I_l=\int_1^l \tsfrac{t^{d-1}{\rm d}t}{\psi^s(t)}\le  K_{12} l^{d-1}\int_1^l \tsfrac{|\psi'(t)|{\rm d}t}{\psi^{s+1}(t)}
 =K_{12} l^{d-1} \Big(\tsfrac{1}{s\psi^s(l)}-
 \tsfrac{1}{s\psi^s(1)}\Big)\ll \tsfrac{l^{d-1}} {\psi^s(l)} .
 \]
 and estimate \eqref{Aux_Eq14} follows from the following relation:
 \[
    \tsfrac{l^{d-1}}{\psi^s(l)}\le \sum\limits_{k=1}^l\tsfrac{k^{d-1}}{\psi^s(k)}\le \tsfrac{l^{d-1}}{\psi^s(l)}+\int_1^l \tsfrac{t^{d-1}dt}{\psi^s(t)}\ll \tsfrac{l^{d-1}} {\psi^s(l)}.
 \]

  To prove the above estimate \eqref{Aux_Eq13_001}, we actually used only the fact that $\psi$ belongs to the set ${\mathfrak M}$ and satisfies the condition \eqref{Eq:OE_Sp10}.  Hence, estimate \eqref{Aux_Eq13_001} also holds for any $\psi\in {\mathfrak M}''_{\infty}$. However, since in this case $\frac{\psi(t)}{|\psi'(t)|}=t\alpha(\psi,t)\downarrow 0$, then
 \[
 J_l=\int_l^\infty t^{d-1}\psi^s(t){\rm d}t\ll l^d\psi^s(l)\alpha(\psi,l)\ll l^{d-1}\psi^s(l).
 \]
 As shown above, the function $h(t)=t^{d-1}\psi^s(t)$   decreases for all $t$ such that $\alpha(\psi,t)<\tsfrac{s}{d-1}$. Therefore,
 \[
 \sum\limits_{k=l+1}^{\infty}{k^{d-1}}\psi^s(k)\ll (l+1)^{d-1}\psi^s(l+1)+J_{l+1}\ll l^{d-1}\psi^s(l+1),
\]
and since
\[
 \sum\limits_{k=l+1}^{\infty}{k^{d-1}}\psi^s(k)\ge l^{d-1}\psi^s(l+1),
\]
relation \eqref{Aux_Eq15} is also true.
  \bend

\subsection{Proof of Theorem \ref{Th:_OE_Sp1}}

First, find the estimates for the quantities $\sigma_m({\mathcal F}_{q,r}^{\psi})_{_{\scriptstyle {\mathcal S}^p}}$.

\emph{Case $0<q\le p<\infty$.}  By virtue of \eqref{Eq:ExV_BnA01}, \eqref{Eq:ExV_BnA07} and \eqref{Eq:ExV_BnA11}, we have
 \begin{equation}\label{Eq:PrTh2_01}
  \sigma_m^p({\mathcal F}_{q,r}^{\psi})_{_{\scriptstyle {\mathcal S}^p}}
  =  \sup\limits_{l> m}(l-m)
  \Big(\sum_{n=1}^{n_l-1} \,\tsfrac {\nu_n}{\psi^q(n)}+
  \tsfrac{l-V_{n_l-1}}{\psi^q(n_l)}\Big)^{-\frac pq},
  \end{equation}
where  the number $n_l$ is defined by \eqref{Eq:ExV_BnA12}  for $s=l$.

By virtue of \eqref{Eq:ExV_BnA13}   and \eqref{Eq:OE_Sp01}, we have
   \begin{equation}\label{Eq:PrTh2_02}
   n_l\asymp l^\frac 1d\qquad \mbox{and}\qquad \psi(n_l) \asymp \psi(l^{\frac 1d}).
  \end{equation}
In view of  \eqref{Aux_Eq003_1}, we conclude that
\begin{equation}\label{Eq:PrTh2_03}
 \sigma_m^p({\mathcal F}_{q,r}^{\psi})_{_{\scriptstyle {\mathcal S}^p}}\asymp \sup\limits_{l> m}
  \tsfrac {l-m}{\big(\tsfrac {n_l^d}{\psi^q(n_l)}\big)^{\frac pq}}
  \asymp\sup\limits_{l> m} \tsfrac{\psi^p(l^{\frac 1d}) (l-m)}{l^{\frac pq}}\ll
  \psi^p(m^{\frac 1d}) \sup\limits_{l> m} \tsfrac{l-m}{l^{\frac pq}}.
\end{equation}

For $t>0$, $m\in {\mathbb N}$ and $s\in (1,\infty)$, the function $h(t)=h(t,s)=\tsfrac{t-m}{t^s}$ attains its maximal value at the point $t_*=\tsfrac{sm}{s-1}$, and
\begin{equation}\label{Eq:PrTh2_04}
 h(t_*,s)=\big(1- \tsfrac 1s\big)^s m^{1-s}.
\end{equation}
If $s=1$, then the function $h(t)=h(t;1)$ is non-decreasing and tends to $1$ as $t$ increases. Therefore,
 \begin{equation}\label{Eq:PrTh2_05}
 \sup\limits_{t>0} h(t;1)= \sup\limits_{t>0}
 \tsfrac {t-m}{t}=\lim\limits_{t\to+\infty} \tsfrac{t-m}{t}=1.
 \end{equation}
Combining  \eqref{Eq:PrTh2_03}--\eqref{Eq:PrTh2_05}, we obtain necessary upper estimate:
 \[
  \sigma_m^p({\mathcal F}_{q,r}^{\psi})_{_{\scriptstyle {\mathcal S}^p}}\ll
  {\psi^p(m^{\frac 1d})}{m^{1-\frac pq}}.
 \]
Taking into account  \eqref{Eq:PrTh2_03} and  \eqref{Eq:OE_Sp01}, we also obtain the lower estimate:
 \[
 \sigma_m^p({\mathcal F}_{q,r}^{\psi})_{_{\scriptstyle {\mathcal S}^p}}\gg
 \tsfrac{\psi^p ((2m)^{\frac 1d}) {(2m-m)}}{(2m)^{\frac pq}}\asymp
 {\psi^p(m^{\frac 1d})}{m^{1-\frac pq}}.
 \]


 \emph{Case $0<p<q<\infty$.}   Let $\Psi=\{\Psi_k\}_{k \in {\mathbb Z}^d}$ be a sequence of complex numbers
 such that  there exists a non-increasing rearrangement $\bar\Psi= \{\bar \Psi_j\}_{j=1}^\infty$ of the
 number system   $\{|\Psi_k|\}_{k \in {\mathbb Z}^d}$. For a fixed $m\in {\mathbb N}$ and    $l\in {\mathbb N}$, $l>m$,  consider the functional
  \[
 Q_m(\Psi,l):=\tsfrac {l-m}{\sum\limits_{j=1}^l \bar\Psi_j^{-q}}.
 \]
 We have
 $$
 Q_m(\Psi,l+1)=Q_m(\Psi,l)+ \tsfrac {\bar\Psi_{l+1}^{-q}}
 {\sum\limits_{j=1}^{l+1} \bar\Psi_j^{-q}}\Big(\bar\Psi_{l+1}^{q}-Q_m(\Psi,l)\Big)
 $$
and
$$
 \bar\Psi_{l+1}^{q}=Q_m(\Psi,l+1)+ \tsfrac {\sum\limits_{j=1}^{l} \bar\Psi_j^{-q}}
 {\sum\limits_{j=1}^{l+1} \bar\Psi_j^{-q}}\Big(\bar\Psi_{l+1}^{q}-Q_m(\Psi,l)\Big)
 $$
Therefore, in view of monotonicity of  $\bar\Psi$ and  relation \eqref{Eq:ExV_BnA03}, we conclude that
 \[
 Q_m(\Psi,l)>Q_m(\Psi,l+1)>\Psi^q(l+1)\qquad  \forall l\ge l_m
 \]
 and
 \[
 Q_m(\Psi,l)\le Q_m(\Psi,l+1)\le \Psi^q(l+1) \qquad  \forall l\in [m,l_m).
 \]
 This yields that
\begin{equation}\label{Eq:PrTh2_06}
\sup\limits_{l>m}  Q_{n}(\Psi,l)=Q_m(\Psi,l_m).
 \end{equation}
According to \eqref{Eq:ExV_BnA03}, we get $\Psi(l_m+1)>\Psi(l_m)$. Hence, if the  system
  $\{|\Psi_k|\}_{k \in {\mathbb Z}^d}$  has the   stepwise form \eqref{Eq:ExV_BnA07},  then
\begin{equation}\label{Eq:PrTh2_07}
l_m=V_{n_{l_m}}=\sum_{i=0}^{n_{l_m}}\nu_i,
\end{equation}
where $n_{l_m}$ is defined in \eqref{Eq:ExV_BnA12} for $s=l_m$.  Thus,   we have
 \begin{equation}\label{Eq:PrTh2_08}
  \sigma_m({\mathcal F}_{q,r}^{\psi})_{_{\scriptstyle {\mathcal S}^p}}= \bigg(  (l_m-m)^{\frac{q}{q-p}}\Big(\sum _{n=1}^{n_{l_m}}\tsfrac{\nu_n}{\psi^q(n)}\Big)^{\frac{p}{p-q}}
  +    \sum\limits_{n=n_{l_m}+1}^\infty
{\nu_n}{\psi^\frac{pq}{q-p}(n)} \bigg)^{\frac{q-p}{pq}}  ,
  \end{equation}
 where
\begin{equation}\label{Eq:PrTh2_09}
 \psi^{-q}({n_{l_m}})\le \tsfrac 1{l_m-m}\sum _{j=1}^{n_{l_m}}\tsfrac{\nu_n}{\psi^q(n)}<\psi^{-q}({n_{l_m}+1}).
\end{equation}

It follows from \eqref{Eq:PrTh2_06} that
 \begin{equation}\label{Eq:PrTh2_10}
 \sup\limits_{l> m}(l-m) \bigg(\sum_{n=1}^{n_l-1} \,\tsfrac {\nu_n}{\psi^q(n)}+ \tsfrac{l-V_{n_l-1}}{\psi^q(n_l)}\bigg)^{-1}=(l_m-m)\bigg(\sum _{n=1}^{n_{l_m}}\tsfrac{\nu_n}{\psi^q(n)}\bigg)^{-1}.
\end{equation}
Then similarly to the case $0<p\le q<\infty$, we show that
\begin{equation}\label{Eq:PrTh2_11}
(l_m-m)\bigg(\sum _{n=1}^{n_{l_m}}\tsfrac{\nu_n}{\psi^q(n)}\bigg)^{-1}\asymp \psi^q(m^{\frac 1d}).
\end{equation}
Taking into account \eqref{Eq:PrTh2_09}, \eqref{Eq:PrTh2_11} and \eqref{Eq:OE_Sp01}, we see that
\begin{equation}\label{Eq:PrTh2_12}
\psi(n_{l_m})\asymp \psi(m^{\frac 1d}).
\end{equation}
Since  $\psi\in B$, then in view of \eqref{Eq:ExV_BnA11}, \eqref{Eq:OE_Sp02} and \eqref{Aux_Eq003_2}, we conclude that
\begin{equation}\label{Eq:PrTh2_13}
\sum\limits_{n=l+1}^\infty{\nu_n}\psi^\frac{pq}{q-p}(n)\asymp \sum\limits_{n=l+1}^{\infty}{n^{d-1}}\psi^\frac{pq}{q-p}(n)\asymp l^d\psi^\frac{pq}{q-p}(l).
\end{equation}

Let us also show that
\begin{equation}\label{Eq:PrTh2_14}
 n_{l_m}\asymp m^\frac 1d.
\end{equation}
Indeed, by virtue of \eqref{Eq:PrTh2_07} and \eqref{Eq:PrTh2_02}, we  see that $
   \widetilde{m}:=K_{13} m^\frac 1d\le K_{13} l_m^\frac 1d\le n_{l_m}$.
 On the other hand, integrating each part of \eqref{Eq:OE_Sp06} in the range from
 $\widetilde{m}$ to $n_{l_m}$,  $\widetilde{m}>t_0$,   we obtain
$\tsfrac{\psi(\widetilde{m})}{\psi(n_{l_m})} \ge (\tsfrac{n_{l_m}}{\widetilde{m}})^{\beta }$.
Therefore, in view of \eqref{Eq:PrTh2_12} and \eqref{Eq:OE_Sp01}, we see that $\widetilde{m}\gg n_{l_m}$ and relation \eqref{Eq:PrTh2_14} is true.


Combining  \eqref{Eq:PrTh2_08}, \eqref{Eq:PrTh2_11},
\eqref{Aux_Eq003_1},  \eqref{Eq:PrTh2_12},  \eqref{Eq:PrTh2_13}
and \eqref{Eq:PrTh2_14} we obtain  \eqref{Eq:OE_Sp04}:
 \[
 \sigma_m({\mathcal F}_{q,r}^{\psi})_{_{\scriptstyle {\mathcal S}^p}}
 \asymp
 \Big(\psi^{\frac {q^2}{q-p}}(m^{\frac 1d}) \tsfrac {m}{\psi^q(m^{\frac 1d})}
  +     m\psi^\frac{pq}{q-p}(m^{\frac 1d}) \Big)^{\frac{q-p}{pq}} \asymp \psi(m^{\frac 1d}) m^{\frac 1p-\frac 1q}.
 \]


\emph{Case  $0<q<p=\infty$.}  By virtue of \eqref{Eq:ExV_BnA04}, \eqref{Eq:ExV_BnA07} and \eqref{Eq:ExV_BnA11}, we have
  \[
 \sigma_m({\mathcal F}_{q,r}^{\psi})_{_{\scriptstyle {\mathcal S}^p}}=
\Big(\sum_{n=1}^{n_{m+1}-1} \,\tsfrac {\nu_n}{\psi^q(n)}+
  \tsfrac{m+1-V_{n_{m+1}-1}}{\psi^q(n_{m+1})}\Big )^{-\frac
  1q}.
 \]
 Taking into account \eqref{Eq:ExV_BnA11},  \eqref{Aux_Eq003_1} and   \eqref{Eq:PrTh2_02}, we see that
  \[
 \sigma_m({\mathcal F}_{q,r}^{\psi})_{_{\scriptstyle {\mathcal S}^p}}\asymp
\Big (\sum_{n=1}^{n_{m+1}} \,\tsfrac {\nu_n}{\psi^q(n)} \Big )^{-\frac 1q}\asymp \psi(m^{\frac 1d})m^{-\frac 1q}.
 \]


\emph{Case  $0<p<q=\infty$.} Similarly to the proof of \eqref{Eq:OE_Sp08} and \eqref{Eq:PrTh2_13}, we can show that  condition \eqref{Eq:OE_Sp02} with $\beta=d(\tsfrac 1p-\frac 1q)=\tsfrac dp$ guarantees the convergence of the series
 \[
 \sum_{k\in {\mathbb Z}^d} |\Psi_k|^p=\sum\limits_{j=1}^{\infty} \bar{\Psi}^{p}_j=\sum\limits_{n=1}^\infty{\nu_n}\psi^p(n)
 \]
 and
 \begin{equation}\label{Eq:PrTh2_15}
  \sum\limits_{n=l+1}^\infty{\nu_n}\psi^p(n)\asymp l^d\psi^p(l).
\end{equation}
Then using  \eqref{Eq:ExV_BnA05}, \eqref{Eq:ExV_BnA07}, \eqref{Eq:ExV_BnA11}, \eqref{Eq:PrTh2_15} and   \eqref{Eq:PrTh2_02}, we obtain \eqref{Eq:OE_Sp04}:
  \[
  \sigma_m({\mathcal F}_{q,r}^{\psi})_{_{\scriptstyle {\mathcal S}^p}}=
 \Big((V_{n_{m+1}}-m-1){\psi^p(n_{m+1})}+\!\!\!\sum\limits_{n=n_{m+1}+1}^\infty
{\nu_n}{\psi^p(n)} \Big)^{\frac 1p} \asymp \psi(m^{\frac 1d})m^{\frac 1p}.
 \]

\emph{Case  $p=q=\infty$.} In this case,  estimate \eqref{Eq:OE_Sp04} follows from \eqref{Eq:ExV_BnA06}, \eqref{Eq:ExV_BnA07} and \eqref{Eq:PrTh2_02}
 \[
  \sigma_m({\mathcal F}_{q,r}^{\psi})_{_{\scriptstyle {\mathcal S}^p}}=
  \bar{\Psi}_{m+1}=\psi(n_{m+1})\asymp  \psi(m^{\frac 1d}).\]

Finally, consider the quantities ${\mathscr D}_m({\cal F}_{q,r}^{\psi})_{_{\scriptstyle {\mathcal S}^p }} $.
In the case $0<q\le p\le \infty$,   estimate \eqref{Eq:OE_Sp03} follows from \eqref{Eq:ExV_BaW05}, \eqref{Eq:ExV_BnA07} and  \eqref{Eq:PrTh2_02}:
  \[
 {\mathscr D}_m({\mathcal F}_{q,r}^{\psi})_{_{\scriptstyle {\mathcal S}^p}} =\psi(n_{m+1})\asymp  \psi(m^{\frac 1d}).
 \]
For  $0<p<q<\infty$, estimate  \eqref{Eq:OE_Sp03}  follows from \eqref{Eq:ExV_BaW07},  \eqref{Eq:ExV_BnA07}, \eqref{Eq:PrTh2_13} and  \eqref{Eq:PrTh2_02}:
\[
 {\mathscr D}_m({\mathcal F}_{q,r}^{\psi})_{_{\scriptstyle {\mathcal S}^p}}
    \asymp\Big(\sum\limits_{n=n_{m+1}+1}^\infty
{\nu_n}{\psi^{\frac{p\,q}{q-p}}(n)} \Big)^{\frac{q-p}{p\,q}}\asymp  \psi(n_{m+1})n_{m+1}^{\frac{q-p}{d\,p\,q}}
\asymp
\psi(m^{\frac 1d})m^{\frac 1p-\frac 1q},
\]
and for  $0<p<q=\infty$, it similarly    follows from  \eqref{Eq:ExV_BaW11}, \eqref{Eq:ExV_BnA07}, \eqref{Eq:PrTh2_15}  and  \eqref{Eq:PrTh2_02}.

 \bend

\subsection{Proof of Theorem \ref{Th:_OE_Sp2}}\label{Pr_Th_OE_Sp2}

First, consider the quantities ${\mathscr D}_m({\cal F}_{q,r}^{\psi})_{_{\scriptstyle {\mathcal S}^p }}$.
 In the case   $0<q\le p\le \infty$, estimate \eqref{Eq:OE_Sp13} follows from the corresponding results of Theorems A and \ref{Th:ApChFS01}, relations \eqref{Eq:ExV_BnA12}, \eqref{Eq:ExV_BnA13} and \eqref{Eq:OE_Sp12}, taking into account Proposition \ref{Aux_Prop1}:
     \begin{equation}\label{Eq:PrTh3_01}
     {\mathscr D}_m({\cal F}_{q,r}^{\psi})_{_{\scriptstyle {\mathcal S}^p }}=\psi(n_{m+1})\asymp \psi(\widetilde{n}_m).
     \end{equation}
 In the case  $0<p<q\le \infty$, estimate \eqref{Eq:OE_Sp13} similarly follows from the results
 of Theorems A and \ref{Th:ApChFS01}, relations \eqref{Eq:ExV_BnA11}, \eqref{Aux_Eq05}, \eqref{Eq:ExV_BnA13} and \eqref{Eq:OE_Sp12}, taking into account Proposition \ref{Aux_Prop1}:
 \[
  {\mathscr D}_m({\cal F}_{q,r}^{\psi})_{_{\scriptstyle {\mathcal S}^p }}\asymp\Big(\sum\limits_{n=n_{m+1}+1}^\infty
{\nu_n}{\psi^{\frac{p\,q}{q-p}}(n)} \Big)^{\frac{q-p}{p\,q}}\asymp
\Big(\sum\limits_{n=n_{m+1}+1}^\infty
{n^{d-1}}{\psi^{\frac{p\,q}{q-p}}(n)} \Big)^{\frac{q-p}{p\,q}}
 \]
 \begin{equation}\label{Eq:PrTh3_02}
 \asymp
 \Big(n_{m+1}^{d}\psi^{\frac{p\,q}{q-p}}(n_{m+1})\alpha(\psi,n_{m+1})\Big)^{\frac{q-p}{p\,q}}\asymp
 \psi(\widetilde{n}_m) (m \alpha(\psi,\widetilde{n}_m))^{\frac 1p-\frac 1q}.
    \end{equation}

Now, find the estimates for the quantities $\sigma_m({\mathcal F}_{q,r}^{\psi})_{_{\scriptstyle {\mathcal S}^p}}$. Note that in the \emph{cases  $p=q=\infty$ and $0<p<q=\infty$,} estimate \eqref{Eq:OE_Sp14} is obtained similarly to the
estimates \eqref{Eq:PrTh3_01} and \eqref{Eq:PrTh3_02}, using the corresponding results of Theorem B.

In the \emph{case  $0<q<p=\infty$,} it  follows from relations \eqref{Eq:ExV_BnA04}, \eqref{Eq:ExV_BnA12},
\eqref{Eq:ExV_BnA13}, \eqref{Eq:OE_Sp12} by using  Propositions \ref{Aux_Prop3} and \ref{Aux_Prop1}:
 \[
 \sigma_m({\mathcal F}_{q}^{\Psi})_{_{\scriptstyle {\mathcal S}^p}}\asymp
 \Big(\sum\limits_{n=1}^{n_{m+1}} \tsfrac {\nu_n}{\psi^q(n)}\Big)^{-\frac 1q}
 \asymp
 \Big(\sum\limits_{n=1}^{n_{m+1}} \tsfrac {n^{d-1}}{\psi^q(n)}\Big)^{-\frac 1q}
 \]
 \[
  \asymp\Big(\tsfrac{n_{m+1}^{d}\alpha(\psi,n_{m+1})} {\psi^q(n_{m+1})}\Big)^{-\frac 1q}
  \asymp \psi(\widetilde{n}_m) \big(m\alpha(\psi,\widetilde{n}_m)\big)^{-\frac 1q}.
 \]

\emph{Case $0<q\le p<\infty$.} Let  $M_{r,d}$ be the number defined by  \eqref{Eq:ExV_BnA09} and $c$ be a fixed real number. For $m\in {\mathbb N}$ and any $l\ge m$, consider the function
 \[
 W_m(l,c):= \tsfrac{l-m}{\big(\int\limits_{1}^{k(l)}
  \tsfrac{t^{d-1}{\rm d}t}{\psi^q(t)}\big)^{\frac pq}},\qquad \mbox{where}\ \ k(l)=k(l,c)=(l/M_{r,d})^\frac 1d+c.
 \]
 This function is continuously differentiable,  its derivative   has the form
 \begin{equation}\label{Eq:PrTh3_03}
 W'_m(l,c)=\Big (\int_{1}^{k(l)}
  \tsfrac{t^{d-1}{\rm d}t}{\psi^q(t)}-\tsfrac
  {p (l-m)}{q\psi^q(k(l))}\cdot \tsfrac{k^{d-1}(l)}{dM_{r,d}^{\frac 1d}l^{(d-1)/d}}\Big )\Big ( \int_{1}^{k(l)}
  \tsfrac{t^{d-1}{\rm d}t}{\psi^q(t)}\Big )^{\frac pq-1},
 \end{equation}
where
  \begin{equation}\label{Eq:PrTh3_04}
  \lim\limits_{l\to\infty} \tsfrac{k^{d-1}(l)}{l^{(d-1)/d}}=
  \lim\limits_{l\to\infty} \Big (\tsfrac{(l/M_{r,d})^{\frac 1d}+c}{l^{\frac 1d}}\Big )^{d-1}
  =M_{r,d}^{\frac{1-d}{d}}.
 \end{equation}
 There exists a unique point $l_m>m$ such that $W'_m(l_m,c)=0$. The function $W_m(l,c)$ increase for $m\le l<l_m$ and decrease for $l>l_m$.  Moreover, it follows from   \eqref{Eq:PrTh3_03}, \eqref{Eq:PrTh3_04}  and \eqref{Aux_Eq10} that
   \begin{equation}\label{Eq:PrTh3_05}
   (l_m-m)\Big (\int_{1}^{k(l_m)}
  \tsfrac{t^{d-1}{\rm d}t}{\psi^q(t)}\Big )^{-1}=\tsfrac{dq M_{r,d}^{\frac 1d}l^{(d-1)/d}}
    {k^{d-1}(l)}  \psi^q(k(l_m))  \asymp \psi^q(k(l_m)).
 \end{equation}

Due to \eqref{Eq:PrTh3_05} and \eqref{Aux_Eq04}, for any  $\psi\in {\mathfrak M}'_{\infty}\cup {\mathfrak M}^c_{\infty}$, we obtain
 $$
 \sup\limits_{l\ge m}W_m(l,c)=W_m(l_m,c)
=\tsfrac{l_m-n}{\int\limits_{1}^{k(l_m)}
  \tsfrac{t^{d-1}{\rm d}t}{\psi^r(t)}}\cdot
  \Big (\int_{1}^{k(l_m)}
  \tsfrac{t^{d-1}{\rm d}t}{\psi^q(t)}\Big )^{1- \frac pq}
 \asymp
 $$
 \begin{equation}\label{Eq:PrTh3_06}
  \asymp  \psi^q(k(l_m))
  \Big ( \tsfrac{k^{d}(l_m)\alpha(\psi,k(l_m))}
 {\psi^q(k(l_m))}\Big )^{   1-\frac pq}     = \tsfrac{\psi^p(k(l_m))}{ \big (k^{d}(l_m)\alpha(\psi,k(l_m))\big )^{\frac pq-1}} .
 \end{equation}
 Since $l_m\ge m$ and $k(l_m)\ge (m/M_{r,d})^\frac 1d-|c|=\widetilde{n}_m-|c|$, for any  $\psi\in {\mathfrak M}'_{\infty}\cup {\mathfrak M}^c_{\infty}$, taking into account Proposition \ref{Aux_Prop1},  we have
  \[
 \tsfrac{\psi^p(k(l_m))}{ \big (k^{d}(l_m)\alpha(\psi,k(l_m))\big )^{\frac pq-1}}\!\le
 \tsfrac{\psi^p(\widetilde{n}_m-|c|)}{ \big ((\widetilde{n}_m-|c|)^{d-1} \frac{\psi(k(l_m))}{|\psi'(k(l_m))|}\big )^{\frac pq-1}}\ll
  \tsfrac{\psi^p(\widetilde{n}_m)}{ \big (\widetilde{n}_m^{d-1} \frac{\psi(\widetilde{n}_m)}{|\psi'(\widetilde{n}_m)|}\big )^{\frac pq-1}}=
  \tsfrac{\psi^p(\widetilde{n}_m)}{ \big (m\alpha(\psi,\widetilde{n}_m)\big )^{\frac pq-1}}.
 \]
 Combining  the last estimate and \eqref{Eq:PrTh3_06}, we conclude that for any $\psi\in {\mathfrak M}'_{\infty}\cup {\mathfrak M}^c_{\infty}$
  \begin{equation}\label{Eq:PrTh3_07}
 \sup\limits_{l\ge m}W_m(l,c)\ll  \tsfrac{\psi^p(\widetilde{n}_m)}
 { \big (m\alpha(\psi,\widetilde{n}_m)\big )^{\frac pq-1}}.
 \end{equation}
By virtue of \eqref{Eq:PrTh2_01}, \eqref{Eq:ExV_BnA13}  and \eqref{Eq:PrTh3_07}, we obtain the upper estimate in \eqref{Eq:OE_Sp14}:
 \[
  \sigma_m^p({\mathcal F}_{q,r}^{\psi})_{_{\scriptstyle {\mathcal S}^p}}
  \le \sup\limits_{l> m} \tsfrac{l-m}{\big(\sum\limits_{n=1}^{n_l-1} \,\tsfrac {\nu_n}{\psi^q(n)} \big)
  ^\frac pq}
  \ll
  \sup\limits_{l> m} \tsfrac{l-m}
  {\big(\sum\limits_{n=1}^{n_l-1} \,\tsfrac {n^{d-1}}{\psi^q(n)} \big)^\frac pq}
  \]
   \[
   \ll \sup\limits_{l> m} \tsfrac{l-m}
  {\big(\int\limits_1^{n_l-1}\tsfrac{t^{d-1}{\rm d}t}{\psi^q(t)} \big)^\frac pq}
  \ll \sup\limits_{l> m} W_m(l,c_{r,d}+1)\ll
  \tsfrac{\psi^p(\widetilde{n}_m)}{ \big (m\alpha(\psi,\widetilde{n}_m)\big )^{\frac pq-1}}.
  \]

Let us also find the lower estimate. Due to \eqref{Eq:PrTh2_01}, \eqref{Eq:ExV_BnA13}, \eqref{Aux_Eq04} and \eqref{Aux_Eq01}, we have
 \[
  \sigma_m^p({\mathcal F}_{q,r}^{\psi})_{_{\scriptstyle {\mathcal S}^p}}
  \gg \sup\limits_{l> m} \tsfrac{l-m}{\big(\sum\limits_{n=1}^{n_l}  \tsfrac {n^{d-1}}{\psi^q(n)} \big)^\frac pq}\gg
  \sup\limits_{l> m, l\in {\mathbb N}} \tsfrac{\psi^p(n_l)(l-m)}{\big(n_l^{d}\alpha(\psi,n_l) \big)^\frac pq}
  \]
\begin{equation}\label{Eq:PrTh3_08}
   \gg
  \sup\limits_{l> m, l\in {\mathbb N}} \tsfrac{\psi^p((l/M_{r,d})^{\frac 1d} )(l-m)}{\big(l\alpha(\psi,(l/M_{r,d})^{\frac 1d}) \big)^\frac pq}=:\sup\limits_{l> m, l\in {\mathbb N}} R_m(l).
\end{equation}
Consider the function
\begin{equation}\label{Eq:PrTh3_09}
g(t)=g(\psi;t):=\psi((t/M_{r,d})^{\frac 1d}),\quad t\ge 1.
\end{equation}
 It is easy to see that $g\in {\mathfrak M}$ and
  \[
  \alpha(g,t)=\tsfrac {g(t)}{t|g'(t)|}=
  dM_{r,d}^{\frac 1d}\,\alpha(\psi,(t/M_{r,d})^{\frac 1d})\downarrow 0
  \]
   for any   $\psi\in {\mathfrak M}'_{\infty}\cup{\mathfrak M}_\infty^c$. Therefore, similar to $\psi$, the function $g\in {\mathfrak M}_\infty^+\subset F$. Hence, by virtue of \eqref{Aux_Eq02} and \eqref{Aux_Eq03}, for the quantity $\eta(t)=\eta(g,t)$, we have
\begin{equation}\label{Eq:PrTh3_10}
 \eta(\eta(t))-\eta(t)\asymp   \eta(t)-t\asymp \tsfrac {g(t)}{|g'(t)|} \asymp
t \alpha(\psi,(t/M_{r,d})^{\frac 1d}).
\end{equation}
Set $l_0=[\eta(m)]+1$. Then by Proposition \ref{Aux_Prop1}, Remark \ref{Aux_Rem1}, relations
\eqref{Eq:PrTh3_09}, \eqref{Eq:PrTh3_10} and \eqref{Eq:OE_Sp12},
 \begin{equation}\label{Eq:PrTh3_11}
 R_m(l_0)\gg \tsfrac{\psi^p((\eta(g,m)/M_{r,d})^{\frac 1d})
 (\eta(g,m)-m)}{( \eta(g,m)\alpha(\psi,(\eta(g,m)/M_{r,d})^{\frac 1d})
 )^\frac pq} \gg  \tsfrac{g^p(m)}{( \eta(g,m)-m)
 )^{ \frac pq-1}}
 \gg
 \tsfrac{\psi^p(\widetilde{n}_m)}
  {(m\alpha(\psi,\widetilde{n}_m)\,
  )^{\frac pq-1}}.
\end{equation}
Combining relations \eqref{Eq:PrTh3_08} and \eqref{Eq:PrTh3_11}, we obtain the necessary lower bound:
 \[
 \sigma_m^p({\mathcal F}_{q,r}^{\psi})_{_{\scriptstyle {\mathcal S}^p}}
 \gg
 \tsfrac{\psi^p(\widetilde{n}_m)}
  {(m\alpha(\psi,\widetilde{n}_m)\,
  )^{\frac pq-1}}
\]
which complies the proof of relation \eqref{Eq:OE_Sp14} in the case  $0<q\le p<\infty$.

\emph{Case $0<p<q<\infty$.} Let us use relation \eqref{Eq:PrTh2_08}--\eqref{Eq:PrTh2_10}, which holds for any positive non-increasing function $\psi$. Similarly to the previous cases, we prove that
  \begin{equation}\label{Eq:PrTh3_12}
 \sup\limits_{l> m}\tsfrac{l-m}{\sum\limits_{n=1}^{n_l-1} \,\tsfrac {\nu_n}{\psi^q(n)}+ \tsfrac{l-V_{n_l-1}}{\psi^q(n_l)}}=\tsfrac{l_m-m}{\sum\limits_{n=1}^{n_{l_m}} \,\tsfrac {\nu_n}{\psi^q(n)}}\asymp \psi^q(\widetilde{n}_m)\asymp \psi^q(n_{l_m}),
\end{equation}
where $\widetilde{n}_m$ and $n_{l_m}$ are defined by \eqref{Eq:OE_Sp12} and \eqref{Eq:PrTh2_09}.
 Further, let us show that
  \begin{equation}\label{Eq:PrTh3_13}
 l_m\asymp \widetilde{n}_m.
 \end{equation}
 Since $l_m>m$,   by \eqref{Eq:ExV_BnA13}  we have $n_{l_m}> n_m \ge \widetilde{n}_m-c_2$.
 If $n_{l_m}\le \widetilde{n}_m$, then relation \eqref{Eq:PrTh3_13} holds.
 Assume that $n_{l_m}> \widetilde{n}_m$.  Since $\alpha(\psi,t)\downarrow 0$,  then
 $\tsfrac 1t\le K_{14} \tsfrac{|\psi'(t)|}{\psi(t)}$  for all $t\ge 1$.
 Integrating both parts of this inequality in the range from $\widetilde{n}_m$ to $n_{l_m}$, we obtain
  \[
   \ln \tsfrac {n_{l_m}}{\widetilde{n}_m}\le K_{14} \ln \tsfrac{\psi(\widetilde{n}_m)}{\psi(n_{l_m})}\qquad \mbox{and}\qquad
   n_{l_m}\ll \widetilde{n}_m.
  \]
 Thus, relation \eqref{Eq:PrTh3_13} is true indeed.

 By virtue of \eqref{Eq:ExV_BnA10}, \eqref{Aux_Eq04}, \eqref{Eq:PrTh3_13}  and \eqref{Aux_Eq01},
  \begin{equation}\label{Eq:PrTh3_14}
 \sum _{n=1}^{n_{l_m}}\tsfrac{\nu_n}{\psi^q(n)}\asymp  \sum _{n=1}^{n_{l_m}}\tsfrac{n^{d-1}}{\psi^q(n)}\asymp \tsfrac{n_{l_m}^{d}\alpha(\psi,n_{l_m})} {\psi^q(n_{l_m})}\asymp \tsfrac{m\alpha(\psi,\widetilde{n}_m)} {\psi^q(\widetilde{n}_m)}.
 \end{equation}
 Finally, due to \eqref{Eq:ExV_BnA10}, \eqref{Aux_Eq05},  \eqref{Eq:PrTh3_13}
 and \eqref{Aux_Eq01}, we get
\[
 \sum\limits_{n=n_{l_m}+1}^\infty
{\nu_n}{\psi^\frac{pq}{q-p}(n)}\asymp \sum\limits_{n=n_{l_m}+1}^\infty
{n^{d-1}}{\psi^\frac{pq}{q-p}(n)}
\]
   \begin{equation}\label{Eq:PrTh3_15}
   \asymp  n_{l_m}^{d}\psi^\frac{pq}{q-p}(n_{l_m})\alpha(\psi,n_{l_m})
   \asymp  m\psi^\frac{pq}{q-p}(\widetilde{n}_m)\alpha(\psi,\widetilde{n}_m).
  \end{equation}

 Combining relations \eqref{Eq:PrTh2_08}, \eqref{Eq:PrTh3_12},  \eqref{Eq:PrTh3_14} and
 \eqref{Eq:PrTh3_15}, we obtain \eqref{Eq:OE_Sp14}:
  \[
  \sigma_m({\mathcal F}_{q,r}^{\psi})_{_{\scriptstyle {\mathcal S}^p}}\asymp
  \Big(\psi^{\frac {q^2}{q-p}}(\widetilde{n}_m)  \tsfrac{m\alpha(\psi,\widetilde{n}_m)} {\psi^q(\widetilde{n}_m)}
  +    m\psi^\frac{pq}{q-p}(\widetilde{n}_m)\alpha(\psi,\widetilde{n}_m) \Big)^{\frac{q-p}{pq}}
  \asymp\psi(\widetilde{n}_m) \big(m \alpha(\psi,\widetilde{n}_m)\big)^{\frac 1p-\frac 1q}.
   \]
  \bend

\subsection{Proof of Theorem \ref{Th:_OE_Sp3}}\label{Pr_Th:_OE_Sp4}

First, we find the estimates for the quantities $\sigma_m({\mathcal F}_{q,r}^{\psi})_{_{\scriptstyle {\mathcal S}^p}}$.

 \emph{{\rm (i)} Case $0<p<q<\infty$.}  Let us use the representation \eqref{Eq:PrTh2_08},
 in which  $n_{l_m}$ is defined in \eqref{Eq:ExV_BnA12} for $s=l_m$, and $l_m$ satisfies
 relations \eqref{Eq:PrTh2_09} and \eqref{Eq:PrTh2_07}.
 By virtue of \eqref{Eq:ExV_BnA11} and \eqref{Aux_Eq15},
  \begin{equation}\label{Eq:PrTh4_01}
 \sum\limits_{n=n_{l_m}+1}^\infty
{\nu_n}{\psi^\frac{pq}{q-p}(n)}\asymp \sum\limits_{n=n_{l_m}+1}^\infty
{n^{d-1}}{\psi^\frac{pq}{q-p}(n)}
\asymp
n_{l_m}^{d-1}{\psi^\frac{pq}{q-p}(n_{l_m}+1)}.
 \end{equation}
Since $l_m\ge m+1$ and $m\in   [V_{s-1},V_s)$, it follows from \eqref{Eq:PrTh2_07} and \eqref{Eq:ExV_BnA12} that $l_m=V_{n_{l_m}}\ge V_{n_{m+1}}$.

 By virtue of \eqref{Aux_Eq14}, \eqref{Eq:ExV_BnA10} and \eqref{Eq:OE_Sp16}, for sufficiently large $s$ we have
 \begin{equation}\label{Eq:PrTh4_02}
\sum _{j=1}^{s}\tsfrac{\nu_j}{\psi^q(j)}\asymp \sum _{j=1}^{s}\tsfrac{j^{d-1}}{\psi^q(j)}\asymp \tsfrac{s^{d-1}}{\psi^q(s)}<\psi^{-q}(s+1).
\end{equation}
Therefore, for   such $s$ and $m\in   [V_{s-1},V_s)$, the following relation holds
 \[
   \tsfrac {V_{s}-m}{\sum\limits _{j=1}^{s}\tsfrac{\nu_j}{\psi^q(j)}}
   \asymp \tsfrac {V_{s}-m}{\sum\limits _{j=1}^{s}\tsfrac{j^{d-1}}{\psi^q(j)}} \asymp \tsfrac {(V_{s}-m)\psi^{q}(s)}{s^{d-1} }\gg\psi^{q}(s+1),
 \]
which implies that $l_m\le V_{s}$, and  therefore, $l_m=V_{s}$ and $n_{l_m}=s$.
 Combining this equality, relations \eqref{Eq:PrTh2_08},  \eqref{Eq:PrTh4_01} and \eqref{Eq:PrTh4_02},
 and taking into account relations \eqref{Eq:ExV_BnA13} and \eqref{Eq:OE_Sp16},
 we see that  relation \eqref{Eq:OE_Sp17} is true indeed:
 \[
  \sigma_m({\mathcal F}_{q,r}^{\psi})_{_{\scriptstyle {\mathcal S}^p}}\asymp \bigg(
  \Big(\tsfrac {(V_{s}-m)\psi^{q}(s)}{s^{d-1} }\Big)^\frac{q}{q-p}
  \tsfrac{s^{d-1}} {\psi^q(s)}
  +   s^{d-1}{\psi^\frac{pq}{q-p}(s+1)}\bigg)^{\frac{q-p}{pq}}
  \]
 \[
  \asymp  \tsfrac {(V_{s}-m)^p\psi(s)}{
  s^{\frac{d-1}{q}}}
  \bigg(1 +  \Big(\tsfrac{s^{\frac {d-1}{p}}\psi (s+1)}{\psi(s)(V_{s}-m)^p}\Big)^\frac{pq}{q-p}\bigg)^{\frac{q-p}{pq}}
   \asymp \psi(s) \tsfrac {(V_{s}-m)^p}{
  m^{\frac{d-1}{dq}}}.
 \]
\emph{{\rm (ii)} Case $0<p<q=\infty$.} Due to \eqref{Eq:ExV_BnA05} and $m\in   [V_{s-1},V_s)$, we have
 \[
 \sigma_m^p({\mathcal F}_{q,r}^{\psi})_{_{\scriptstyle {\mathcal S}^p}}= (V_{s}-m){\psi^p(s)}+ \sum\limits_{n=s+1}^\infty
  {\nu_n}{\psi^p(n)}.
 \]
Taking into account \eqref{Eq:ExV_BnA11} and \eqref{Aux_Eq15}, we obtain
\[
  \sigma_m^p({\mathcal F}_{q,r}^{\psi})_{_{\scriptstyle {\mathcal S}^p}}
  \asymp (V_{s}-m){\psi^p(s)}+s^{d-1} \psi^p(s+1),
 \]
and by virtue of \eqref{Eq:OE_Sp16}, we see that for $0<p<q=\infty$,  relation \eqref{Eq:OE_Sp17} also holds.

\emph{Case $0<q\le p<\infty$.} By virtue of \eqref{Eq:PrTh2_01}  and  \eqref{Aux_Eq14},   we have
   \[
   \sigma_m^p({\mathcal F}_{q,r}^{\psi})_{_{\scriptstyle {\mathcal S}^p}}
   =  \sup\limits_{l> m} \frac{l-m}{
  \Big(\sum\limits_{n=1}^{n_l-1} \,\tsfrac {n^{d-1}}{\psi^q(n)}+
  \tsfrac{l-V_{n_l-1}}{\psi^q(n_l)}\Big)^{\frac pq} }
  \asymp\sup\limits_{l>m,\, l\in {\mathbb N}} \frac{l-m}{
  \Big( \tsfrac {(n_l-1)^{d-1}}{\psi^q(n_l-1)}+
  \tsfrac{l-V_{n_l-1}}{\psi^q(n_l)}\Big)^{\frac pq}}.
  \]
 In this case, condition \eqref{Eq:OE_Sp16} is also satisfied for    $\beta=\tsfrac{d-1}{p}$.    Therefore,
 \begin{equation}\label{Eq:PrTh4_03}
   \sigma_m^p({\mathcal F}_{q,r}^{\psi})_{_{\scriptstyle {\mathcal S}^p}}\asymp \sup\limits_{l\in {\mathbb N},\,l>m}\tsfrac{\psi^p(n_l)(l-m)}{(l-V_{n_l-1})^{\frac pq}}=
   \sup\limits_{j>m}   \Big( \psi^p(n_j)
    \max\limits_{l\in I_j} \tsfrac{l-m}{(l-V_{n_j-1})^{\frac pq}}\Big),
    \end{equation}
 where $I_{m+1}$ denotes the sets of  positive integers from the half-interval
 $(m,V_s]$ and  $I_{j}$, $j=m+2, m+3,\ldots$, denotes the sets of  positive integers $l$ from the
 half-intervals $(V_{n_{j}-1},V_{n_{j}}]$.

By virtue of \eqref{Eq:ExV_BnA10}, \eqref{Eq:ExV_BnA12} and \eqref{Eq:OE_Sp16}, for $j>m$  we have
\[
  \psi^p(n_{j})\max\limits_{l\in I_j}\tsfrac{l-m}{l-V_{n_j-1}}
  = \tsfrac{\psi^p(n_{j})(V_{n_{j}}-m)}{V_{n_{j}}-V_{n_{j}-1}}\asymp
 \tsfrac{\psi^p(n_{j})(V_{n_{j}}-m)}{n_{j}^{d-1}}
 \gg \psi^p(n_j+1)\max\limits_{l\in I_{j+1}}\tsfrac{l-m}{l-V_{n_j}}.
 \]
Thus, if $p=q$, then the supremum in \eqref{Eq:PrTh4_03} is attained at the point $j=m+1$, and in this case, relation \eqref{Eq:OE_Sp17} holds:
 \[
   \sigma_m({\mathcal F}_{p,r}^{\psi})_{_{\scriptstyle {\mathcal S}^p}}\asymp \psi(n_{m+1})\tsfrac{(V_{n_{m+1}}-m)^\frac 1p}{n_{m+1}^\frac{d-1}{p}}\asymp \psi(n_{m+1})\tsfrac{(V_{n_{m+1}}-m)^\frac 1p}{m^\frac{d-1}{dp}}
   \asymp \psi(s)\tsfrac{(V_{s}-m)^\frac 1p}{m^\frac{d-1}{dp}}.
 \]

Now, let $0<q<p<\infty$.   For   fixed  $m$ and $v>m$,  the function
 \begin{equation}\label{Eq:PrTh4_04}
 f_m(l,v)=\tsfrac{l-m}{(l-v)^{\frac pq}},\qquad l>v,
 \end{equation}
 attains its maximal value at the point $l^*_v=\tsfrac {p m -qv}{p-q}$,
 \begin{equation}\label{Eq:PrTh4_05}
  f_m(l^*_v,v)
   =\tsfrac{q(p-q)^{\frac pq-1}}{p^{\frac pq}}\tsfrac 1{(m-v)^{\frac pq-1}}\asymp
   \tsfrac 1{(m-v)^{\frac pq-1}},
\end{equation}
$ f_m(l,v)$ increases for $l\in (v,l_v^*)$ and decreases for $l>l_v^*$.

If $m>V_{s-1} $,  according to   \eqref{Eq:ExV_BnA10}, \eqref{Eq:OE_Sp16}  and \eqref{Eq:PrTh4_05},  we have
\[
  \psi^p(s)\max\limits_{l\in I_{m+1}}\tsfrac{l-m}{(l-V_{s-1})^{\frac pq}}\gg
  \tsfrac{\psi^p(s) }{(V_{s}-V_{s-1})^{\frac pq}}
  \asymp \tsfrac{\psi^p(s) }{s^{\frac {(d-1)p}{q}}}\gg \psi^p(s+1)
\]
\[
\gg \tsfrac {\psi^p(s+1)}{(m-V_{s})^{\frac pq-1}}=\psi^p(s+1)
\sup\limits_{l>V_{s}} f_m(l,V_{s})\ge
\psi^p(n_j)\max\limits_{l\in I_{j}}\tsfrac{l-m}{(l-V_{n_j-1})^{\frac pq}}
\]
for any $j>m+1$.  Thus, the supremum in \eqref{Eq:PrTh4_03} is attained at $j=m+1$,  and
\[
  \sigma_m^p({\mathcal F}_{q,r}^{\psi})_{_{\scriptstyle {\mathcal S}^p}}\asymp
  \psi^p(s)\max\limits_{l\in  I_{m+1}}  \tsfrac{(l-m)}{(l-V_{s-1})^{\frac pq}}
  \asymp  \psi^p(s)\max\limits_{l\in  (m,V_{s}]}
  f(l,v_m),
 \]
 where $f(l,v)$ is defined by \eqref{Eq:PrTh4_04} and $v_m=V_{s-1}$.
 Further, if condition \eqref{Eq:OE_Sp19} holds, the point $l^*_{v_m}=\tsfrac {p m -qV_{s-1}}{p-q}$
 belongs to   $(m,V_{s}]$, and taking into account \eqref{Eq:PrTh4_05}, we obtain \eqref{Eq:OE_Sp18}:
  \[
  \sigma_m ({\mathcal F}_{q,r}^{\psi})_{_{\scriptstyle {\mathcal S}^p}}\asymp   \psi (s)
  f^\frac 1p(l^*_{r_m},v_m)\asymp
   \tsfrac {\psi (s)}{(m+1-V_{s-1})^{\frac 1q-\frac 1p}}.
 \]
  If condition \eqref{Eq:OE_Sp19} is not satisfied, then   $ f(l,v_m)$ increases on $(m,V_{s}]$. Therefore, taking into account \eqref{Eq:ExV_BnA10}, we obtain \eqref{Eq:OE_Sp17}:
 \[
  \sigma_m ({\mathcal F}_{q,r}^{\psi})_{_{\scriptstyle {\mathcal S}^p}}\asymp   \psi (s)
  f^\frac 1p(V_{s},r_m)\asymp
   \tsfrac{\psi (s)(V_{s}-m)^\frac 1p}{(V_{s}-V_{s-1})^{\frac 1q}}
   \asymp \tsfrac{\psi (s)(V_{s}-m)^\frac 1p}{ m^{\frac {d-1}{dq}}}.
 \]
Similarly, if  $m=V_{s-1}$,  for any $j>m$ we have
 \[
  \psi^p(s)\max\limits_{l\in I_{m+1}}\tsfrac{l-m}{(l-V_{s-1})^{\frac pq}}
  =\psi^p(s)\max\limits_{l\in I_{m+1}}(l-V_{s-1})^{-\frac pq+1}=\psi^p(s)
  \]
\[
\gg \tsfrac {\psi^p(s+1)}{(m+1-V_{s})^{\frac pq-1}}=\psi^p(s+1)
\sup\limits_{l>V_{s}} f_m(l,V_{s})\ge
\psi^p(n_j)\max\limits_{l\in I_{j}}\tsfrac{l-m}{(l-V_{n_j-1})^{\frac pq}}.
\]
 Therefore, in this case relation \eqref{Eq:OE_Sp18} holds:
\[
   \sigma_m ({\mathcal F}_{q,r}^{\psi})_{_{\scriptstyle {\mathcal S}^p}}\asymp  \psi (s)=
    \tsfrac {\psi (s)}{(m+1-V_{s-1})^{\frac 1q-\frac 1p}}.
\]

\emph{{\rm (iii)} Case $0<q<p=\infty$.} Due to \eqref{Eq:ExV_BnA04},  \eqref{Eq:ExV_BnA07} and \eqref{Eq:PrTh4_02}, we have
\[
  \sigma_m^q({\mathcal F}_{q,r}^{\psi})_{_{\scriptstyle {\mathcal S}^\infty}}
  =  \Big( \sum_{n=1}^{s-1} \,\tsfrac {\nu_n}{\psi^q(n)}+
  \tsfrac{m+1-V_{s-1}}{\psi^q(s)}\Big)^{-\frac 1q} \asymp
  \Big(\tsfrac {(s-1)^{d-1}}{\psi^q(s-1)}+
  \tsfrac{m+1-V_{s-1}}{\psi^q(s)}\Big)^{-\frac 1q}
  \]
  \[
   \asymp \tsfrac{\psi (s)}{(m+1-V_{s-1})^\frac 1q} \Big(1+
   \tsfrac{s^{d-1}\psi^q(s)}{\psi^q(s-1)}\tsfrac {1}{m+1-V_{s-1}}
  \Big)^{-\frac 1q}.
  \]
  Since  $\psi$ satisfies  \eqref{Eq:OE_Sp16} with   $\beta =\tsfrac{d-1}{q}$, we conclude that relation \eqref{Eq:OE_Sp18} also holds.

 {\rm (iv)} Find estimates for  ${\mathscr D}_m({\cal F}_{q,r}^{\psi})_{_{\scriptstyle {\mathcal S}^p }}$
 in the \emph{case $0<p<q\le\infty$.}  Due to \eqref{Eq:ExV_BaW07} and \eqref{Eq:ExV_BnA07},
 \[
  {\mathscr D}_m^{\frac{q-p}{p\,q}}({\mathcal F}_{q,r}^{\psi})_{_{\scriptstyle {\mathcal S}^p}}
  = (V_{s}-m){\psi^{\frac{p\,q}{q-p}}(s)}+ \sum\limits_{n=s+1}^\infty
  {\nu_n}{\psi^{\frac{p\,q}{q-p}}(n)}.
 \]
Taking into account \eqref{Eq:ExV_BnA11} and \eqref{Aux_Eq15}, we obtain
\[
  {\mathscr D}_m^{\frac{q-p}{p\,q}}({\mathcal F}_{q,r}^{\psi})_{_{\scriptstyle {\mathcal S}^p}}
  \asymp (V_{s}-m){\psi^{\frac{p\,q}{q-p}}(s)}+(s+1)^{d-1} \psi^{\frac{p\,q}{q-p}}(s+1).
 \]
By  virtue of \eqref{Eq:OE_Sp16}, we see that for $0<p<q<\infty$ relation \eqref{Eq:OE_Sp20}  is valid.
 \bend


\subsection{Proof of Theorem \ref{Th:_OE_Lp1}}\label{Pr_Th:_OE_Lp1}

 \textit{Upper estimates.} In the case $1\le p\le 2$, due to \eqref{Eq:ApCh12}, \eqref{Eq:OE_Lp02}  and \eqref{Eq:ApCh13}, we have
 \begin{equation}\label{Eq:PT_Lp01}
  \sigma_m({\mathcal F}_{q,r}^{\psi})_{_{\scriptstyle L_p}}\ll
  \sigma_m^\perp({\mathcal F}_{q,r}^{\psi})_{_{\scriptstyle L_p}}\ll
  G_m({\mathcal F}_{q,r}^{\psi})_{_{\scriptstyle L_p}}\ll
  G_m({\mathcal F}_{q,r}^{\psi})_{_{\scriptstyle L_2}}\ll
  \sigma_m({\mathcal F}_{q,r}^{\psi})_{_{\scriptstyle {\mathcal S}^{2}}}.
\end{equation}
Thus,  to obtain the necessary upper bound, it suffices to use    \eqref{Eq:OE_Sp04} for ${\mathcal S}^p={\mathcal S}^{2}$.

In the case  $2\le p<\infty$, using relations \eqref{Eq:ApCh12}, \eqref{Eq:OE_Lp01}   and  \eqref{Eq:ApCh13}, we obtain
 \begin{equation}\label{Eq:PT_Lp02}
   \sigma_m({\mathcal F}_{q,r}^{\psi})_{_{\scriptstyle L_p}}\ll
   \sigma_m^\perp({\mathcal F}_{q,r}^{\psi})_{_{\scriptstyle L_p}}\ll
   G_m({\mathcal F}_{q,r}^{\psi})_{_{\scriptstyle L_p}} \ll
   G_m({\mathcal F}_{q,r}^{\psi})_{_{\scriptstyle {\mathcal S}^{p'}}}\asymp
  \sigma_m({\mathcal F}_{q,r}^{\psi})_{_{\scriptstyle {\mathcal S}^{p'}}},
 \end{equation}
 and the upper bound similarly follows from relation \eqref{Eq:OE_Sp04}.

\textit{ Lower estimate.}  Let ${\mathscr T}_n$, $n\in {\mathbb N}$, denote the set of all polynomials of the form
 \[
 \tau_n=\sum\limits_{|k|_{\infty}\le n}\widehat{\tau}_n(k)e_k,
 \]
  and let ${\mathcal A}_q({\mathscr T}_n)$, $0<q\le\infty$, denote the subset of
  all polynomials $\tau_n\in {\mathscr T}_n$ such that $\|\tau_n\|_{_{\scriptstyle {\mathcal S}^{q}}}\le 1$.
  From Theorem 5.2 of \cite{DeVore_Temlyakov_1995}, it follows that for any
  $0<q\le \infty$, $1\le p\le\infty$, $n=1,2,\ldots$ and $m=\tsfrac{(2n+1)^d-1}2$,
 \[
 \sigma_m({\mathcal A}_q({\mathscr T}_n))_{_{\scriptstyle L_p}}\ge K m^{\frac 12-\frac 1q}.
 \]
 For a fixed $n\in {\mathbb N}$, consider the set
 $\psi(dn){\mathcal A}_q({\mathscr T}_n)=\{\tau \in {\mathscr T}_n:
 \|\tau\|_{_{\scriptstyle {\mathcal S}^{q}}}\le \psi(dn)\}.$
 Due to monotonicity $\psi$, for any   $\tau\in \psi(dn){\mathcal A}_q({\mathscr T}_n)$ and $0<r\le \infty$,  we have
$$
\sum\limits_{k\in {\mathbb Z}^d} \left|\tsfrac{\widehat{\tau}(k)}{\psi(|k\|_{r})} \right|^q\le \sum\limits_{|k|_{\infty}\le n} \left|\tsfrac {\widehat{\tau}(k)}{\psi(d|k|_{\infty})} \right|^q \le\sum\limits_{|k|_{\infty}\le n} \left|\tsfrac{\widehat{\tau}(k)}{\psi(dn)} \right|^q\le
1
$$
Therefore, $\psi(dn){\mathcal A}_q({\mathscr T}_n)$ is contained in the set ${\mathcal F}_{q,r}^{\psi}$. In view of definition of the set $B$, for all $n=1,2,\ldots$ and $m=\tsfrac{(2n+1)^d-1}2$, we obtain
 $$
 \sigma_m({\mathcal F}_{q,r}^{\psi})_{_{\scriptstyle L_p}}\ge  \sigma_m(\psi(dn){\mathcal A}_q({\mathscr T}_n))_{_{\scriptstyle L_p}}\gg  \psi(dn)m^{\frac 12- \frac 1q}\gg  \psi(m^\frac 1d)m^{\frac 12-\frac 1q}.
 $$
Taking into account the relation \eqref{Eq:ApCh12}, monotonicity of the quantity $\sigma_m$ and inclusion
$\psi\in B$, we see that for all $1\le p\le\infty$,
 \[
 G_m({\mathcal F}_{q,r}^{\psi})_{_{\scriptstyle L_p}}\gg \sigma_m^\perp({\mathcal F}_{q,r}^{\psi})_{_{\scriptstyle L_p}}\gg \sigma_m({\mathcal F}_{q,r}^{\psi})_{_{\scriptstyle L_p}}\gg
  \psi(m^\frac 1d)m^{\frac 12-\frac 1q}.
 \]

In the case $2< p<\infty$, for the quantities $\sigma_m^\perp({\mathcal F}_{q,r}^{\psi})_{_{\scriptstyle L_p}}$ and $G_m({\mathcal F}_{q,r}^{\psi})_{_{\scriptstyle L_p}},$ this estimate can be improved. For this purpose, consider the function
\[
h_{1}=\sum\limits_{|k|_r\le n_m} \widehat{h}_1(k)e_k={\mathfrak h}_1\sum\limits_{|k|_r\le n_m} e_k,
\]
 where ${\mathfrak h}_1^{-q}(m):=\sum_{|j|_r\le n_m}\psi^{-q}(|j|_r)$, $n_m:=n_m(r,d)=[(\tsfrac {2m}{M_{r,d}})^{\frac 1d}]$
 and  $M_{r,d}$ is a constant defined in \eqref{Eq:ExV_BnA09}.
It is obviously that $f_1\in {\mathcal F}_{q,r}^{\psi}$ and by  \eqref{Aux_Eq003_1} and \eqref{Eq:OE_Sp01}
  \[
  {\mathfrak h}_1^{-q}(m)\asymp \sum\limits_{j=1}^{n_m}\tsfrac {j^{d-1}}{\psi^{q}(j)}\asymp
  \tsfrac {n_m^{d}}{\psi^q(n_m)}\asymp \tsfrac m {\psi^q(m^{\frac  1d})}.
  \]
 For any collection $\gamma_n\subset {\mathbb Z}^d$, using Nikol'skii's inequality \cite{Nikol'skii_1951} and \eqref{Eq:ExV_BnA10}, we obtain
 \[
 \bigg|\bigg|h_1-\sum\limits_{k\in \gamma_m} \widehat{h}_1(k) e_k\bigg|\bigg|_{_{\scriptstyle L_p}}
 \gg {\mathfrak h}_1(m) m^{-\frac 1p}\bigg|\bigg| \sum\limits_{|k|_1\le n_m,\,k\notin \gamma_m}    e_k\bigg|\bigg|_{_{\scriptstyle L_{\infty}}}  \asymp \psi(m^{\frac 1d})m^{1-\frac 1p-\frac 1q}.
 \]
 Therefore, for all $2\le p<\infty$, the following estimates are true:
 \[
   G_m({\mathcal F}_{q,r}^{\psi})_{_{\scriptstyle L_p}}\gg
   \sigma_m^\perp({\mathcal F}_{q,r}^{\psi})_{_{\scriptstyle L_p}}\gg
   \sigma_m^\perp(h_1)_{_{\scriptstyle L_p}}\gg \psi(m^{\frac 1d})m^{1-\frac 1p-\frac 1q}.
 \]
\bend

\subsection{ Proof of Theorem \ref{Th:_OE_Lp2}}
Proof of Theorem \ref{Th:_OE_Lp2} is  similar to the proof of the upper estimates in
Theorem 6.1 \cite{DeVore_Temlyakov_1995}. It uses the following lemma from \cite{DeVore_Temlyakov_1995}.

\begin{lemma}{\rm \cite{DeVore_Temlyakov_1995}}
 For any $0<q\le \infty$,  $n=1,2,\ldots$, and $1\le m\le (2n+1)^d$,
 \begin{equation}\label{Eq:PT_Lp03}
\sigma_m({\mathcal A}_q({\mathscr T}_n))_{_{\scriptstyle L_\infty}}\le C m ^{\frac 12-\frac 1q} L(n^d/m), \quad 0<q\le 1,
\end{equation}
where $L(x)=(1+(\ln x)_+)^{1/2}$ and
 \begin{equation}\label{Eq:PT_Lp04}
\sigma_m({\mathcal A}_q({\mathscr T}_n))_{_{\scriptstyle L_\infty}}\le C n^{d-\frac dq} m^{-\frac 12} L(n^d/m), \quad 1<q\le \infty,
\end{equation}
with $C$ depending only on $q$ and $d$.
\end{lemma}

 For any $f\in {\mathcal F}_{q,\infty}^{\psi}$, we use the decomposition
 \[
  f=\sum\limits_{j=0}^\infty f_j,\qquad \mbox{where}\quad
  f_j:=\sum_{2^{j-1}\le |k|_\infty<2^j} \widehat f(k)e_k,\quad j\ge 1,\quad f_0:=\widehat f(0).
  \]
  We note that
\begin{equation}\label{Eq:PT_Lp05}
f_j/\psi(2^{j-1})\in {\mathcal A}_q({\mathscr T}_{2^j}),\quad j=1,2,\ldots
\end{equation}

For any  $N=1,2,\ldots$, we approximate $f$ as follows.
Let $N_0$ be the largest integer $j$ such that $m_j:=[(j-N)^{-2}2^{Nd}]\ge 1$,
i.e. $N_0=[2^{\frac{Nd}2}+N]$. If $j\le N$, we set $P_j:=f_j$. If $N<j\le N_0$, then by virtue of
\eqref{Eq:PT_Lp03}, \eqref{Eq:PT_Lp04} and \eqref{Eq:PT_Lp05}, there is a polynomial $P_j\in \Sigma_{m_j}$  such that
 \begin{equation}\label{Eq:PT_Lp06}
\|f_j-P_j\|_{_{\scriptstyle L_p}}\ll  m_j ^{\frac 12-\frac 1q} L({2^{jd}}/{m_j})\psi(2^{j-1}), \quad 0<q\le 1.
\end{equation}
and
 \begin{equation}\label{Eq:PT_Lp07}
\|f_j-P_j\|_{_{\scriptstyle L_p}}\ll   2^{j(d-\frac dq)} m_j^{-\frac 12} L({2^{jd}}/{m_j})\psi(2^{j-1}), \  1<q<\infty.
\end{equation}

Set $P=\sum\limits_{j=0}^{N_0}P_k$. Since
 $$
 (2\cdot 2^N+1)^d+\sum\limits_{j=N+1}^{N_0}(j-N)^{-2}2^{Nd}\le a2^{Nd},
 $$
where $a$ depends only on $d$, then $P$ is a linear combination of at most $a2^{Nd}$ exponentials $e_k$. Hence, $P$ is in $\Sigma_{a2^{Nd}}$. We also have
 \begin{equation}\label{Eq:PT_Lp08}
\|f-P\|_{_{\scriptstyle L_p}}\le \sum\limits_{j=N+1}^{N_0} \|f_j-P_j\|_{_{\scriptstyle L_p}}+\sum\limits_{j=N_0+1}^\infty \|f_j\|_{_{\scriptstyle L_p}}=:S_1+S_2.
\end{equation}
For all  $x\ge 1$, we have $[x]\ge x/2$. Therefore, for sufficiently large $N$ and  $N<j\le N_0$, from the definition of $L(x)$, we have
 \begin{equation}\label{Eq:PT_Lp09}
L({2^{jd}}/{m_j})\le (1+\ln (2^{d(j-N)+1}(j-N)^2))^{\frac 12}\ll (j-N)^{\frac 12}.
\end{equation}

First, consider the case $0<q\le 1$. Reasoning similar to the proof of upper estimate in \eqref{Aux_Eq003_2},
it is easy to show  that if the function $\psi$ belongs to the set $B$ and  for all\ \
$t$, larger than a certain number  $t_0$, $\psi$ is convex  and  it satisfies  condition
\eqref{Eq:OE_Sp02}  with a fixed $\beta\ge 0$,  then for any $\alpha\in {\mathbb R}$ and sufficiently large $t>N$, the function $h_{\alpha,\beta}(t):=2^{\beta t}  (t-N)^{\alpha} \psi(2^{t-1})$ decreases to zero, as well as
 \begin{equation}\label{Eq:PT_Lp10}
\sum\limits_{j=N+1}^\infty 2^{\beta j}  (j-N)^{\alpha} \psi(2^{j-1})\ll 2^{\beta N} \psi(2^N).
\end{equation}

 In this case, $\beta=0$.   By virtue of \eqref{Eq:PT_Lp06}, \eqref{Eq:PT_Lp09} and \eqref{Eq:PT_Lp10}, we obtain the estimate of the first sum $S_1$ in (\ref{Eq:PT_Lp08}):
$$
S_1\ll \sum\limits_{j=N+1}^\infty (j-N)^{2(\frac 1q-\frac 12)} 2^{-Nd(\frac 1q-\frac 12)}(j-N)^{\frac 12}\psi(2^{j-1})\ll
$$
 \begin{equation}\label{Eq:PT_Lp11}
\ll 2^{-N(\frac dq-\frac d2)}\sum\limits_{j=N+1}^\infty (j-N)^{\frac 2q-\frac 12} \psi(2^{j-1})\ll 2^{-N(\frac dq-\frac d2)}\psi(2^{N}).
\end{equation}

To estimate $S_2$, we note that from \eqref{Eq:PT_Lp05}
 $$
 S_2\le
   \sum\limits_{j=N_0+1}^\infty \|f_j\|_{_{\scriptstyle L_\infty}}\le  \sum\limits_{j=N_0+1}^\infty \Big(\sum_{2^{j-1}\le \|k\|_\infty<2^j} |\widehat f(k)|\Big) $$
   $$
  \le  \sum\limits_{j=N_0+1}^\infty  \|f_j\|_{_{\scriptstyle {\mathcal S}^q}}\ll
  \sum\limits_{j=N_0+1}^\infty \psi(2^{j-1}) \ll \psi(2^{N_0}).
 $$

 Further,  note that if  for all\ \  $t$, larger than a certain number  $t_0$, $\psi$
 is convex  and satisfies \eqref{Eq:OE_Sp02}, then for any $\alpha>0$, we have $
\psi(2^{N(\alpha+1)})\ll \psi(2^N) 2^{-N\alpha}$.

From the definition of $N_0$, we have $N_0\ge N+2^{\frac{Nd}2}-1$. It follows that if $N$ is sufficiently large (depending only on $d$ and $q$), then
$N_0\ge N(1+\frac dq-\frac d2)$. Hence,
 \[
 S_2\ll \psi(2^{N(1+\frac dq-\frac d2)})\ll 2^{-N(\frac dq-\frac d2)}\psi(2^N) .
 \]

Using this and \eqref{Eq:PT_Lp11} in \eqref{Eq:PT_Lp08}, we find that
\begin{equation}\label{Eq:PT_Lp12}
 \sigma_{a 2^{Nd}}(f)_{_{\scriptstyle L_p}}\le \|f-P\|_{_{\scriptstyle L_p}}\ll 2^{-N(\frac dq-\frac d2)}\psi(2^N).
\end{equation}

In the case $1<q<\infty$, condition \eqref{Eq:OE_Sp02}  is satisfied
with $\beta=d-\frac dq$. By virtue of \eqref{Eq:PT_Lp07}, \eqref{Eq:PT_Lp09} and \eqref{Eq:PT_Lp10}, we have
\begin{equation}\label{Eq:PT_Lp13}
S_1
 \ll 2^{-\frac N2}\sum\limits_{j=N+1}^{\infty} 2^{j(d-\frac dq)}  (j-N)^{\frac 32}\psi(2^{j-1})\ll \psi(2^N)2^{-N(\frac dq-\frac d2)}.
\end{equation}
To estimate $S_2$, we use H\"{o}lder's inequality, \eqref{Eq:PT_Lp05}, \eqref{Eq:PT_Lp10} and the inequalities
$N_0\ge N+2^{\frac{Nd}2}-1$ and $h_{\alpha,\beta}(N_0)\le h_{\alpha,\beta}(N+1)$ with
$\alpha=1$ and $\beta=d-\tsfrac dq$,
 $$
  S_2\le
   \sum\limits_{j=N_0+1}^\infty \|f_j\|_{_{\scriptstyle L_\infty}}\le
 \sum\limits_{j=N_0+1}^\infty \psi(2^{j-1})\Big(\sum_{2^{j-1}\le \|k\|_\infty<2^j}
\Big|\tsfrac{\widehat f(k)}{ \psi(2^{j-1})}\Big|\Big)\le   $$
   $$
  \le  \sum\limits_{j=N_0+1}^\infty \psi(2^{j-1}) 2^{(j-1)(d-\frac dq)} \ll
  \psi(2^{N_0}) 2^{N_0(d-\frac dq)}\ll
  \psi(2^{N}) 2^{N(\frac d2-\frac dq)}.
 $$
Using this and \eqref{Eq:PT_Lp13} in \eqref{Eq:PT_Lp08}, we see that in this case, relation \eqref{Eq:PT_Lp12} is also true.
 Therefore, the upper estimate in \eqref{Eq:OE_Sp02}  follows from the monotonicity of $\sigma_m$ and
 inclusion  $\psi\in B$. \bend

\subsection{ Proof of Theorem \ref{Th:_OE_Lp_3}}

 \textit{Upper estimates.} Let $m=m(s)=V_{s}-c_s$, $1\le c_s\le c$. Then due to
 relations \eqref{Eq:PT_Lp01} and \eqref{Eq:PT_Lp02} and estimate \eqref{Eq:OE_Sp24}, we see that in this case
 for any  $0<q\le\infty$ and $1\le p<\infty$
    \begin{equation}\label{Eq:PT_Lp14}
    \sigma_m^\perp({\cal F}_{q,r}^{\psi})_{_{\scriptstyle L_p}}\ll
    G_m({\cal F}_{q,r}^{\psi})_{_{\scriptstyle L_p}}\ll
    \tsfrac{\psi(s)}{m^{\frac{d-1}{qd}}}.
    \end{equation}
If $m=m(s)=V_{s-1}+c_s$, $0\le c_s\le c$, and   $0<q< p'<\infty$ or if $m=m(s)=V_{s-1}$ and
  $0<p'=q<\infty$, then the  upper estimate similarly follows from
  \eqref{Eq:PT_Lp01}, \eqref{Eq:PT_Lp02} and  \eqref{Eq:OE_Sp29}
     \begin{equation}\label{Eq:PT_Lp15}
    \sigma_m^\perp({\cal F}_{q,r}^{\psi})_{_{\scriptstyle L_p}}\ll
    G_m({\cal F}_{q,r}^{\psi})_{_{\scriptstyle L_p}}\ll  \psi(s).
    \end{equation}

 \textit{Lower estimates.} Let $\{k_1^*,k_2^*,\ldots\}$ be a rearrangement
 of vectors from    ${\mathbb Z}^d$ such that
 \begin{equation}\label{Eq:PT_Lp16}
  |\Psi_{k_j^*}|=\bar \Psi_j, \quad j=1,2,\ldots,
 \end{equation}
 where, as above,  $\bar\Psi= \{\bar \Psi_j\}_{j=1}^\infty$  is a non-increasing rearrangement
 of the  system $\{|\Psi_k|\}_{k \in {\mathbb Z}^d}$.

 Consider the function
  \[
 h_{2}=\sum\limits_{j=1}^{m+1}\widehat{h}_3(k_j^*) e_{k_j}=
 {\mathfrak h}_2(m)\sum\limits_{j=1}^{m+1}e_{k_j},\qquad
 \mbox{\rm where } \quad {\mathfrak h}_2^{-q}(m)=\sum\limits_{i=1}^{m+1}|\Psi_{k_i^*}|^{-q}.
 \]
 It is easy to see that $h_2\in {\cal F}_{q,r}^{\psi}$  and due to \eqref{Eq:PT_Lp16}  and \eqref{Eq:ExV_BnA07},
 for   $m\in [V_{s-1}, V_s)$ we have
 \begin{equation}\label{Eq:PT_Lp17}
 {\mathfrak h}_2^{-q}(m)=\sum\limits_{j=1}^{m+1}\bar \Psi_j^{-q}
 \asymp \sum\limits_{k=1}^{s-1}\tsfrac {V_k-V_{k-1}}{\psi^{q}(k)}+\tsfrac{m+1-V_{s-1}}{\psi^{q}(s)},
 \end{equation}
where in view of   \eqref{Eq:ExV_BnA11} and \eqref{Aux_Eq14},
 \begin{equation}\label{Eq:PT_Lp18}
 \sum\limits_{k=1}^{s-1}\tsfrac {V_k-V_{k-1}}{\psi^{q}(k)}\asymp
 \sum\limits_{k=1}^{s-1}\tsfrac {k^{d-1}}{\psi^{q}(k)}\asymp  \tsfrac {(s-1)^{d-1}}{\psi^{q}(s-1)}.
 \end{equation}
 Combining  \eqref{Eq:PT_Lp17} and \eqref{Eq:PT_Lp18},  taking into account \eqref{Eq:OE_Sp16}, we conclude that
 \[
 {\mathfrak h}_2^{-q}(m)  \asymp
\tsfrac{m+1-V_{s-1}}{\psi^{q}(s)} \Big(1+ \tsfrac {(s-1)^{d-1}\psi^q(s)}{\psi^{q}(s-1)(m+1-V_{s-1})}\Big)\asymp
 \tsfrac{m+1-V_{s-1}}{\psi^{q}(s)}.
 \]
 Then for any collection $\gamma_m\in \Gamma_m$, we have
 \[
 \Big\|h_2-\sum\limits_{k\in \gamma_m} \widehat{h}_2(k) e_k\Big\|_{_{\scriptstyle L_p}}
 ={\mathfrak h}_2(m) \Big\|\sum\limits_{j\in [1,m+1]\setminus\gamma_m}
  e_{k_j^*}\Big\|_{_{\scriptstyle L_p}}\asymp  \tsfrac{\psi(s)}{(m+1-V_{s-1})^\frac 1q}.
 \]
Therefore, for all $1\le p<\infty$, the following estimate holds:
   \begin{equation}\label{Eq:PT_Lp19}
 G_m({\cal F}_{q,r}^{\psi})_{_{\scriptstyle L_{p}}}\gg \sigma_n^\perp({\cal F}_{q,r}^{\psi})_{_{\scriptstyle L_{p}}}\gg
 \sigma_n^\perp(h_2)_{_{\scriptstyle L_p}}\gg   \tsfrac{\psi(s)}{(m+1-V_{s-1})^\frac 1q}.
 \end{equation}
 If $m=m(s)=V_{s-1}+c_s$, $0\le c_s\le c$, then the  quantities on the right-hand sides of relations
 \eqref{Eq:PT_Lp19} and \eqref{Eq:PT_Lp15} are equivalent.
 If $m=m(s)=V_{s}-c_s$, $1\le c_s\le c$, then the equivalence of quantities on the right-hand sides of relations
 \eqref{Eq:PT_Lp19} and \eqref{Eq:PT_Lp14} follows from \eqref{Eq:ExV_BnA11}, \eqref{Eq:ExV_BnA12}
 and \eqref{Eq:ExV_BnA13}.  \bend


\subsection{ Proof of Theorem \ref{Th:_OE_Lp_4}}

Theorem \ref{Th:_OE_Lp_4} can be proven similarly to Theorem \ref{Th:_OE_Lp_3}.
\textit{Upper estimates} for
 $\sigma_m^\perp({\cal F}_{q}^{\psi})_{_{\scriptstyle L_{p}({\mathbb T}^1)}}$ and
 $G_n({\cal F}_{q}^{\psi})_{_{\scriptstyle L_p({\mathbb T}^1)}}$ follows from
 \eqref{Eq:PT_Lp01} and \eqref{Eq:PT_Lp02}, $(\ref{Eq:OE_Sp14}')$,   taking into account
 \eqref{Eq:OE_Sp11}.

 Similarly to \eqref{Eq:PT_Lp01} and \eqref{Eq:PT_Lp02}, we have
\begin{equation}\label{Eq:PT_Lp20}
 {\mathscr D}_m^\perp({\cal F}_{q,r}^{\psi})_{_{\scriptstyle L_{p}}}\le
 {\mathscr D}_m^\perp({\cal F}_{q,r}^{\psi})_{_{\scriptstyle L_{2}}}
 ={\mathscr D}_m({\cal F}_{q,r}^{\psi})_{_{\scriptstyle {\mathcal S}^{2}}}, \qquad 1\le p\le 2,
\end{equation}
 and
\begin{equation}\label{Eq:PT_Lp21}
 {\mathscr D}_m^\perp({\mathcal F}_{q,r}^{\psi})_{_{\scriptstyle L_p}}\le
 {\mathscr D}_m^\perp({\mathcal F}_{q,r}^{\psi})_{_{\scriptstyle {\mathcal S}^{p'}}}
 ={\mathscr D}_m({\mathcal F}_{q,r}^{\psi})_{_{\scriptstyle {\mathcal S}^{p'}}},\qquad 2\le p<\infty.
\end{equation}
 Therefore, to obtain upper estimates for
 ${\mathscr D}_m^\perp({\cal F}_{q}^{\psi})_{_{\scriptstyle L_{p}({\mathbb T}^1)}}$,
 it sufficient to use estimates \eqref{Eq:PT_Lp20}, \eqref{Eq:PT_Lp21},
 \eqref{Eq:OE_Sp03} and $(\ref{Eq:OE_Sp13}')$,   taking into account
 \eqref{Eq:OE_Sp11}.

 To obtain  \textit{lower estimates}, consider a rearrangement $\{k_1^*,k_2^*,\ldots\}$
 of the set ${\mathbb Z}$ such that
 \begin{equation}\label{Eq:PT_Lp22}
  |\Psi_{k_j^*}|=\bar\Psi_j, \quad j=1,2,\ldots,
\end{equation}
where $\bar\Psi= \{\bar \Psi_j\}_{j=1}^\infty$ is a non-increasing rearrangement  of
the system  $\{|\Psi_k|\}_{k \in {\mathbb Z}}$.

For a given $m\in {\mathbb N}$, consider the function
 \[
 h_{3}=\sum\limits_{j=1}^{m+1}\widehat{h}_3(k_j^*) e_{k_j}={\mathfrak h}_3(m)\sum\limits_{j=1}^{m+1}e_{k_j},\qquad
 \mbox{\rm where } \quad {\mathfrak h}_3^{-q}(m)=\sum\limits_{j=1}^{m+1}|\Psi_{k_i^*}|^{-q}.
 \]
It is clear that  $h_3\in {\cal F}_{q}^{\psi}$, and  due to \eqref{Eq:PT_Lp22}, \eqref{Eq:ExV_BnA07},
\eqref{Aux_Eq04} and \eqref{Eq:OE_Sp11}, we see that
 \[
  {\mathfrak h}_3^{-q}(m)\asymp  \sum\limits_{j=1}^{m+1}\bar\Psi^{-q}_j
 \asymp \sum\limits_{k=1}^{[\frac{m+1}{2}]}{\psi^{-q}(k)}\asymp
 \tsfrac{[\frac{m+1}{2}]\alpha(\psi,[\frac{m+1}{2}])} {\psi^s([\frac{m+1}{2}])} \asymp \psi^{-q}([\tsfrac{m+1}{2}]).
 \]
Therefore, for any collection  $\gamma_m\in \Gamma_m$, we have
 \[
\Big\|h_3-\sum\limits_{k\in \gamma_m} \widehat{h}_3(k) e_k\Big\|_{_{\scriptstyle L_p({\mathbb T}^1)}} =
 {\mathfrak h}_3(m)\Big\|\mathop{\sum\limits_{j=1}^{n+1}}\limits_{k_j^*\notin \gamma_m}
  e_{k_j^*}\Big\|_{_{\scriptstyle L_p({\mathbb T}^1)}}\asymp \psi([\tsfrac{m+1}{2}])
 \]
 and the necessary lower estimate is valid:
 \[
  \sigma_m^\perp({\cal F}_{q}^{\psi})_{_{\scriptstyle L_{p}({\mathbb T}^1)}}
  \gg \sigma_m^\perp(h_3)_{_{\scriptstyle L_p({\mathbb T}^1)}}\gg \psi([\tsfrac{m+1}{2}])\asymp
  \psi(\tsfrac{m}{2}).
 \] \bend

 \subsection{ Proof of Theorem \ref{Th:_OE_Lp_5}}  To obtain
 \textit{upper estimates} in \eqref{Eq:OE_Lp12}
 and \eqref{Eq:OE_Lp13}, it is sufficient to use relations \eqref{Eq:PT_Lp20}, \eqref{Eq:PT_Lp21}
 \eqref{Eq:OE_Sp15} and \eqref{Eq:OE_Sp25}. To obtain \textit{lower estimates},
 for any collection $\gamma_m\in \Gamma_m$,
 we choose any vector
 $k_0=k_0({\gamma_m})\in {\mathbb Z}^d\setminus \gamma_m$
 such that
 \[
 |\psi(|k_0|_r)|=\bar{\Psi}_{\gamma_m}(1)= \sup\limits_{k\in {\mathbb Z}^d\setminus \gamma_m}
 |\Psi_{\gamma_m}(k)|
 = \sup\limits_{k\in {\mathbb Z}^d\setminus \gamma_m}|\psi(|k|_r)|,
 \]
 where, as in  Section \ref{Sec:AC_Sp}, the system
 $\Psi_{\gamma_m}=\{\Psi_{\gamma_m}(k)\}_{k\in {\mathbb Z}^d}$ is defined by \eqref{Eq:ExV_BaW03} and
  $\bar{\Psi}_{\gamma_m}=\{\bar{\Psi}_{\gamma_m}(j)\}_{j=1}^\infty$ is  a non-increasing rearrangement
  of the system $\{|\Psi_{\gamma_m}(k)|\}_{k\in {\mathbb Z}^d}$.

 Then the function  $h_4:=\psi(|k_0|_r)e_{k_0}$  belongs to the set  ${\cal F}_{q,r}^{\psi}$ and
  \[
 {\mathscr E}_{\gamma_m}(h_4)_{_{\scriptstyle L_p}}=\big\|\psi(|k_{0}|_r)e_{k_{0}}
 \big\|_{_{\scriptstyle L_p}}=\psi(|k_{0}|_r)=\bar{\Psi}_{\gamma_m}(1).
 \]
 Taking into account \eqref{Eq:ExV_BnA07}, \eqref{Eq:ExV_BnA12}, \eqref{Eq:ExV_BnA13} and \eqref{Eq:OE_Sp01},
  we obtain the necessary lower estimate:
 \[
  {\mathscr D}_m^\perp( {\cal F}_{q,r}^{\psi})_{_{\scriptstyle L_p}}\ge
  \inf\limits_{\gamma_m\in \Gamma_m} {\mathscr E}_{\gamma_m}(h_4)_{_{\scriptstyle L_p}}
 =\inf\limits_{\gamma_m\in \Gamma_m} \bar{\Psi}_{\gamma_m}(1)=\bar{\Psi}_{m+1}=
 \psi(s).
 \] \bend

  \subsection{ Proof of Theorem \ref{Th:_OE_Lp_6}}  \textit{Upper estimates} in \eqref{Eq:OE_Lp14}
  follow from relations \eqref{Eq:PT_Lp20}, \eqref{Eq:PT_Lp21} and estimate  \eqref{Eq:OE_Sp03}.

  In the case when $0<q\le  \tsfrac {p}{p-1}$, the \textit{lower estimate} follows from \eqref{Eq:OE_Lp12},
  taking into account relations \eqref{Eq:ExV_BnA12}, \eqref{Eq:ExV_BnA13} and the  definition
  of the set $B$.  In the case when $\tsfrac p{p-1}<q\le\infty$, the \textit{lower estimate}
  in \eqref{Eq:OE_Lp14} follows from \eqref{Eq:OE_Lp05}. \bend


   \section{Acknowledgments}
 This work was partially supported by a grant from the Simons Foundation (SFI-PD-Ukraine-00014586, A.Sh.).
     

{\footnotesize
     \begin{thebibliography}{19}\parskip-2pt \addcontentsline{toc}{section}{References}

 \bibitem{Abdullayev_Ozkartepe_Savchuk_Shidlich_2019}
        \textsc{Abdullayev~F.\,G.,  \"{O}zkartepe~P.,  Savchuk~V.\,V.,  Shidlich~A.\,L.}
        Exact constants in direct and inverse approximation theorems for functions of
        several variables in the spaces ${\mathcal S}^p$.
        \emph{Filomat}, 33 (5), 1471-1484  (2019).
        DOI: 10.2298/FIL1905471A.

 \bibitem{Abdullayev_Serdyuk_Shidlich_2021}
       \textsc{Abdullayev~F.,  Serdyuk~A.,  Shidlich~A.}
       Widths of functional classes defined by majorants of generalized moduli of smoothness
       in the spaces ${\mathcal S}^p$.
       \emph{Ukrainian Mathematical Journal}, 73 (6), 841-858  (2021).
       DOI: 10.1007/s11253-021-01963-6

 \bibitem{Abdullayev_Chaichenko_Shidlich_2021}
       \textsc{Abdullayev~F.,  Chaichenko~S.,  Shidlich~A.}
       Direct and inverse approximation theorems of functions in the Musielak-Orlicz type spaces.
       \emph{Mathematical Inequalities and Applications,} 24 (2),  323-336 (2021).
       DOI: 10.7153/mia-2021-24-23

 \bibitem{Barron_1993}
       \textsc{Barron~R.}
        Universal approximation bounds for superpositions of a sigmoidal function.
        \emph{IEEE Transactions on Information Theory}, 39 (3), 930-945 (1993).
        DOI: 10.1109/18.256500

 \bibitem{Chaichenko_Savchuk_Shidlich_2020}
       \textsc{Chaichenko~S., Savchuk~V., Shidlich~A.}
        Approximation of functions by linear summation methods in the Orlicz-type spaces.
        \emph{Journal of Mathematical Sciences (United States)}, 249 (5), 705-719 (2020).
        DOI: 10.1007/s10958-020-04967-y

 \bibitem{Chaichenko_Shidlich_Shulyk_2022}
        \textsc{Chaichenko~S.\,O., Shidlich~A.\,L., Shulyk~T.\,V.}
        Direct and inverse approximation theorems in the
        Besicovitch-Musielak-Orlicz spaces of almost periodic functions.
        \emph{Ukrainian Mathematical Journal}, 74 (5),  801-819 (2022).
        DOI: 10.1007/s11253-022-02102-5

 \bibitem{Chaichenko_Shidlich_2024}
         \textsc{Chaichenko~S.\,O., Shidlich~A.\,L.}
         Approximation of functions by linear methods in
         weighted Orlicz type spaces with variable exponent.
         \emph{Researches in Mathematics}, 32 (2), 70-87 (2024).
         DOI: 10.15421/242420

 \bibitem{DeVore_1998}
         \textsc{DeVore~R.\,A.}
         Nonlinear approximation.
         \emph{Acta Numerica}, 7, 51-150  (1998).
         DOI: 10.1017/S0962492900002816

 \bibitem{DeVore_Petrova_Wojtaszczyk_2024}
         \textsc{DeVore~R., Petrova~G., Wojtaszczyk~P.}
         A note on best $n$-term approximation for generalized Wiener classes.
         \emph{In: DeVore, R., Kunoth, A. (eds) Multiscale, Nonlinear and Adaptive Approximation II.}
         Springer, Cham, 2024.
         DOI: 10.1007/978-3-031-75802-7\_11


 \bibitem{DeVore_Temlyakov_1995}
         \textsc{DeVore~R.\,A.,  Temlyakov~V.\,N.}
         Nonlinear approximation by  trigonometric sums.
         \emph{Journal of Fourier Analysis and Applications}, 2 (1), 29-48  (1995).
         DOI: 10.1007/s00041-001-4021-8

 \bibitem{Dung_Temlyakov_Ullrich_2018}
         \textsc{Dung,~D., Temlyakov,~V.\,N. and  Ullrich,~T.}
         \emph{Hyperbolic Cross Approximation}.
         Advanced Courses in Mathematics -- CRM Barcelona,
         Birkh\"{a}auser/Springer, Cham, 2018.
         DOI: 10.1007/978-3-319-92240-9

 \bibitem{Fang_Qian_2006}
         \textsc{Fang~G., Qian~L.}
         Approximation characteristics for diagonal operators in different computational settings.
         \emph{Journal of Approximation Theory}, 140 (2), 178-190 (2006).
         DOI: 10.1016/j.jat.2005.12.005.

 \bibitem{Fang_Qian_2007}
         \textsc{Fang~G., Qian~L.}
         Optimal algorithms for diagonal operators on $N$-widths in different computational setting.
         \emph{Analysis in Theory and Applications}, 23 (2), 180-187 (2007).
         DOI: 10.1007/s10496-007-0180-z.

 \bibitem{Gruber_2007}
         \textsc{Gruber~P.\,M.}
         \emph{Convex and Discrete Geometry.}
         Die Grundlehren der mathematischen Wissenschaften, 336,
         Berlin, Heidelberg / Springer, 2007.
         DOI: 10.1007/978-3-540-71133-9

 \bibitem{Jahn_Ullrich_Voigtlaender_2023}
         \textsc{Jahn~T.,  Ullrich~T.,   Voigtlaender~F.}
         Sampling numbers of smoothness classes via $l^1$-minimization.
         \emph{Journal of Complexity}, 79, 101786 (2023).
         DOI: 10.1016/j.jco.2023.101786.

 \bibitem{Kahane_1970}
        \textsc{Kahane~J.-P.},
        \emph{S\'{e}ries de Fourier absolument convergentes}.
        Springer Berlin, Heidelberg, 1970
        DOI: 10.1007/978-3-662-59158-1

 \bibitem{Kolomoitsev_Lomako_Tikhonov_2023}
         \textsc{Kolomoitsev~Y., Lomako~T., Tikhonov~S.}
         Sparse grid approximation in weighted Wiener spaces.
         \emph{Journal of Fourier Analysis and Applications}, 29 (2), 19 (2023).
         DOI: 10.1007/s00041-023-09994-2

 \bibitem{Li_2010}
        \textsc{Li~R.\,S., Liu~Y.\,P.}
        Best $m$-term one-sided trigonometric approximation of
        some function classes defined by a kind of multipliers.
        \emph{Acta Mathematica Sinica, English Series}, 26 (5), 975-984 (2010).
        DOI:  10.1007/s10114-009-6478-3

 \bibitem{Moeller_Stasyuk_Ullrich_2024}
        \textsc{Moeller M., Stasyuk S.,  Ullrich T.}
        High-dimensional sparse trigonometric approximation in the uniform norm
        and consequences for sampling recovery.
        \emph{Arxiv: 2407.15965}, 2024.
        DOI:  10.48550/arXiv.2407.15965.

 \bibitem{Moeller_Stasyuk_Ullrich_2025}
       \textsc{Moeller M., Stasyuk S.,  Ullrich T.}
       Best $m$-term trigonometric approximation in weighted Wiener spaces and applications.
       \emph{Arxiv: 2508.07336v2}, 2025.
       DOI: 10.48550/arXiv.2508.07336.

 \bibitem{Moricz_2006}
        \textsc{M\'{o}ricz~F.}
        Absolutely convergent Fourier series and function classes.
        \emph{Journal of Fourier Analysis and Applications}, 324 (2), 1168-1177  (2006).
        DOI: 10.1016/j.jmaa.2005.12.051

 \bibitem{Moricz_2008}
        \textsc{M\'{o}ricz~F.}
        Absolutely convergent Fourier series and function classes. II.
        \emph{Journal of Fourier Analysis and Applications}, 342 (2), 1246-1249  (2008).
        DOI: 10.1016/j.jmaa.2007.12.055

 \bibitem{Moricz_2008_1}
        \textsc{M\'{o}ricz~F.}
        Absolutely convergent multiple Fourier series and multiplicative Lipschitz
        classes of functions.
        \emph{Acta Mathematica Hungarica}, 121 (1-2), 1-19  (2008).
        DOI: 10.1007/s10474-008-7164-0

 \bibitem{Moricz_2010}
        \textsc{M\'{o}ricz~F., S\'{a}f\'{a}r~Z.}
        Absolutely convergent double Fourier series, enlarged Lipschitz and Zygmund
        classes of functions of two variables.
        \emph{East Journal on Approximations}, 16 (1), 1-24  (2010).

 \bibitem{Nguyens_2022}
         \textsc{Nguyen~V.\,K., Nguyen~V.\,D.}
         Best $n$-term approximation of diagonal operators and application
         to function spaces with mixed smoothness.
        \emph{Analysis Mathematica}, 48 (4) (2022), 1127-1152.
        DOI: 10.1007/s10476-022-0169-z

\bibitem{Nikol'skii_1951}
         \textsc{Nikol'skii~S.\,M.}
         Inequalities for entire functions of finite degree and their application
         in the theory of differentiable functions of several variables.
        \emph{American Mathematical Society Translations: Series 2}, 80, 1-38  (1969),
        Transl. from   \emph{Trudy Matematicheskogo Instituta imeni V. A. Steklova},
            38, 244-278 (1951).




 \bibitem{Pinkus_1985}
         \textsc{Pinkus,~A.}
        \emph{$n$-Widths in  approximation theory.}
        Springer-Verlag, Berlin-Heidelberg-New York-Tokyo, 1985.
        DOI: 10.1007/978-3-642-69894-1


 \bibitem{Romanyuk_2012}
        \textsc{Romanyuk,~A.\,S.}
        \emph{Approximation characteristics of classes of periodical functions of several variables}.
        Proceedings of Institute of Mathematics NAS Ukraine, 40,
        Institute of Mathematics NAS Ukraine, Kyiv, 2012.


 \bibitem{Romanyuk_2015}
        \textsc{Romanyuk~V.\,S.}
        Nonlinear approximation in spaces of multiple sequences.
        \emph{In: Approximation theory of functions and related problems [in Ukrainian],
        Transactions of Institute of Mathematics of NAS of Ukraine}, 12 (4), 253-261 (2015).

 \bibitem{Savchuk_Shidlich_2014}
        \textsc{Savchuk~V.\,V., Shidlich~A.\,L.}
        Approximation of functions of several variables by linear methods in the space $S^p$.
        \emph{Acta Scientiarum Mathematicarum (Szeged),} 80 (3-4), 477-489 (2014).
        DOI: 10.14232/actasm-012-837-8

 \bibitem{Serdyuk_Shidlich_2022}
         \textsc{Serdyuk~A.\,S., Shidlich~A.\,L.}
         Direct and inverse theorems on the approximation
         of almost periodic functions in Besicovitch-Stepanets spaces.
        \emph{Carpathian Mathematical Publications}, 13 (3), 687-700 (2021).
         DOI: 10.15330/cmp.13.3.687-700

 \bibitem{Serdyuk_Shidlich_2025}
         \textsc{Serdyuk A.\,S., Shidlich A.\,L.}
         Actual problems of the theory of approximations in metrics of discrete
         spaces on the sets of summable periodic and almost periodic functions.
         \emph{Ukrainian Mathematical Journal}, 76 (11),  1858-1900  (2025).
         DOI: 10.1007/s11253-025-02428-w


   \bibitem{Shidlich_2016}
        \textsc{Shidlich~A.\,L.}
         Nonlinear approximation of the classes ${\mathcal F}^{\psi}_{q,r}$ of
         functions of several variables in the integral metrics.
         \emph{In: Mathematical problems of mechanics and computational mathematics,
          Transactions of Institute of Mathematics of NAS of Ukraine},
          13 (3), 256-274 (2016).


  \bibitem{Shydlich_UMZh2009}
        \textsc{Shydlich~A.\,L.}
         Order equalities for some functionals and their application to the estimation of the
         best $n$-term approximations and widths.
         \emph{Ukrainian Mathematical Journal}, 61 (10),  1649-1671  (2009).
         DOI: 10.1007/s11253-010-0304-z

 \bibitem{Shidlich_2011}
         \textsc{Shydlich~A.\,L.}
         Order estimates of the best $n$-term orthogonal trigonometric
         approximations of classes ${\mathcal F}^\psi_{q,\infty}$ in spaces $L_p({\mathbb T}^d)$.
         \emph{In: Problems of Approximation Theory of Functions and Related Problems [in Ukrainian],
         Transactions of Institute of Mathematics of NAS of Ukraine},
          8 (1), 224-243 (2011).


   \bibitem{Shidlich_2013}
          \textsc{Shydlich~A.\,L.}
          Order estimates for some approximation characteristics,
         \emph{In: Approximation theory of functions and related problems [in Ukrainian].
          Transactions of Institute of Mathematics of NAS of Ukraine}, 10 (1),
          304-337 (2013).

   \bibitem{Shidlich_2014}
        \textsc{Shydlich~A.\,L.}
        Order estimates of functionals in terms of which the
        best $n$-term approximations of the classes ${\mathcal F}^{\psi}_{q,r}$ are expressed.
        \emph{In: Approximation theory of functions and related problems [in Ukrainian],
        Transactions of Institute of Mathematics of NAS of Ukraine}, 11 (3), 287-314  (2014).


 \bibitem{Stechkin_1951}
          \textsc{Stechkin~S.\,B.}
          On absolute convergence of orthogonal series. I.
          \emph{Matematicheskii Sbornik (Novaya Seriya)}, 29 (71), (1), 225-232 (1951).


 \bibitem{Stechkin_1953}
          \textsc{Stechkin~S.\,B.}
          On absolute convergence of Fourier series.
          \emph{Izvestiya Akademii Nauk SSSR. Seriya Matematicheskaya}, 17 (2), 87-98 (1953).

 \bibitem{Stechkin_1955_IZV}
          \textsc{Stechkin~S.\,B.}
          On absolute convergence of Fourier series. II.
          \emph{Izvestiya Akademii Nauk SSSR. Seriya Matematicheskaya}, 19, 4  (1955), 221-246.


  \bibitem{Stechkin_1955}
          \textsc{Stechkin~S.\,B.}
          On absolute convergence of orthogonal series.
          \emph{Doklady Akademii Nauk SSSR (Novaya Seriya)}, 102, 37-40 (1955).


 \bibitem{Stepanets_1999}
          \textsc{Stepanets~A.\,I.}
          Several statements for convex functions.
          \emph{Ukrainian Mathematical Journal}, 51 (5), 764-780  (1999).
           DOI: 10.1007/BF02591710


  \bibitem{Stepanets_2001}
          \textsc{Stepanets~A.\,I.}
          Approximation characteristics of the spaces  $S^p_\varphi$ in different metrics.
          \emph{Ukrainian Mathematical Journal}, 53 (10), 1340-1374 (2001).
          DOI: 10.1023/A:1013307912783

 \bibitem{Stepanets_2005}
          \textsc{Stepanets~A.\,I.}
          \emph{Methods of approximation theory}.
          VSP, Leiden--Boston, 2005.
          DOI: 10.1515/9783110195286

 \bibitem{Stepanets_UMZh2006_1}
          \textsc{Stepanets~A.\,I.},
          Problems of approximation theory in linear spaces.
          \emph{Ukrainian Mathematical Journal}, 58 (1), 54-102  (2006).
          DOI: 10.1007/s11253-006-0052-2

 \bibitem{Stepanets_Serdyuk_UMZh2002}
        \textsc{Stepanets~A.\,I.,  Serdyuk~A.\,S.}
        Direct and inverse theorems in the theory of the approximation of
        functions in the space ${\mathcal S}^p$.
        \emph{Ukrainian Mathematical Journal}, 54  (1), 322-338 (2002).
        DOI: 10.1023/A:1019701805228

 \bibitem{Stepanets_Shidlich_2010}
        \textsc{Stepanets~A.\,I., Shidlich~A.\,L. }
        Best approximations of integrals by integrals of finite rank.
        \emph{Journal of Approximation Theory}, 162 (2), 323-348  (2010).
        DOI: 10.1016/j.jat.2009.05.007


\bibitem{Stepanets_Shidlich_IZV_2010}
        \textsc{Stepanets~A.I., Shidlich~A.\,L.}
        Extremal problems for integrals of nonnegative functions,
        \emph{Izvestiya Mathematics}, 74 (3), 607-660 (2010).
        DOI: 10.1070/IM2010v074n03ABEH002500

 \bibitem{Stepanets_Shydlich_2003}
        \textsc{Stepanets'~O.\,I., Shydlich~A.\,L.}
        Best $n$-term approximations by $\Lambda$-methods in the spaces $S^\phi_p$.
        \emph{Ukrainian Mathematical Journal}, 55 (8), 1334-1358 (2003).
        DOI: 10.1023/B:UKMA.0000010763.01827.6f

 \bibitem{Szasz_1942}
        \textsc{Sz\'{a}sz~O.}
        On Convergence and summability of trigonometric series.
        \emph{American Journal of Mathematics}, 61 (1), 575-591 (1942).
        DOI: 10.2307/2371705

 \bibitem{Sunouchi_1951}
        \textsc{Sunouchi~G.}
        A convergence criterion for Fourier series.
        \emph{Tohoku Mathematical Journal (Second Series)}, 3 (2), 216-219 (1951).
        DOI: 10.2748/tmj/1178245525


 \bibitem{Sunouchi_1952}
        \textsc{Sunouchi~G.}
        Convergence criteria for Fourier series.
        \emph{Tohoku Mathematical Journal (Second Series)}, 4 (2),  187-193 (1952).
        DOI: 10.2748/tmj/1178245421

 \bibitem{Temlyakov_1998}
        \textsc{Temlyakov~V.\,N.}
        Greedy algorithm and $m$-term trigonometric approximation.
        \emph{Constractive Approximation}, 14 (4), 569-587  (1998).
        DOI: 10.1007/s003659900090

 \bibitem{Temlyakov1993}
        \textsc{Temlyakov~V.\,N.}
        \emph{Approximation of Periodic Functions}.
        Computational Mathematics and Analysis Series,
        Commack, New York, Nova Science Publ., 1993.


 \bibitem{Temlyakov_2011}
        \textsc{Temlyakov~V.\,N.}
        \emph{Greedy approximation}.
        Cambridge Monographs on Applied and Computational Mathematics, 20.
        Cambridge University Press, Cambridge, 2011.
        DOI: 10.1017/CBO9780511762291

 \bibitem{Timan_Shavrova_2007}
         \textsc{Timan~M.\,P., Shavrova~O.\,B.}
         Linear Marcinkiewicz type operators for periodic functions of several
         variables in spaces  $S^p_m$ and $L_p^m$.
        \emph{In: Problems of approximation theory and related problems [in Ukrainian],
        Transactions of Institute of Mathematics of NAS of Ukraine}, 4 (1), 352-375 (2007).

\bibitem{Vakarchuk_2004}
        \textsc{Vakarchuk~S.\,B.}
        Jackson-type inequalities and exact values of widths of classes of
        functions in the spaces  $S^p$, $1\leq p< \infty$.
        \emph{Ukrainian Mathematical Journal}, 56 (5), 718-729  (2004).
        DOI: 10.1007/s11253-005-0070-5


 \bibitem{Vakarchuk_Shchitov_2006}
        \textsc{Vakarchuk~S.\,B., Shchitov~A.\,N.}
        On some extremal problems in the theory of approximation of functions
        in the spaces $S^p$, $1\leq p< \infty $.
        \emph{Ukrainian Mathematical Journal}, 58 (3) (2006), 340-356.
        DOI: 10.1007/s11253-006-0070-0



 \bibitem{Voigtlaender_2022}
       \textsc{Voigtlaender~F.}
        $L_p$-sampling numbers for the Fourier-analytic Barron space.
        \emph{arXiv: 2208.07605}, 2022.
        DOI: 10.48550/arXiv.2208.07605


 \bibitem{Volosivets_2025}
         \textsc{Volosivets~S.}
         Approximation in Orlicz-Stepanets spaces defined by Vilenkin-Fourier coefficients.
        \emph{Journal of Mathematical Sciences (United States)}, 290 (2), 287-300  (2025).
         DOI: 10.1007/s10958-025-07607-5

 \bibitem{Wiener_1933}
        \textsc{Wiener~N.}
        \emph{The Fourier Integral and Certain of Its Applications}.
        Cambridge University Press, Cambridge, 1933
        DOI: 10.1017/CBO9780511662492


 \bibitem{Wang_2005}
        \textsc{Wang~X.}
        Volumes of Generalized Unit Balls.
        \emph{Mathematics Magazine}, 78 (5) (2005), 390-395.
         DOI: 10.2307/30044198




     \end{thebibliography}}
\enddocument